\documentclass[aap,MSNbibl,citesort,dvips]{arximspdf}
\usepackage{graphicx}

\doi{10.1214/10-AAP757}
\volume{21}
\issue{6}
\pubyear{2011}
\firstpage{2263}
\lastpage{2314}

\makeatletter

\newcommand{\cc}{{\mathcal C}}
\newcommand{\ce}{{\mathcal E}}
\newcommand{\ci}{{\mathcal I}}
\newcommand{\cm}{{\mathcal M}}
\newcommand{\cp}{{\mathcal P}}
\newcommand{\ct}{{\mathcal T}}
\newcommand{\cu}{{\mathcal U}}

\newcommand{\rme}{\mathrm{e}}
\newcommand{\Co}{\mathcal{C}}
\newcommand{\dlangle}{\langle\! \langle}
\newcommand{\drangle}{\rangle\! \rangle}

\newcommand{\D}{\mathbb{D}}
\newcommand{\N}{\mathbb{N}}
\newcommand{\PP}{\mathbb{P}}
\newcommand{\R}{\mathbb{R}}
\newcommand{\E}{\mathbb{E}}
\newcommand{\U}{\mathcal{U}}
\newcommand{\T}{\mathbb{T}}
\newcommand{\w}{\widetilde}
\newcommand{\ind}{\mathbf{1}}

\newproclaim{ex}{Example}

\newtheorem{theorem}{Theorem}[section]
\newtheorem{prop}[theorem]{Proposition}
\newtheorem{lemme}[theorem]{Lemma}
\newtheorem{corol}[theorem]{Corollary}

\newproclaim{definition}[theorem]{Definition}

\newtheorem{hyp}[theorem]{Assumption}

\newproclaim{rque}[theorem]{Remark}

\makeatother

\begin{document}
\begin{frontmatter}

\title{Limit theorems for Markov processes indexed by continuous time
Galton--Watson trees\thanksref{T1}}
\runtitle{Markov processes on Galton--Watson trees}

\thankstext{T1}{Supported by the ANR MANEGE (ANR-09-BLAN-0215),
the ANR A3 (ANR-08-BLAN-0190), the ANR Viroscopy (ANR-08-SYSC-016-03)
and the Chair ``Mod\'{e}lisation Math\'{e}matique et Biodiversit\'{e}''
of Veolia Environnement-Ecole Polytechnique-Museum National
d'Histoire Naturelle-Fondation Ec. Polytechnique.}

\begin{aug}
\author[A]{\fnms{Vincent} \snm{Bansaye}\ead[label=e1]{vbansaye@gmail.com}\ead[label=u1,url]{http://www.cmapx.polytechnique.fr/\textasciitilde bansaye/}},
\author[B]{\fnms{Jean-Fran\c{c}ois} \snm{Delmas}\ead[label=e2]{delmas@cermics.enpc.fr}\ead[label=u2,url]{http://cermics.enpc.fr/\textasciitilde delmas/}},
\author[C]{\fnms{Laurence} \snm{Marsalle}\ead[label=e3]{laurence.marsalle@math.univ-lille1.fr}\ead[label=u3,url]{http://math.univ-lille1.fr/\textasciitilde marsalle/}}
and
\author[A,C]{\fnms{Viet Chi} \snm{Tran}\corref{}\ead[label=e4]{chi.tran@math.univ-lille1.fr}\ead[label=u4,url]{http://math.univ-lille1.fr/\textasciitilde tran/}}
\runauthor{Bansaye, Delmas, Marsalle and Tran}
\affiliation{Ecole Polytechnique, Ecole Nationale des Ponts et Chauss\'
{e}es,
Universit\'{e} des Sciences et Technologies Lille 1, and Ecole
Polytechnique and
Universit\'{e} des Sciences et Technologies Lille 1}
\address[A]{V. Bansaye\\
V. C. Tran\\
Centre de Math\'{e}matiques Appliqu\'{e}es (CMAP)\\
UMR CNRS 7641 Ecole Polytechnique\\
Route de Saclay\\
91128 Palaiseau Cedex\\
France\\
\printead{e1}\\
\printead{u1}}
\address[B]{J.-F. Delmas\\
CERMICS, Ecole Nationale\\
\quad des Ponts et Chauss\'{e}es\\
6 et 8, Avenue Blaise Pascal\\
Cit\'{e} Descartes---Champs-sur-Marne\\
77455 Marne-la-Vall\'{e}e Cedex 2\\
France\\
\printead{e2}\\
\printead{u2}}
\address[C]{L. Marsalle\\
V. C. Tran\\
Laboratoire Paul Painlev\'{e}\\
UFR de Math\'{e}matiques\\
Universit\'{e} des Sciences et Technologies Lille 1\\
Cit\'{e} scientifique\\
UMR CNRS 8524\\
59655 Villeneuve d'Ascq Cedex\\
France\\
\printead{e3}\\
\hphantom{E-mail: }\printead*{e4}\\
\printead{u3}\\
\hphantom{URL: }\printead*{u4}}
\end{aug}

% HISTORY:
\received{\smonth{11} \syear{2009}}
\revised{\smonth{9} \syear{2010}}

% ABSTRACT
%
\begin{abstract}
We study the evolution of a particle system whose genealogy is given by
a supercritical continuous time Galton--Watson tree. The particles move
independently according to a Markov process and when a branching event
occurs, the offspring locations depend on the position of the mother
and the number of offspring. We prove a~law of large numbers for the
empirical measure of individuals alive at time $t$. This relies on a
probabilistic interpretation of its intensity by mean of an auxiliary
process. The latter has the same generator as the Markov process along
the branches plus additional jumps, associated with branching events of
accelerated rate and biased distribution. This comes from the fact that
choosing an individual uniformly at time $t$ favors lineages with more
branching events and larger offspring number. The central limit
theorem is considered on a special case. Several examples are
developed, including applications to splitting diffusions, cellular
aging, branching L\'{e}vy processes.
\end{abstract}

% KEYWORDS
%

\end{frontmatter}

%s1 ###
\section{Introduction and main results}\label{intro}

We consider a continuous time Galton--Watson tree $\T$, that
is, a tree where each branch lives during an independent exponential
time of
mean $1/r$,
then splits into a random number of new branches
given by an independent random variable (r.v.) $\nu$ of law $(p_k,
k\in
\N)$, where\vadjust{\goodbreak} $\N=\{0,1,\ldots\}$. We are interested in the following
process indexed by this
tree. Along the edges of the tree, the process evolves as a c\`{a}dl\`{a}g
strong Markov process $(X_t, t\geq 0)$ with values in a Polish space~%
$E$ and with infinitesimal generator~$L$ of domain $D(L)$. The branching
event is nonlocal; the positions of the offspring are described by a
random vector $(F^{(k)}_1(x,\Theta), \ldots, F_{k}^{(k)}(x,\Theta
))$, which depends on the position $x$ of the mother just before the
branching event and on the number $\nu=k$ of offspring; the randomness
of these positions is modeled via the random variable~$\Theta
$, which is uniform on $[0,1]$.
Finally, the newborn branches evolve independently from each other.

This process is a branching Markov process for which there has been
vast literature.
We refer to Asmussen and Hering \cite{asmussenhering}, Dawson
\cite{dawson} and Dawson, Gorostiza and Li \cite{dawsongorostizali} for
nonlocal branching processes similar to those considered here. Whereas
the literature often deals with limit theorems that consider
superprocess limits corresponding to high densities of small and
rapidly branching
particles (see, e.g., Dawson \cite{dawson}, Dynkin \cite{dynkin},
Evans and
Steinsaltz \cite{evanssteinsaltz}), here we stick with the discrete
tree in continuous time.% which we aim at characterizing.

Let us also mention some results in the discrete time case. Markov
chains indexed by a binary tree have been studied in the symmetric
independent case (see, e.g., Athreya and Kang \cite{atarbre}, Benjamini and
Peres \cite{benj}) where for every $x$, $F^{(2)}_1(x,\Theta)$ and
$F^{(2)}_2(x,\Theta)$ are i.i.d. A motivation for considering asymmetric
branching comes from models for cell division. For instance, the binary tree
can be used to describe a dividing cell genealogy in discrete time. The
Markov chain indexed by this binary tree then indicates the evolution of
some characteristic of the cell, such as its growth rate, its quantity
of proteins or parasites$\ldots$ and depends on division events.
Experiments (Stewart et al.~\cite{tad}) indicate that the transmission
of this characteristic in the two daughter cells may be asymmetric. See
Bercu, De Saporta and G\'egout-Petit~\cite{bercugegoutpetitsaporta} or
Guyon \cite{guyon} for
asymmetric models for cellular aging and Bansaye~\cite{bansaye} for
parasite infection. In Delmas and Marsalle \cite{delmasmarsalle} a
generalization of these models where there might be 0, 1 or 2 daughters
is studied. Indeed under stress conditions, cells may divide less or
even die. The branching Markov chain, which in their case represents the
cell's growth rate, is then restarted for each daughter cell at a value
that depends on the mother's growth rate and on the total number of
daughters.

We investigate the continuous time case and allow both asymmetry and
random number of offspring. To illustrate this model, let us give two
simple examples
related to parasite infection problems. In the first case, the cell
divides in two
daughter cells after an exponential time and a random fraction of
parasites goes in one of the daughter cells, whereas the rest go in the
second one. In the second case, the cell divides in $k$ daughter cells
and the process $X$ is equally shared between each of the $k$ daughters:
$\forall j\in\{ 1,\ldots, k\}, F_j^{(k)}(x,\Theta)=x/k$. Notice that
another similar model has been investigated\vadjust{\goodbreak} in Evans and Steinsaltz
\cite{evanssteinsaltz} where the evolution of damages in a population of
dividing cells is studied, but with a superprocess point of view. The
authors assume that the cell's death rate depends on the damage of the
cell, which evolves as a~diffusion between two fissions. When a division
occurs, there is an unbalanced transmission of damaged material that
leads to the consideration of nonlocal births. Further examples are
developed in Section \ref{sectionexamples}.

Our main purpose is to characterize the empirical distribution of this
process. More precisely, if we denote by $N_t$ the size of the living
population~$V_t$ at time $t$, and if $(X_t^u, u\in V_t)$ denotes the
values of the Markov process for the different individuals of $V_t$, we
will focus on the following probability measure which describes the
state of the population:
\[
\frac{\ind_{\{N_t>0\}}}{N_t}\sum_{u\in V_t}\delta
_{X^u_t}(dx),\qquad t\in\R_+.
\]
This is linked to the value of the process for an individual chosen
uniformly at time~$t$, say $U(t)$, as we can see
from this simple identity,
\[
\E\biggl[\frac{\ind_{\{N_t>0\}}}{N_t}\sum_{u\in V_t} f(X^u_t)
\biggr]=\E\bigl[\ind_{\{N_t>0\}}f\bigl(X^{U(t)}_t\bigr)\bigr].
\]
We show that the distribution of the path leading from the ancestor to
a~uniformly chosen individual can be approximated by means of an
\textit{auxiliary Markov process} $Y$ with infinitesimal generator
characterized by $\forall f\in D(L)$,
%
%e1 ###
%
\begin{equation}\label{generateurdeY}
Af(x)=Lf(x)+rm \sum_{k=1}^{+\infty} \frac{p_k}{m} \int_0^1 \sum
_{j=1}^k \bigl( f\bigl(F_j^{(k)}(x,\theta)\bigr) -f(x)\bigr)\,d\theta,
\end{equation}
where we recall that $r$ denotes the particle branching
rate and where we introduce $m=\sum_{k=1}^{+\infty} k p_k$ the mean
offspring number. In this paper, we will be interested in the
supercritical case $m>1$, even if some remarks are made for the
critical and subcritical cases. The auxiliary process has the same
generator $L$ as the Markov process running along the branches, plus jumps
due to the branching. However, we can observe a bias phenomenon: the
resulting jump rate $rm$ is equal to the original rate $r$ times the mean
offspring number $m$ and the resulting offspring distribution is the
size-biased distribution $(k p_k/m,k\in\N)$. For $m>1$, for instance,
this is heuristically explained by the fact that when one chooses an
individual uniformly in the population at time $t$, an individual
belonging to a lineage with more generations or with prolific ancestors
is more likely to be chosen. Such biased phenomena have already been
observed in the field of branching processes (see, e.g., Chauvin,
Rouault and Wakolbinger
\cite{chauvinrouaultwakolbinger}, Hardy and Harris \cite
{hardyharris3}, Harris and Roberts~\cite{harrisroberts}). Here, we allow nonlocal births, prove pathwise
results and establish laws of large numbers when $Y$ is ergodic. Our
approach is entirely based on a probabilistic interpretation
via the auxiliary process\vadjust{\goodbreak} $Y$.

In case $Y$ is ergodic, we prove the laws of large numbers stated in
Theorems \ref{thLGNannonceintro} and \ref{thLGNpathannonceintro},
where $W$ stands for the renormalized asymptotic size of the number of
individuals at time $t$ (e.g., Athreya and Ney \cite{athreyaney},
Theorems 1 and 2, page 111),
\[
W:=\lim_{t\rightarrow+\infty}N_t/\E[N_t] \qquad\mbox{a.s.}
\quad\mbox{and}\quad \{W>0\}=\{\forall t\geq0, N_t>0\} \qquad\mbox{a.s.}
\]

\begin{theorem}\label{thLGNannonceintro}
If the auxiliary process $Y$ is ergodic with invariant measure~$\pi$,
we have
for any real continuous bounded function $f$ on $E$,
%
%e2 ###
%
\begin{equation}
\lim_{t\rightarrow\infty} \frac{\ind_{\{N_t>0\}}}{N_t}\sum_{u\in
V_t}f(X^u_t)= \ind_{\{W>0\}}\int_{E} f(x)\pi(dx)
\qquad\mbox{in probability}.
\end{equation}
\end{theorem}

This result in particular implies that for such function $f$,
%
%e3 ###
%
\begin{equation}\label{equationintro2}
\lim_{t\rightarrow+\infty} \E\bigl[f\bigl(X_t^{U(t)}\bigr) | N_t>0
\bigr]=\int_{E} f(x) \pi(dx),
\end{equation}
where $U(t)$ stands for a particle taken at random in the
set $V_t$ of living particles at time $t$.

Theorem \ref{thLGNannonceintro} is a consequence of Theorem \ref{thLGN}
(which gives similar results under weaker hypotheses) and of Remark
\ref{rquecontientcontinuborne}. The convergence is proved using $L^2$
techniques.

Theorem \ref{thLGNannonceintro} also provides a limit theorem for the
empirical distribution of the tree indexed Markov process.
\begin{corol}\label{corolintro}Under the assumption of Theorem \ref
{thLGNannonceintro},
%
%e4 ###
%
\begin{equation}
\lim_{t\rightarrow\infty} \frac{\ind_{\{N_t>0\}}}{N_t}\sum_{u\in
V_t}\delta_{X^u_t}(dx)= \ind_{\{W>0\}} \pi(dx) \qquad\mbox{in probability,}
\end{equation}
where the space $\mathcal{M}_F(E)$ of finite measures on $E$ is
embedded with the weak convergence topology.
\end{corol}

We also give in Propositions \ref{propeqfluctu} and \ref
{propconvergencefluctuations} a result on the associated
fluctuations. Notice that contrary to the discrete case treated in
Delmas and Marsalle \cite{delmasmarsalle}, the fluctuation process is a Gaussian process
with a finite variational part.

In addition, we generalize the result of Theorem \ref
{thLGNannonceintro} to ancestral paths of particles (Theorem \ref
{thLGNpathannonceintro}).
\begin{theorem}
\label{thLGNpathannonceintro}
Assume that $Y$ is ergodic with invariant measure $\pi$ and that for
any bounded measurable function $f$,
\[
\lim_{t\rightarrow+\infty}\E_x[f(Y_t)]=\int_{E}f(x)\pi(dx).
\]
Then for any real bounded measurable function $\varphi$ on the
Skorohod space $\D([0,T],E)$, we have the following convergence in probability:
\[
\lim_{t\rightarrow\infty} \frac{\ind_{\{N_t>0\}}}{N_{t}}\sum
_{u\in
V_{t}}\varphi(X^u_{s}, t-T\leq s<t)
=
\E_\pi[\varphi(Y_s, s< T) ]\ind_{\{W \neq0\}},
\]
where, for simplicity, $X^u_s$ denotes the value of the tree indexed
Markov process at time $s$ for the ancestor of $u$ living at this time.
\end{theorem}

Biases that are typical to all renewal problems have been known
for a~long time in the literature (see, e.g., Feller \cite{fellerlivre},
Volume 2, Chapter 1). Size biased trees are linked with the
consideration of Palm measures, themselves related to the problem of
building a population around the path of an individual picked uniformly
at random from the population alive at a certain time $t$. In Chauvin,
Rouault and Wakolbinger \cite{chauvinrouaultwakolbinger} and in Hardy
and Harris \cite{hardyharris3}, a~spinal decomposition is obtained for
continuous time branching processes. Their result states that along the
chosen line of descent, which constitutes a~bridge between the initial
condition and the position of the particle chosen at time $t$, the
birth times of the new branches form a homogeneous Poisson point
process of intensity $rm$ while the reproduction law that is seen along
the branches is given by $(kp_k/m,k\in\N)$. Other references for Palm
measures, spinal decomposition and size-biased Galton--Watson can be
found in discrete time in Kallenberg \cite{kallenberg}, Liemant,
Mattes and Wakolbinger \cite{liemant} and for the continuous time we
mention Gorostiza, Roelly and Wakolbinger~\cite
{gorostizaroellywakolbinger}, Geiger and Kauffmann \cite
{geigerkauffmann}, Geiger \cite{geiger} or Olofsson \cite{olofsson}.
Notice also that biases for an individual chosen uniformly in a
continuous time tree had previously been observed by Samuels \cite
{samuels} and Biggins \cite{bigginsaap76}. In the same vein, we refer
to Nerman and Jagers \cite{nermanjagers} for consideration of the
pedigree of an individual chosen randomly at time $t$ and to Lyons,
Pemantle and Peres \cite{lyonspemantleperes}, Geiger \cite{geiger2}
for spinal decomposition of size biased discrete-time Galton--Watson processes.

Other motivating topics for this kind of results come from
branching random walks (see, e.g., Biggins \cite{bigginsbrw}, Rouault
\cite{rouaultbrw}) and homogeneous fragmentation (see Bertoin \cite
{bertoinfragmentation,bertoinrandomfrag}). We refer to the examples in
Section \ref{sectionexamples} for more details.
%For example, consider the case where each fragment divides at a
%constant rate $r\in(0,\infty)$ in two fragments of the same size. We
%then recover the rate of decrease of a fragment chosen uniformly among
%fragments in time $t$, which can be derived from
%decrease of a fragment obtained by keeping randomly with the same
%probability one of the two fragments at each division as well as the
%one of a tagged fragment are equal to $r$.

The law of large numbers that we obtain belongs to the family of
law of large numbers (LLN) for branching processes and superprocesses.
We mention Benjamini and Peres \cite{benj} and Delmas and Marsalle
\cite{delmasmarsalle} in discrete time, with spatial motion for the
second reference. In continuous time, LLNs have been obtained by
Georgii and Baake \cite{georgiibaake} for multitype branching
processes. Finally, in the more different setting of superprocesses
(obtained by renormalization in large population and where individuals
are lost), Engl\"{a}nder and Turaev \cite{englanderturaev}, Engl\"
{a}nder and Winter \cite{englanderwinter} and Evans and Steinsaltz
\cite{evanssteinsaltz} have proved similar results. Here, we work in
continuous time, discrete population, with spatial motion and nonlocal
branching. This framework allows us to trace individuals which may be
interesting for statistical applications. Our results are obtained by
means of the auxiliary process $Y$, while the other approaches involve
spectral techniques and changes of measures via martingales.

In Section \ref{sec2}, we define our Markov process indexed by a continuous time
Galton--Watson tree. We start with the description of the tree and then
provide a measure-valued description of the process of interest. In
Section~\ref{sec:MT1}, we build an \textit{auxiliary process} $Y$ and prove that
its law is deeply related to the distribution of the lineage of an
individual drawn uniformly in the population. In Section \ref{sectionlgn}, we
establish the laws of large numbers mentioned in Theorems \ref
{thLGNannonceintro} and \ref{thLGNpathannonceintro}. Several examples
are then investigated in Section~\ref{sectionexamples}: splitting diffusions indexed by a
Yule tree, a model for cellular aging generalizing Delmas and Marsalle \cite
{delmasmarsalle} and an application to nonlocal branching random walks.
Finally, a central limit theorem is considered for splitting diffusions
in Section \ref{sec6}.

%s2 ###
\section{Tree indexed Markov processes}\label{sec2}
We first give a description of the continuous time Galton--Watson trees
and preliminary estimates in Section~\ref{sectionGWcontinu}.
Section~\ref{sectiondescriptiontreeindexedmarkovprocess} is devoted to the
definition of tree indexed Markov processes.

%s2.1 ###
\subsection{Galton--Watson trees in continuous time}\label{sectionGWcontinu}

In a first step we recall some definitions about discrete trees. In a
second step, we introduce continuous time and finally, in a third step,
we give the definition of the Galton--Watson tree in continuous time.
For all this section, we refer mainly to~\cite{duquesnelegall,harris,lambert}.

\subsubsection*{Discrete trees} Let
%
%e5 ###
%
\begin{equation}
\mathcal{U}=\bigcup_{m=0}^{+\infty} (\N^*)^m,
\end{equation}
where $\N^*=\{1,2, \ldots\}$ with the convention $(\N^*)^0 =
\{\varnothing\}$. For $u\in(\N^*)^m$, we define $|u|=m$ the generation of
$u$. If $u=(u_1, \ldots, u_n)$ and $v=(v_1, \ldots, v_p)$ belong
to~$\mathcal{U}$, we write $uv=(u_1, \ldots, u_n,v_1, \ldots, v_p)$ for
the concatenation of $u$ and $v$. We identify both
$\varnothing u$ and $u\varnothing$ with $u$. We also introduce the following
order relation: $u \preceq v$ if there exists $w\in\mathcal{U} $ such
that $v =uw$; if, furthermore, $w\neq\varnothing$, we write $u \prec
v$. Finally, for $u$ and $v$ in $\mathcal{U}$ we define their most
recent common ancestor (MRCA), denoted by $u \wedge v$, as the element
$w \in\mathcal{U}$ of highest generation such that $w \preceq u$ and $w
\preceq v$.
\begin{definition}
A rooted ordered tree $\mathcal{T}$ is a subset of $\mathcal{U}$ such that:
\begin{longlist}
\item$\varnothing\in\mathcal{T}$,
\item if $v \in\mathcal{T}$ then $u \preceq v$ implies $u\in
\mathcal{T}$,\vadjust{\goodbreak}
\item for every $u \in\mathcal{T}$, there exists a number
$\nu_u\in\N$ such that if $\nu_u=0$ then $v\succ u$ implies
$v\notin\mathcal{T}$, otherwise $uj \in\mathcal{T}$ if and only if
$1 \le j \le\nu_u$.
\end{longlist}
\end{definition}

Notice that a rooted ordered tree $\mathcal{T}$ is
completely defined by the sequence $(\nu_u, u \in\mathcal{U})$ which
gives the number of children for every individual. To obtain a
continuous time tree, we simply add the sequence of lifetimes.

\subsubsection*{Continuous time discrete trees} For a sequence $(l_u,
u \in\mathcal{U})$ of nonnegative real numbers, let us define% For
%$u \in\mathcal{U}$, we define
%
%e6 ###
%
\begin{equation}\label{defalphabeta}
\forall u\in\mathcal{U}\qquad \alpha(u)= \sum_{v \prec u} l_v
\quad\mbox{and}\quad
\beta(u)= \sum_{v \preceq u} l_v = \alpha(u)+ l_u
\end{equation}
with the convention $\alpha(\varnothing)=0$. The variable $l_u$ stands
for the lifetime of individual $u$ while $\alpha(u)$ and $\beta(u)$
are its birth and death times. Let
%
%e7 ###
%
\begin{equation}
\mathbb{U}=\mathcal{U} \times[0, + \infty).
\end{equation}

\begin{definition} A continuous time rooted discrete tree (CT)
is a subset~$\T$ of $\mathbb{U}$ such that:
\begin{longlist}
\item $(\varnothing, 0) \in\T$,
\item the projection of $\T$ on $\mathcal{U}$, $\mathcal
{T}$, is a discrete rooted ordered tree,
\item there exists a sequence of nonnegative real numbers
$(l_u,u\in
\mathcal{U})$ such that for every $u \in\mathcal{T}$, $(u,s) \in\T
$ if and only if $\alpha(u) \le s < \beta(u)$, where $\alpha(u)$ and
$\beta(u)$ are defined by (\ref{defalphabeta}).
\end{longlist}
\end{definition}

Let $\T$ be a CT. The set of individuals of $\T$ living at
time $t$ is denoted by~$V_t$,
%
%e8 ###
%
\begin{equation} \label{defVt}
V_t=\{u\in\mathcal{U}\dvtx(u,t)\in\T\}=\{u \in\mathcal{T}\dvtx\alpha
(u)\leq t < \beta(u)\}.
\end{equation}
The number of individuals alive at time $t$ is $N_t=\mbox{Card}(V_t)$.
We denote by~$D_t$ the number of individuals which have died before
time $t$,
%
%e9 ###
%
\begin{equation}\label{defDt}
D_t = \mbox{Card} \{u \in\mathcal{T}\dvtx\beta(u) <
t\}.
\end{equation}
For $(u, s) \in\T$ and $t \le s$, we call $u(t)$ the ancestor of $u$
living at time $t$,
%
%e10 ###
%
\begin{equation}\label{u(s)}
u(t) = v \qquad\mbox{if } \bigl(v \preceq u \mbox{ and }
(v,t) \in\T\bigr).
\end{equation}
Eventually, for $(u,s) \in\T$, we define the shift of $\T$ at
$(u,s)$ by
$\theta_{(u,s)}\T=\{(v,t) \in\mathbb{U}\dvtx(uv,s+t) \in\T\}$. Note that
$\theta_{(u,s)}\T$ is still a CT.

\subsubsection*{Continuous time Galton--Watson trees}
Henceforth, we work on some probability space denoted by $(\Omega,
\mathcal{F}, \PP)$.
\begin{definition}
We say that a random CT on $(\Omega, \mathcal{F}, \PP)$
is a continuous time Galton--Watson tree with offspring distribution
$p=(p_k,k\in\N)$\vadjust{\goodbreak} and exponential lifetime with mean $1/r$ if:
\begin{longlist}
\item The sequence of the numbers of offspring, $(\nu_u, u \in
\mathcal{U})$, is a sequence of independent random variables with
common distribution $p$.
\item The sequence of lifetimes $(l_u, u \in\mathcal{U})$ is
a sequence of independent exponential random variables with mean $1/r$.
\item The sequences $(\nu_u, u \in\mathcal{U})$ and $(l_u,
u \in\mathcal{U})$ are independent.
\end{longlist}
\end{definition}

We suppose that the offspring distribution $p$ has finite
second moment. We introduce
%
%e11 ###
%
\begin{equation}
m=\sum_{k \ge0} k p_k \quad\mbox{and}\quad \varsigma^2=\sum
_{k\ge0}(k-m)^2 p_k,
\end{equation}
its expectation and variance. The offspring distribution is critical
(resp., supercritical, resp., subcritical) if $m=1$ (resp., $m>1$,
resp., $m<1$). In this work, we mainly deal with the supercritical case.

%f1 ###
%
\begin{figure}

\includegraphics{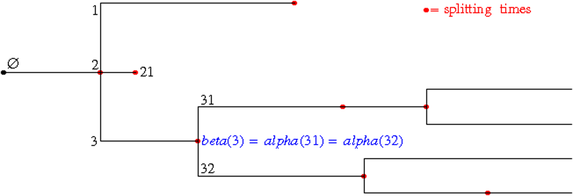}

\caption{Continuous time Galton--Watson tree.}
\label{fig_galton_watson_continu}
\end{figure}

We end Section \ref{sectionGWcontinu} with some estimates on $N_t$ and
$D_t$. To begin with, the following lemma gives an equivalent for $N_t$.
\begin{lemme}\label{lemmeNt}For $t\in\R_+$, we have
%
%e13 ###
%e12 ###
%
\begin{eqnarray}\label{equivdeterministe}
\mathbb{E}[N_t] & = & \rme^{r(m-1)t},\\
\label{solutioncarre}
\E[N_t^2] & = &\cases{
\displaystyle \rme^{r(m-1)t}+\biggl(\frac{\varsigma^2}{m-1}+m\biggr)\bigl(\rme
^{2r(m-1)t}-\rme^{r(m-1)t}\bigr), &\quad if $m\not=1$,\cr
\displaystyle 1+\varsigma^2rt, &\quad if $m=1$.}\hspace*{-28pt}
\end{eqnarray}
If $m>1$, there exists a nonnegative random variable $W$ such that $\{
W>0\}=\{\forall t>0, N_t>0\} $ a.s., $\PP(W > 0) >0$ and
%
%e14 ###
%
\begin{equation}
\label{eq:CV-N-W}
\lim_{t\rightarrow+\infty}\frac{N_t}{\E[N_t]} =
W\qquad \mbox{a.s. and in } L^2.
\end{equation}
\end{lemme}
\begin{pf}
The process $(N_t, t \ge0)$ is a continuous time Markov branching
process so that the expectation and the variance of $N_t$ are
well known (see~\cite{athreyaney}, Chapter~III, Section 4). Almost sure
convergence toward $W$ is stated again in~\cite{athreyaney} (Theorems
1\vadjust{\goodbreak}
and 2, Chapter III, Section 7). Finally, since the martingale $(N_t
\rme^{-r(m-1)t},t \ge0)$ is bounded in $L^2$, we obtain the $L^2$
convergence (e.g., Theorem 1.42, page 11 of \cite{jacodshiryaev}).
\end{pf}

We give the asymptotic behavior of $D_t$, the number of
deaths before $t$.
\begin{lemme}\label{lemmeDt}If $m>1$, the following convergence holds
a.s. and in $L^2$:
%
%e15 ###
%
\begin{equation}\label{convergenceDt}
\lim_{t \rightarrow+\infty}\frac{D_t}{\E[D_t]} = W
\end{equation}
with
%
%e16 ###
%
\begin{equation}
\label{eq:EDt}
\E[D_t]= (m-1)^{-1}\bigl(\rme^{r(m-1)t}-1\bigr)
\end{equation}
and $W$ defined by (\ref{eq:CV-N-W}).
\end{lemme}
\begin{pf}
First remark that
$(D_t,t\geq0)$ is a counting process with compensator $(\int_0^t rN_s
\,ds, t\geq0)$.
We set $\Delta N_t=N_t-N_{t-}$ so that $dN_t= \Delta N_t \,dD_t$. To
prove~(\ref{convergenceDt}), it is sufficient to prove that
$\rme^{-r(m-1)t} I_t$ goes to $0$ a.s. and in~$L^2$, where $I_t=(m-1)D_t
- N_t$. Since $I=(I_t,t\geq0)$ satisfies the following stochastic
equation driven by $(D_t,t\geq0)$:
%
%e17 ###
%
\begin{equation}
dI_t = (m-1-\Delta N_t) \,dD_t,
\end{equation}
we get that $I$ is an $L^2$ martingale. We deduce that $d\langle I
\rangle_t= \varsigma^2 r N_t \,dt$ and% that is $\langle I\rangle_t=
%
%e18 ###
%
\begin{eqnarray}
\E[ I^2_t]&=&1+\E[ \langle I\rangle_t]=1+\varsigma
^2 r \int_0^t \rme^{r(m-1)s}\,ds\nonumber\\[-8pt]\\[-8pt]
&=&1+\frac{\varsigma^2}{m-1}\bigl(\rme
^{r(m-1)t}-1\bigr),\nonumber
\end{eqnarray}
which implies the $L^2$ convergence of $\rme^{-r(m-1)t} I_t$ to $0$.
Besides, the process $(\rme^{-r(m-1)t} I_t, t\ge0)$ is a supermartingale
bounded in $L^2$ and hence, the convergence also holds almost surely.
% We deduce easily from \reff{equivdeterministe} and \reff{eq:CV-N-W}
%that on
%$\{W>0\}$, $\langle I \rangle_t\sim\frac{\varsigma^2}{m-1}e^{r(m-1)t}
%W$, and from martingale theory (e.g., \cite{jacodshiryaev}
%???????) that
%a.s. $e^{-r(m-1)t} I_t$ converges to $0$. On $\{W=0\}$, $I$ is bounded
%and a.s. $e^{-r(m-1)t} I_t$ converges to $0$. This implies that
%a.s. $e^{-r(m-1)t} I_t$ converges to $0$.
\end{pf}
\begin{ex}[(Yule tree)] \label{ex:YuleNt}
The so-called Yule tree is a continuous time Galton--Watson tree with
a deterministic offspring distribution: each individual of the
population gives birth to 2 individuals, that is, $p_2 =1$ (i.e.,
$p=\delta_2$, the Dirac mass at~2). The Yule tree is thus a binary
tree whose edges have independent exponential lengths with mean $1/r$.
In that case, $W$ is exponential with mean 1 (see, e.g.,
\cite{athreyaney}, page 112). We deduce from Lemma~\ref{lemmeNt}
that, for $t\in\R_+$,
%
%e19 ###
%
\begin{equation}\label{exempleYuletree}
\mathbb{E}[N_t]= \rme^{rt} \quad\mbox{and}\quad \E[N_t^2]= 2\rme
^{2rt}-\rme^{rt}.
\end{equation}
Notice that (\ref{exempleYuletree}) is also a consequence of the
well-known fact that $N_t$ is geometric with parameter $\rme^{-rt}$
(see, e.g., \cite{harris}, page 105).
\end{ex}

%s2.2 ###
\subsection{Markov process indexed by the continuous time Galton--Watson
tree}\label{sectiondescriptiontreeindexedmarkovprocess}

In this section, we define the Markov process $X_\T=(X^u_t, (u,t)\in
\mathbb{T})$\vadjust{\goodbreak} indexed by the continuous time Galton--Watson tree
$\mathbb{T}$ and with initial condition $\mu$. Branching Markov
processes have already been the object of abundant literature (e.g.,
\cite{asmussenhering,asmussenheringbook,athreyaney,ethierkurtz,dawson}).
The process that we consider jumps at branching
times (nonlocal branching property) but these jumps may be dependent.

Let $(E,\ce)$ be a Polish space. We denote by $\mathcal{P}(E)$ the set
of probability measures on $(E,\ce)$.\vspace*{-3pt}

\begin{definition}
\label{def:XT}
Let $X=(X_t,t\geq0)$ be a c\`{a}dl\`{a}g $E$-valued strong Markov process.
Let $\widetilde F=(F^{(k)}_j , 1\leq j \leq k,k\in\N^*)$ be a\vspace*{1pt} family of
measurable functions from $E\times[0,1]$ to $E$. The continuous time branching
Markov (CBM) process $X_\T=(X_t^u,(u,t)\in\T)$ with offspring distribution
$p$, exponential lifetimes with mean $1/r$, offspring position $\widetilde
F$, underlying motion $X$ and starting
distribution $\mu\in\cp(E)$, is defined recursively as follows:
\begin{longlist}
\item$\T$ is a continuous time Galton--Watson tree with offspring
distribution $p$ and exponential lifetimes with mean $1/r$.
\item Conditionally on $\T$, $X^\varnothing=(X^\varnothing_t,t\in
[0,\beta(\varnothing)))$ is distributed as $(X_t,t\in
[0,\beta(\varnothing)))$ with $X_0$ distributed as $\mu$.
\item Conditionally on $\T$ and $X^\varnothing$, the initial
positions of
the first genera\-tion offspring $(X^u_{\alpha(u)},1\leq u\leq
\nu_\varnothing)$ are given by
\mbox{$(F^{(\nu_\varnothing)}_u(X^\varnothing_{\beta(\varnothing)-},\Theta),
1\leq
u \leq\nu_\varnothing)$} where $\Theta$ is a uniform random variable on
$[0,1]$.% independent of $X^\varnothing$ and $\mathbb{T}$.
\item Conditionally\vspace*{2pt} on $X^\varnothing$, $\nu_\varnothing$,
$\beta_\varnothing$ and $(X^u_{\alpha(u)},1\leq u\leq\nu_\varnothing)$,
the tree-indexed Markov processes $(X^{uv}_{\alpha(u)+t}, (v,t)\in
\theta_{(u,\alpha(u))}\mathbb{T})$ for $1\leq u \leq\nu_\varnothing$
are independent and, respectively, distributed as $X_\T$ with
starting distribution the Dirac mass at $X^u_{\alpha(u)}$.\vspace*{-3pt}
\end{longlist}
\end{definition}

For $x \in E$, we define $\PP_x(A)=\PP(A |X_0^{\varnothing}=x)$ for all
$A \in\mathcal{F}$ and denote by~$\E_x$ the corresponding
expectation. For $\mu\in\mathcal{P}(E)$ we set in a classical manner
$\PP_\mu(A) = \int_E \PP_x(A) \mu(dx)$ and write $\E_\mu$ for the
expectation w.r.t. $\PP_\mu$.

%f2 ###
%
\begin{figure}

\includegraphics{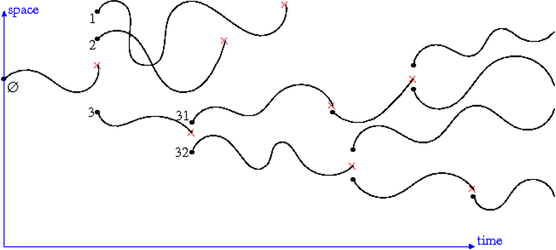}

\caption{Continuous time Markov process indexed by the Galton--Watson
tree of Figure~\protect\ref{fig_galton_watson_continu}.}
\label{fig_processus_branchement_non_local2}
\vspace*{-3pt}
\end{figure}

For $u \in\mathcal{T}$, we extend the definition of $X_t^u$ when $t
\in[0, \alpha(u))$ as follows: $X_t^u=X_t^{u(t)}$, where $u(t)$,
defined by
(\ref{u(s)}), is the ancestor of $u$ living at\vadjust{\goodbreak} time~$t$.

Notice that for $u \in\ct$, $(X_t^u, t \in[0, \beta(u)))$ does not
encode the information about the genealogy of $u$. We remedy this by
introducing the following process $(\Lambda^u_t,t\geq0)$ for $u \in
\mathcal{U}$:
\[
\Lambda^u_t=\sum_{v\prec u(t)} \log(\nu_{v}).
\]
This process provides the birth times of the ancestors of $u$ as well
as their offspring numbers. Notice that it is well defined for all
$u\in\mathcal{U}$ contrary to $X_t^u$. Indeed, the state of $u$ at
its birth time, $X_{\alpha(u)}^u$, is well\vspace*{1pt} defined only for $u \in\ct
$, since it depends on the state of the parent and the number of
its\vadjust{\goodbreak}
offspring.

For $u \in\mathcal{U}$, the process $(\Lambda^u_t, t\in[0,
\beta(u)))$ is a compound Poisson process with rate $r$ for the
underlying Poisson process $(S^u_t, t\geq0)$ and increments
distributed as $\log(\nu)$
with $\nu$ distributed as $p$, stopped at its $(|u|+1)$th
jump.\vspace*{1pt}

In the sequel, we denote by $\w{X}^u_t=(X^u_t, \Lambda^u_t)$ the
couple containing the information on the position and genealogy of the
particle $u$.

%s2.3 ###
\vspace*{-3pt}\subsection{Measure-valued description}\label{sec2.3}
Let $\mathcal{B}_b(E,\R)$ be the set of
real-valued measurable bounded functions on $E$ and $\cm_F(E)$ the set
of finite measures on $E$ embedded with the topology of weak
convergence. For $\mu\in\cm_F(E)$ and $f\in\mathcal{B}_b(E,\R)$ we
write $\langle\mu,f\rangle=\int_E f(x)\mu(dx)$.

We introduce the following measures to represent the population at $t$:
%
%e20 ###
%
\begin{equation}\label{defz}
\bar{Z}_t=\sum_{u\in V_t}\delta_{(u,X_t^u)}
\quad\mbox{and}\quad {Z}_t=\sum_{u\in V_t}\delta_{X_t^u},
\end{equation}
where $V_t$ has been defined in (\ref{defVt}). Note that $\langle Z_t, f
\rangle= \sum_{u\in V_t} f(X^u_t)$. Since $X$ is c\`{a}dl\`{a}g, we
get that the
process $Z=(Z_t, t\geq0)$ is a c\`{a}dl\`{a}g measure-valued Markov
process of $\mathbb{D}(\R_+,\mathcal{M}_F(E))$.

Following Fournier and M\'{e}l\'{e}ard
\cite{fourniermeleard}, we can describe the evolution of~$Z$ in terms of stochastic differential equations (SDE). Let
$\rho(ds,du,dk, d\theta)$ be a
Poisson point measure of intensity
$r \,ds\otimes n(du)\otimes p( dk)\otimes d\theta$ where $ds$ and
$d\theta$ are Lebesgue measures on $\R_+$ and $[0,1]$, respectively,
$n(du)$ is the counting measure on $\mathcal{U}$ and $p(dk)=\sum
_{i\in
\N} p_i \delta_i (dk)$ is the
offspring distribution. This measure
$\rho$ gives the information on the branching events. Let $L$ be the
infinitesimal generator of $X$. If $\Co^{1,0}_b(\R_+\times E,\R)$ denotes
the space of continuous bounded functions that are $\Co^1$ in time with
bounded derivatives, then for test functions
$f\dvtx(t,x)\mapsto f_t(x)$ in $\Co^{1,0}_b(\R_+\times E,\R)$ such that
$\forall t\in\R_+,f_t\in D(L)$, we have
%
%e21 ###
%
\begin{eqnarray}\label{martingalegdepop}
\langle Z_t,f_t\rangle &=& f_0(X^\varnothing_0)+\int_0^t\int_{\R
_+}\bigl(Lf(x) +\partial_s
f(x)\bigr)Z_s(dx) \,ds + W^f_t\nonumber\\[-2pt]
&&{} + \int_0^t
\int_{\mathcal{U}\times\N\times[0,1]} \ind_{\{u\in V_{s-}\}}
\Biggl(\sum_{j=1}^k
f_s\bigl(F^{(k)}_j(X^u_{s_-},\theta)\bigr)-f_s(X^u_{s_-})\Biggr)\\[-2pt]
&&\hphantom{{} +\int_0^t
\int_{\mathcal{U}\times\N\times[0,1]}}
{}\times\rho(ds,du,
dk,d\theta),\nonumber\vadjust{\goodbreak}
\end{eqnarray}
where $W^f_t$ is a martingale. Explicit expressions of this martingale
and of the infinitesimal generator of $(Z_t,t\ge0)$ can be obtained
when the form of the generator $L$ is given.
\begin{ex}[(Splitted diffusions)]\label{exfragmenteddiffusion}
The case when the
Markov process $X$ is a real diffusion ($E=\R$) is an interesting
example. Let $L$ be given by
%
%e22 ###
%
\begin{equation}\label{generateurexemple1}
Lf(x)=b(x)f'(x)+\frac{\sigma^2(x)}{2}f''(x),
\end{equation}
where we assume that $b$ and $\sigma$ are bounded and Lipschitz
continuous. In this case, we can consider the following class of
cylindrical functions from~$\mathcal{M}_F(\R)$ into $\R$
defined by $\phi_f(Z)=\phi(\langle Z,f\rangle)$ for $f\in\Co
^2_b(\R,\R)$ and
$\phi\in\Co^2_b(\R)$ which is known to be convergence determining on
$\mathcal{P}(\mathcal{M}_F(\R))$ (see, e.g., \cite{dawson}, Theorem
3.2.6).\vadjust{\goodbreak} We can define the infinitesimal generator $\mathcal{L}$ of
$(Z_t)_{t\geq
0}$ for these functions as
%
%e23 ###
%
\begin{equation}
\mathcal{L}\phi_f(Z)=\mathcal{L}_1 \phi_f(Z)+\mathcal{L}_2\phi_f(Z),
\end{equation}
where $\mathcal{L}_1$ and $\mathcal{L}_2$ correspond to the branching
and motion parts. Such decompositions were already used in Dawson \cite{dawson}
(Section 2.10) and in Roelly and Rouault \cite{roellyrouault}, for
instance. The generator
$\mathcal{L}_1$ is defined by
\begin{eqnarray*}
\mathcal{L}_1\phi_f(Z)
&=& r \int_E \int_0^1 \sum_{k\in
\N}\Biggl(\phi\Biggl(\langle
Z,f\rangle
+\sum_{j=1}^{k}f\bigl(F^{(k)}_j(x,\theta)\bigr)-f(x)\Biggr)-\phi
_f(Z)\Biggr)\\
&&\hphantom{r \int_E \int_0^1 \sum_{k\in\N}}
{}\times p_k
\,d\theta \,Z(dx)
\end{eqnarray*}
with the convention that the sum over $j$ is zero when $k=0$. The
generator~$\mathcal{L}_2$ is given by
\[
\mathcal{L}_2\phi_f(Z)= \langle Z, Lf \rangle\phi'(\langle
Z,f\rangle) +\langle Z, \sigma(x)f'^2(x)\rangle\phi''(\langle
Z,f\rangle).
\]

For a test function $f\dvtx(t,x)\mapsto f_t(x)$ in $\Co^{1,2}_b(\R
_+\times\R,\R)$, the evolution of $(Z_t$, $t\ge0)$ can then be
described by the following SDE:
%
%e24 ###
%
\begin{eqnarray}\label{pbmsplitteddiff}
\langle Z_t,f_t\rangle &=& f_0(X^\varnothing_0)
+ \int_0^t
\int_{\mathcal{U}\times\N\times[0,1]} \ind_{\{u\in V_{s-}\}}
\Biggl(\sum_{j=1}^k
f_s\bigl(F^{(k)}_j(X^u_{s_-},\theta)\bigr)-f_s(X^u_{s_-})\Biggr)\nonumber\\
&&\hphantom{f_0(X^\varnothing_0)
+ \int_0^t
\int_{\mathcal{U}\times\N\times[0,1]}}
{}\times\rho(ds,du,
dk,d\theta)\nonumber\\[-8pt]\\[-8pt]
&&{} + \int_0^t\int_{\R} \bigl(Lf_s(x)+\partial_s f_s(x)
\bigr)Z_s(dx) \,ds\nonumber\\
&&{} +\int_0^t \sum_{u\in V_{s}}\sqrt{2} \sigma(X^u_s) \partial_x
f_s(X^u_s)\,dB_s^u,\nonumber
\end{eqnarray}
where $(B^u)_{u\in\mathcal{U}}$ is a family of independent standard
Brownian motions. In~\cite{bansayetran}, such splitted diffusions are
considered to describe a multi-level population. The cells, which
correspond to the individuals in the present setting, undergo binary
divisions and contain a continuum of parasites that evolves as a Feller
diffusion with drift $b(x)=(b-d)x$ and diffusion $\sigma^2(x)=2\sigma
^2 x$. At the branching time $s$ for the individual $u$, each daughter
inherits a random fraction of the value of the mother. The daughters
$u1$ and $u2$ start, respectively, at $F^{(2)}_1(X_{s_-}^u,\theta
)=G^{-1}(\theta)X_{s_-}^u$ and $F^{(2)}_2(X_{s_-}^u,\theta
)=(1-G^{-1}(\theta))X_{s_-}^u$, where~$G^{-1}$ is the generalized
inverse of $G$, the cumulative distribution function of the random
fraction. %This model describes the limiting evolution of a large
%population of parasites in dividing cells with accelerated births and
%deaths for the parasites.
\end{ex}

%s3 ###
\section{The auxiliary Markov process and Many-to-One formulas}
\label{sec:MT1}

In this section, we are interested in the distribution of the path of
an individual picked at random in the population at time $t$.
%There is a \textit{size bias} phenomenon. Indeed by c
By choosing
uniformly among the individuals present at time $t$, we give a more
important weight to branches where there have been more divisions and
more children since the proportion of the corresponding offspring will
be higher. Our pathwise approach generalizes \cite{delmasmarsalle}
(discrete time) and \cite{bansayetran} (continuous time Yule
process). As mentioned in
the \hyperref[intro]{Introduction}, this size bias has already been observed by
\cite{samuels,bigginsaap76} for the tree structure when
considering marginal distributions and by
\cite{chauvinrouaultwakolbinger,hardyharris3} for local branching Markov
process.

In Section \ref{defprocauxiliaire}, we introduce an auxiliary Markov
process which approximates the distribution of an individual picked at
random among all the individuals living at time $t$.
%of all possible trees. This appears in (\ref{defYt}) for an individual
%living at time $t$.
The relation between $X_{\mathbb{T}}$ and the auxiliary process also
appears when summing the contributions of all individuals of $\mathbb
{T}$ (Section
\ref{sectionmainpties}) and of all pairs of individuals
(Section \ref{identititesforforks}).

%s3.1 ###
\subsection{Auxiliary process and Many-to-One formula at fixed
time}\label{defprocauxiliaire}

We focus on the law of an individual picked at random and show that it
is obtained from an auxiliary Markov process.
%in the set of living people at fixed time.
This auxiliary Markov process $\w{Y}=(Y,\Lambda)$ has two components.
The component $Y$ describes the motion on the space~$E$. The second
component $\Lambda$ encodes a virtual genealogy and $Y$ can then be
seen as the motion along a random lineage of this genealogy. More
precisely,~$\Lambda$ is a compound Poisson process with rate $rm$; its
jump times provide the branching events of the chosen lineage and its
jump sizes are related to the offspring number $H$ whose distribution
is the size biased distribution of $p$.
As for the motion, $Y$ behaves like $X$ between two jumps of $\Lambda
$. At these jump times, $Y$ starts from a new position given by
$F_{J}^{(H)}(\cdot,\Theta)$ where $J$ is uniform on $\{1,\ldots, H\}$ and
$\Theta$ is uniform on $[0,1]$.

For the definition of $\Lambda$, we shall consider the logarithm of
the offspring number as this is the quantity that is involved in the
Girsanov formulas. Notice that we cannot recover all the jump times
from $\Lambda$ unless there is no offspring number equal to 1, that
is, $p_1=0$. This can, however, always be achieved by changing the
value of the jump rate $r$ and adding the jumps related to $F_1^{(1)}$
to the process $X$. Henceforth, we assume without loss of generality
the following.
\begin{hyp}
The offspring distribution satisfies $p_1=0$.
\end{hyp}

By convention for a function $f$ defined on an interval $I$, we set
$f_{\zeta}=(f(t), t\in\zeta)$ for any $\zeta\subset I$.
\begin{definition}
\label{def:Y}
Let $X_\T$ be as in Definition \ref{def:XT} with starting distribution
$\mu\in\cp(E)$. The corresponding auxiliary process
$\w{Y}=(Y,\Lambda)$, with $Y=(Y_t,t\geq0)$ and $\Lambda=(\Lambda_t,
t\geq
0)$, is an $E\times\R$-valued c\`{a}dl\`{a}g Markov process. The process
$(Y,\Lambda)$ and $\ci=(I_k, k\in\N^*)$, a sequence of random
variables, are defined as follows:
\begin{longlist}
\item$\Lambda$ is a compound Poisson process:
$\Lambda_t=\sum_{k=1}^{S_t} \log(H_k)$, where $S=(S_t, t \ge0)$ is a
Poisson process with intensity $rm$ and $(H_k, k\in\N^*)$ are
independent random variables independent of $S$ and with common
distribution the size biased distribution of $p$ $(hp_h/m, h\in\N^*)$.
\item Conditionally on $\Lambda$, $(I_k, k\in\N^*)$ are
independent random variables and $I_k$ is uniform on $\{1, \ldots,
H_k\}$.
\item Conditionally on $(\Lambda,\ci)$, $\tau_1=\inf\{t\ge
0 ; S_t\neq S_0\}$ is known and
the process $Y_{[0, \tau_1)}$ is distributed as $X_{[0,\tau_1)}$.
\item Conditionally on $(\Lambda,\ci,Y_{[0, \tau_1)})$,
$Y_{\tau_1}$ is distributed as $F^{(H_1)}_{I_1}(Y_{\tau_1-},\Theta
)$, where
$\Theta$ is an independent uniform random variable on $[0,1]$.
\item The distribution of $(Y_{\tau_1+t}, t\geq0)$
conditionally on $(\Lambda,\ci,Y_{[0, \tau_1]})$ is equal to the
distribution of $Y$ conditionally on $(\Lambda_{\tau_1+t}-\Lambda
_{\tau_1}, t\geq0)$ and $(I_{1+k},
k\in\N^*)$ and started at $Y_{\tau_1}$.
\end{longlist}
\end{definition}

We write $\E_\mu$ when we take the expectation with respect to
$(Y,\Lambda, \ci)$ and the starting measure is $\mu$ for the $Y$
component. We also use the same
convention as those described just after Definition \ref{def:XT}.

The formula (\ref{defYt}) in the next proposition is similar to the so-called
Many-to-One theorem of Hardy and Harris \cite{hardyharris3} (Section
8.2) that enables expectation of sums over particles in the branching
process to be calculated in terms of an expectation of an auxiliary
process. Notice that in our setting an individual may have no
offspring with positive probability (if $p_0>0$) which is not the case
in \cite{hardyharris3}.
\begin{prop}[(Many-to-One formula at fixed time)]
\label{propprocessusauxiliaire}
For $t\ge0$ and for any nonnegative measurable function $f\in
\mathcal{B}(\D([0,t],E\times\R),\R_+)$,
%
%e25 ###
%
\begin{equation}
\label{defYt}
\frac{\E_\mu[\sum_{u \in V_t} f(\w{X}^u_{[0,t]})
]}{\E[N_t]}
= \mathbb{E}_\mu\bigl[f\bigl(\w{Y}_{[0,t]}\bigr)\bigr].
\end{equation}
\end{prop}
\begin{rque}
\begin{longlist}[(1)]
\item[(1)] Asymptotically, $N_t$ and $\E[N_t]$ are of same order on $\{W>0\}
$ [see (\ref{eq:CV-N-W})]. Thus, the left-hand side of (\ref{defYt})
can be seen as an approximation, for large~$t$, of the law of an
individual picked at random in $V_t$.
\item[(2)] For $m>1$, a typical individual
living at time $t$ has prolific ancestors with shorter lives.
For $m<1$, a typical individual living at time
$t$ has still prolific ancestors but with longer lives.
\item[(3)] If births are local [i.e., for all $ j\leq k$,
$F^{(k)}_j(x,\theta)=x$],
then $Y$ is distributed as $X$.
\end{longlist}
\end{rque}
\begin{pf*}{Proof of Proposition \ref{propprocessusauxiliaire}}
Let $\Lambda$ be a compound Poisson process as in Definition \ref
{def:Y}(i). Let us show the following
Girsanov formula, for any nonnegative measurable function $g$:
%
%e26 ###
%
\begin{equation}\label{utgirsanov}
\E\bigl[g\bigl(\Lambda_{[0,t]}\bigr)\bigr]=\E\bigl[g\bigl(\Lambda'_{[0,t]}\bigr)\rme^{-r(m-1)t
+ \Lambda'_t} \bigr],\vadjust{\goodbreak}
\end{equation}
where the process $\Lambda'$ is a compound process
with rate $r$ for the underlying Poisson process and increments
distributed as $\log(\nu)$ with $\nu$ distributed as $p$.
Indeed, $g(\Lambda_{[0,t]})$ is a function of $t$, of the times $\tau
_q=\inf\{t\ge0; S_t=q\}-\inf\{t\ge0; S_t=q-1\}$ and of jump sizes
$\log(H_q)$ of $\Lambda$,
\[
g\bigl(\Lambda_{[0,t]}\bigr)= \sum_{q=0}^{+\infty}G_q(t,\tau_1,\ldots, \tau
_q,H_1,\ldots, H_q)\ind_{\{\sum_{i=1}^q\tau_i\leq t<\sum
_{i=1}^{q+1}\tau_i\}}
\]
for some functions $(G_q, q\in\N)$. We deduce that
\begin{eqnarray*}
&&\E\bigl[g\bigl(\Lambda_{[0,t]}\bigr)\bigr] \\
&&\qquad= \sum_{q=0}^{+\infty} \int_{\R_+^q}\Biggl( \sum_{h_1,\ldots,
h_q} (rm)^q \rme^{-rmt} G_q(t,t_1,\ldots, t_q,h_1,\ldots, h_q) \\
&&\hspace*{146pt}\qquad\quad{} \times\prod
_{i=1}^q \frac{p_{h_i}h_i}{m} \ind_{\{\sum_{i=1}^q t_i\leq t\}
}\Biggr)\,dt_1\cdots dt_q\\
&&\qquad= \sum_{q=0}^{+\infty} \int_{\R_+^q}\Biggl(\sum_{h_1,\ldots,
h_q} r^q
\rme^{-rt} G_q(t,t_1,\ldots, t_q,h_1,\ldots, h_q)
\\
&&\hphantom{\sum_{q=0}^{+\infty} \int_{\R_+^q}\Biggl(\sum_{h_1,\ldots,
h_q}}
\qquad\quad{} \times\rme^{-r(m-1)t+\sum_{i=1}^q \log(h_i)}\prod
_{i=1}^q p_{h_i}
\ind_{\{\sum_{i=1}^q t_i\leq t\}} \Biggr)\,dt_1\cdots dt_q\\
&&\qquad= \E\bigl[g\bigl(\Lambda'_{[0,t]}\bigr)\rme^{-r(m-1)t + \Lambda'_t} \bigr].
\end{eqnarray*}
Recall
that $(S_t, t\geq0)$ [resp., $(S^u_t, t\geq0)$] is the underlying
Poisson process of $\Lambda$ (resp., $\Lambda^u$). Notice that if
$|u|=q$, then $\{S^u_t=q\}=\{\alpha(u)\leq t<\beta(u)\}$. We thus
deduce from (\ref{utgirsanov}) that for $q\in\N$, $u\in\cu$ such
that $|u|=q$,
%
%e27 ###
%
\begin{equation}
\label{eq:LY=Lu}
\E\bigl[g\bigl(\Lambda_{[0,t]}\bigr)\ind_{\{S_t=q\}}\bigr]=
\E\bigl[g\bigl(\Lambda^u_{[0,t]}\bigr)\rme^{-r(m-1)t + \Lambda^u_t} \ind
_{\{ \alpha(u)\leq
t<\beta(u)\}} \bigr].
\end{equation}
Let $q\in\N^*$. By construction, conditionally on
$\{\Lambda_{[0,t]}=\lambda_{[0,t]}\}$, $\{S_t=q\}$, $\{(I_1, \ldots, I_q
)=u\}$, $Y_{[0,t]} $ is distributed as $X^u_{[0,t]} $ conditionally on
$\{\Lambda^u_{[0,t]}=\lambda_{[0,t]}\}$. This\vspace*{1pt} holds also for $q=0$ with
the convention that $(I_1, \ldots, I_q )=\varnothing$.
Therefore, we have for any nonnegative measurable functions $g$ and $f$,
\begin{eqnarray*}
&&\E_\mu\bigl[g\bigl(\Lambda_{[0,t]}\bigr) f\bigl(Y_{[0,t]}\bigr)\bigr]\\
&&\qquad=\sum_{u\in\cu}\sum_{q\in\N} \ind_{\{|u|=q\}} \E_\mu
\bigl[g\bigl(\Lambda_{[0,t]}\bigr)
f\bigl(Y_{[0,t]}\bigr)\ind_{\{(I_1, \ldots,
I_{q})=u\}}\ind_{\{S_t=q\}} \bigr] \\
&&\qquad=\sum_{u\in\cu}\sum_{q\in\N} \ind_{\{|u|=q\}}
\E_\mu\bigl[g\bigl(\Lambda_{[0,t]}\bigr)\E_\mu\bigl[f\bigl(X^u_{[0,t]}\bigr)|
\Lambda^u_{[0,t]}\bigr]_{|\Lambda^u_{[0,t]}=\Lambda_{[0,t]}}
\\
&&\hspace*{162pt}\qquad\quad{}\times\ind_{\{(I_1, \ldots,
I_{q})=u\}}\ind_{\{S_t=q\}} \bigr].
\end{eqnarray*}
Using the points (i) and (ii) of Definition \ref{def:Y}, we see that
\begin{eqnarray*}
&&\PP(I_1=u_1,\ldots,I_q=u_q | \Lambda_t, \{S_t=q\})\\
&&\qquad= \prod_{k=1}^q 1/H_k= \rme^{-\Lambda_t} \qquad\mbox{if }
u_1\leq H_1,\ldots, u_q\leq H_q\\
&&\qquad= 0 \qquad\mbox{otherwise}.
\end{eqnarray*}
Hence,
\begin{eqnarray*}
&&\E_\mu\bigl[g\bigl(\Lambda_{[0,t]}\bigr) f\bigl(Y_{[0,t]}\bigr)\bigr]\\
&&\qquad=\sum_{u\in\cu}\sum_{q\in\N} \ind_{\{|u|=q\}}
\E_\mu\bigl[g\bigl(\Lambda_{[0,t]}\bigr)\E_\mu\bigl[f\bigl(X^u_{[0,t]}\bigr)|
\Lambda^u_{[0,t]}\bigr]_{|\Lambda^u_{[0,t]}=\Lambda_{[0,t]}}\\
&&\hspace*{121.6pt}\qquad\quad{} \times
\rme^{-\Lambda_t} \ind_{\{u_1\leq H_1,\ldots, u_q\leq H_q\}}\ind
_{\{S_t=q\}}\bigr] \\
&&\qquad=\sum_{u\in\cu} \sum_{q\in\N} \ind_{\{|u|=q\}}
\E_\mu\bigl[g\bigl(\Lambda^u_{[0,t]}\bigr)\E_\mu\bigl[f\bigl(X^u_{[0,t]}\bigr)|
\Lambda^u_{[0,t]}\bigr]\rme^{-r(m-1)t}\\
&&\hspace*{102.8pt}\qquad\quad{} \times
\ind_{\{u_1\leq\nu_\varnothing,\ldots, u_q\leq\nu_{(u_1, \ldots,
u_{q-1})} \}}\ind_{\{S^u_t=q\}}\bigr],
\end{eqnarray*}
thanks to (\ref{eq:LY=Lu}). Remark that for $u=(u_1 ,\ldots, u_q)$,
\[
\bigl\{u_1\leq\nu_\varnothing,\ldots, u_q\leq\nu_{(u_1 ,\ldots, u_{q-1})}
\bigr\}=\{u\in\mathcal{T}\},
\]
and for such $u$, we have $\{S^u_t=q\}=\{u\in V_t\}$ as noticed before.
As a consequence,
\[
\E_\mu\bigl[g\bigl(\Lambda_{[0,t]}\bigr) f\bigl(Y_{[0,t]}\bigr)\bigr]
=\sum_{u\in\cu}
\E_\mu\bigl[g\bigl(\Lambda^u_{[0,t]}\bigr)f\bigl(X^u_{[0,t]}\bigr)\rme^{-r(m-1)t}\ind_{\{
u\in V_t\}}\bigr].
\]
Finally, we use a monotone class
argument to conclude.
\end{pf*}

%s3.2 ###
\subsection{Many-to-One formulas over the whole tree}\label{sectionmainpties}

In this section we generalize identity (\ref{defYt}) on the link
between the tree indexed process $X_\T$ and the auxiliary Markov process
$Y$ by considering sums over the whole tree.

Let us consider the space $\mathcal{D}$ of nonnegative measurable
functions $f\in\mathcal{B}(\R_+\times\D(\R_+,E\times\R),\R_+)$
such that $f(t,y)=f(t,z)$ as soon as $y_{[0,t)}=z_{[0,t)}$. By
convention, if $y$ is defined at least on $[0,t)$, we will write
$f(t,y_{[0,t)})$ for $f(t,z)$ where $z$ is any function
such that $z_{[0,t)}=y_{[0,t)}$.\vspace*{-3pt}

\begin{prop}[(Many-to-One formula over the whole tree)]
\label{proptraj}
For all nonnegative measurable function $f$ of $\mathcal{D}$, we
have
%
%e28 ###
%
\begin{equation}
\label{formuledualite2}
\E_\mu\biggl[\sum_{u \in\ct} f\bigl(\beta(u), \w{X}^u_{[0,\beta(u))}\bigr)
\biggr]
=r\int_0^{+\infty} ds \,\rme^{r(m-1)s} \E_\mu\bigl[f\bigl(s, \w{Y}_{[0,s)}\bigr)
\bigr].\vspace*{-3pt}
\end{equation}
\end{prop}

Before coming to the proof of Proposition \ref{proptraj}, we introduce
a notation that will be very useful in the sequel.
By convention for two functions $f,g$ defined, respectively, on two
intervals $I_f, I_g$, for $[a,b)\subset I_f$ and $[c,d)\subset I_g$, we
define the concatenation $[f_{[a,b)};g_{[c,d)}]=h_{J}$ where
$J=[a,b+(d-c))$,
\[
h(t)=\cases{
f(t), &\quad if $t\in[a,b)$,\cr
g\bigl(c+(t-b)\bigr), &\quad if $t\in
[b, d-c+b)$.}\vspace*{-3pt}
\]
\begin{pf*}{Proof of Proposition \ref{proptraj}}
We first notice that if $\tau$ is an exponential random
variable with mean $1/r$ ($r>0$), then we have, for any nonnegative
measurable function $g$,
%
%e29 ###
%
\begin{equation}\label{eq:moyenneexp}
\E\biggl[r \int_0^{\tau} g(t) \,dt \biggr] = \E[g(\tau)].
\end{equation}
Besides, we have
\[
\E_\mu\bigl[ \ind_{\{u\in\mathcal{T}\}}f\bigl(\beta(u), \w
{X}^u_{[0,\beta(u))}\bigr) \bigr]= \E_\mu\bigl[ \ind_{\{u\in\mathcal
{T}\}}f\bigl(\beta(u),
\bigl[\w{X}^u_{[0,\alpha(u))};\w{X}_{[0,\beta(u)-\alpha(u))}\bigr]\bigr)\bigr],
\]
where conditionally on $\w{X}^u_{[0,\alpha(u))}$, $\beta(u)$, $\{
u\in\mathcal{T}\}$, $\w{X}=(X,c)$ with $X$ of distribution $\PP
_{X^u_{\alpha(u)}}$ and $c$ the constant process equal to $\Lambda
^u_{\alpha(u)}$. Notice that we have chosen $\w{X}$ independent of
$\beta(u)$.
Thus, conditioning with respect to $[\w{X}^u_{[0,\alpha(u))}; \w
{X}_{[0,+\infty)}]$, $\{u\in\ct\}$ and using (\ref{eq:moyenneexp}),
we get
% \E_\mu[ \ind_{\{u\in\mathcal{T}\}}\E[f(\beta(u), [
% and $\w{X}^{\w{X}^u_{\alpha(u)}}_{[0,+\infty)}$, the integrand is a
%function of the sole variable $\beta(u)$. Applying (
%
\begin{eqnarray*}
&&\E_\mu\bigl[ \ind_{\{u\in\mathcal{T}\}}f\bigl(\beta(u), \w
{X}^u_{[0,\beta(u))}\bigr) \bigr]\\[-2pt]
&&\qquad=r\E\biggl[\ind_{\{u\in\mathcal{T}\}} \int_{0}^{\beta(u) -\alpha
(u)} \,ds \,f\bigl(\alpha(u)+s, \w{X}^u_{[0,\alpha(u)+s)}\bigr) \biggr].
\end{eqnarray*}
We deduce,
\begin{eqnarray*}
\E_\mu\bigl[ \ind_{\{u\in\mathcal{T}\}}f\bigl(\beta(u), \w
{X}^u_{[0,\beta(u))}\bigr) \bigr]
&=&r\E\biggl[\ind_{\{u\in\mathcal{T}\}} \int_{\alpha(u)}^{\beta
(u) } ds \,f\bigl(s, \w{X}^u_{[0,s)}\bigr) \biggr]\\[-2pt]
&=&r\int_0^{+\infty} ds \,\E\bigl[\ind_{\{u\in V_s\}} f\bigl(s,
\w{X}^u_{[0,s)}\bigr) \bigr].
\end{eqnarray*}
%
%where we used \reff{eq:moyenneexp} for the first equality.
% and that $\w{X}^u_{[0,\beta(u))}$ depends on $\beta(u)$ only through
%the interval $[0, \beta(u))$ on which the function $\w{X}^u$ is
%defined.
Using Proposition \ref{propprocessusauxiliaire}, we get
\begin{eqnarray*}
\E_\mu\biggl[\sum_{u \in\ct} f\bigl(\beta(u), \w{X}^u_{[0,\beta(u))}\bigr)
\biggr]
&=& r\int_0^{+\infty} ds \,\E_\mu\biggl[\sum_{u\in V_s} f\bigl(s,
\w{X}^u_{[0,s)}\bigr) \biggr]\\
&=&r\int_0^{+\infty} ds \,\rme^{r(m-1)s} \E_\mu\bigl[f\bigl(s, \w{Y}_{[0,s)}\bigr)
\bigr].
\end{eqnarray*}
\upqed\end{pf*}

The equality (\ref{formuledualite2}) means that adding the contributions
over all the individuals in the Galton--Watson tree corresponds (at least
for the first moment) to integrating the contribution\vspace*{1pt} of the auxiliary
process over time with an exponential weight $\rme^{r(m-1)t}$ which is
the average number of living individuals at time $t$. Notice\vadjust{\goodbreak} the weight is
increasing if the Galton--Watson tree is supercritical and decreasing if
it is subcritical.
%The left hand side of (\ref{formuledualite2}) corresponds
%heuristically to picking an individual uniformly among all the
%individuals of all the possible trees.
%
\begin{rque}
We shall give two alternative formulas for
(\ref{formuledualite2}).
\begin{longlist}[(1)]
\item[(1)]
We deduce from (\ref{formuledualite2}) that, for
all nonnegative measurable function $f$,
%
%e30 ###
%
\begin{equation}
\label{eq:formuledualite3}
\E_\mu\biggl[\sum_{u \in\ct} f\bigl(\beta(u), \w{X}^u_{[0,\beta(u))}\bigr)
\biggr]
=\E_\mu\bigl[f\bigl(\tau, \w{Y}_{[0,\tau)}\bigr)\rme^{rm \tau} \bigr],
\end{equation}
where $\tau$ is an independent exponential random variable of mean
$1/r$. Thus, the right-hand side of equation
(\ref{formuledualite2}) can be read as the expectation of a~functional
of the process $\w{Y}$ up to an independent exponential time
$\tau$ of mean $1/r$ with a weight $\rme^{rm\tau}$.
\item[(2)] Let $\tau_q=\inf\{t\geq0; S_t=q\} $ the time of the $q$th jump
for the compound Poisson process $\Lambda$. Using
(\ref{eq:moyenneexp}), it is easy to check that, for any
nonnegative measurable function $g$,
\[
\frac{1}{m} \sum_{q\geq1}\E_\mu\bigl[g\bigl(\w{Y}_{[0,\tau_q)}, \tau
_q\bigr)\bigr]=r\int_0^{+\infty
} \E_\mu\bigl[g\bigl(\w{Y}_{[0,s)}, s\bigr)\bigr] \,ds.
\]
Therefore, we deduce from (\ref{formuledualite2}) that, for
all nonnegative measurable function~$f$,
%
%e31 ###
%
\begin{equation}
\label{eq:formuledualite4}
\E_\mu\biggl[\sum_{u \in\ct} f\bigl(\beta(u), \w{X}^u_{[0,\beta(u))}\bigr)
\biggr]
=\frac{1}{m} \sum_{q\geq1}\E_\mu\bigl[f\bigl(\tau_q, \w{Y}_{[0,\tau
_q)}\bigr)\rme^{r(m -1)\tau_q} \bigr].
\end{equation}
This formula emphasizes that the jumps of the auxiliary process
correspond to death times in the tree.
\end{longlist}
\end{rque}

%s3.3 ###
\subsection{Identities for forks}
\label{identititesforforks}
In order to compute second moments, we shall need the distribution of two
individuals\vadjust{\goodbreak} picked at random in the whole population and which are not
in the same lineage. As in the
Many-to-One formula, it will involve the auxiliary process.

First, we define the following sets of forks:
%
%e32 ###
%
\begin{equation}\label{defforks}\quad
\mathcal{FU}= \{(u,v) \in\mathcal{U}^2 \dvtx|u \wedge v| <
\min(|u|,|v|)\} \quad\mbox{and}\quad \mathcal{FT}= \mathcal{FU}
\cap\mathcal{T}^2.
\end{equation}
Let $\w{J}_2$ be the operator defined for all nonnegative
measurable function $f$ from $(E\times\R)^2$ to $\R$ by
%
%e33 ###
%
\begin{eqnarray}
\label{defJ2}
\w{J}_2 f(x,\lambda)
&=&\int_0^1 \mathop{\sum_{(a,b) \in(\N
^*)^2}}_{a \neq b}
\sum_{k\geq\max(a,b)}p_k
f\bigl(F^{(k)}_a(x,\theta),\lambda+\log(k),\nonumber\\[-8pt]\\[-8pt]
&&\hphantom{\int_0^1 \mathop{\sum_{(a,b) \in(\N
^*)^2}}_{a \neq b}
\sum_{k\geq\max(a,b)}p_k
f\bigl(}\hspace*{3pt}
F^{(k)}_b(x,\theta),
\lambda+
\log(k)\bigr)\,d\theta.
\nonumber
\end{eqnarray}
Informally, the functional $\w{J}_2$ describes the starting positions
of two siblings.
Notice that we have
%
%e34 ###
%
\begin{eqnarray}
\label{defJ2-2}
\w{J}_2 f(x,\lambda)
&=& m \int_0^1 \E\bigl[(H-1)f\bigl(F^{(H)}_I(x, \theta),
\lambda+\log(H) ,\nonumber\\[-8pt]\\[-8pt]
&&\hphantom{m \int_0^1 \E\bigl[(H-1)f\bigl(}\hspace*{2.6pt}
F^{(H)}_K(x,
\theta), \lambda+\log(H) \bigr)\bigr]\,d\theta,
\nonumber
\end{eqnarray}
where $H$ has the size-biased offspring distribution and conditionally
on $H$, $(I, K)$ is distributed as a drawing of a couple without
replacement among
the integers $\{1, \ldots, H\}$.

For measurable real functions $f$ and $g$ on $E\times\R$, we
denote by $f\otimes g$ the real measurable function on
$(E\times\R) ^2$ defined by: $(f\otimes
g)(\w{x},\w{y})=f(\w{x})g(\w{y})$ for $\w{x},\w{y}\in E\times
\R$.
\begin{prop}[(Many-to-One formula for forks over the whole tree)]
\label{fourche2bis}
For all nonnegative measurable functions $\varphi,\psi\in\mathcal
{D}$, we
have
%
%e35 ###
%
\begin{eqnarray}\label{tatalo3bis}
&&\E_\mu\biggl[\sum_{(u,v) \in\mathcal{FT}}\varphi\bigl(\beta(u),\w{X}^{u}_{[0,
\beta(u)) }\bigr)\psi\bigl( \beta(v),\w{X}^v_{[0,\beta(v))}\bigr)\biggr] \nonumber
\\
&&\qquad= \E_\mu\bigl[ \rme^{rm\tau}
\w{J}_2\bigl({\E}_{\cdot}'\bigl[\varphi\bigl(t+{\tau'}, \bigl[\w{y}_{[0,t)};
\w{Y}'_{[0,{\tau'})}\bigr] \bigr) \rme^{rm{\tau'}} \bigr]_{|t=\tau, \w
{y}=\w{Y}}\\
&&\hphantom{\E_\mu\bigl[ \rme^{rm\tau}
\w{J}_2\bigl(}
\qquad\quad{}\otimes \E_{\cdot}'\bigl[\psi\bigl(t+{\tau'}, \bigl[\w{y}_{[0,t)};
\w{Y}'_{[0,{\tau'})}\bigr]\bigr) \rme^{rm{\tau'}}\bigr]_{|t=\tau, \w{y}=\w{Y}}
\bigr)(\w{Y}_{\tau-})\bigr],\nonumber
\end{eqnarray}
where, under $\E_\mu$, $\tau$ is exponential with mean $1/r$
independent of $\w{Y}$ and, under $\E'_{x,\lambda}$, $(\w{Y}',
{\tau'})$ is distributed as
$(( Y,\Lambda+\lambda), \tau)$ under $\E_x$.
\end{prop}
\begin{pf}
Notice that $\{(u,v) \in\mathcal{FU}\}$ is equal to $\{\exists(w,
\widetilde u, \widetilde v) \in\mathcal{U}^3, \exists(a, b) \in(\N^*)^2, a
\neq b, u=wa\widetilde{u}, v=wb\widetilde{v}\}$. Let $A$
be the left-hand side of (\ref{tatalo3bis}). We have
\begin{eqnarray*}
A&=& \sum_{w\in\U}\mathop{\sum_{a,b\in\N^*}}_{
a\not=b}\sum_{\w{u},\w{v}\in\U}
\E_\mu
\bigl[\varphi\bigl(\beta(w)+\bigl(\beta(wa\w{u})-\beta(w)\bigr),\\
&&\hphantom{\sum_{w\in\U}\mathop{\sum_{a,b\in\N^*}}_{
a\not=b}\sum_{\w{u},\w{v}\in\U}
\E_\mu
\bigl[\varphi\bigl(}
\hspace*{3.4pt}\bigl[
\w{X}^{w}_{[0,\beta(w))}; \w{X}^{wa\widetilde{u}}_{[\beta(w),
\beta(wa\w{u}))}\bigr] \bigr)\ind_{\{wa\w{u}\in\mathcal{T}\}} \\
&&\hphantom{\sum_{w\in\U}\mathop{\sum_{a,b\in\N^*}}_{
a\not=b}\sum_{\w{u},\w{v}\in\U}
\E_\mu
\bigl[}
{}\times
\psi\bigl(\beta(w)+\bigl(\beta(wb\w{v})-\beta(w)\bigr), \\
&&\hspace*{120pt}\bigl[\w
{X}^{w}_{[0,\beta(w))
};\w{X}^{wb\widetilde{v}}_{[\beta(w),
\beta(wb\w{v}))}\bigr] \bigr)\ind_{\{wb\w{v}\in\mathcal{T}\}} \bigr].
\end{eqnarray*}
Using the strong Markov property at time $\beta(w)$, the conditional
independence between descendants and Proposition \ref{proptraj}, we get
%
%e36 ###
%
\begin{eqnarray}\label{etapediffpreuve}
A
&=& \sum_{w\in\U}\mathop{\sum_{a,b\in\N^*}}_
{a\not=b}\E_\mu
\bigl[\E'_{\w{X}
^{wa}_{\alpha(wa)}}\bigl[\varphi\bigl(t+{\tau'}, \bigl[\w{x}_{[0,t)}
; \w{Y}_{[0, {\tau'})} \bigr]\bigr)\rme^{rm{\tau'}}\bigr]_{| t=\beta(w),
\w{x}=\w{X}^w} \ind_{\{wa\in\mathcal{T}\}} \nonumber\\
&&\hspace*{62.01pt}{}\times\E'_{\w{X}
^{wb}_{\alpha(wb)}}\bigl[\psi\bigl(t+{\tau'}, \bigl[\w{x}_{[0,t)}
; \w{Y}_{[0, {\tau'})} \bigr]\bigr)\rme^{rm{\tau'}}\bigr]_{| t=\beta(w),
\w{x}=\w{X}^w} \\
&&\hspace*{255.3pt}{}\times\ind_{\{wb\in\mathcal{T}\}}\bigr],
\nonumber\vadjust{\goodbreak}
\end{eqnarray}
where under $\E'_{x,\lambda}$, $(\w{Y}',
{\tau'})$ is distributed as
$(( Y,\Lambda+\lambda), \tau)$ under $\E_x$.
As $\{wa,wb\in\mathcal{T}\}=\{w\in\mathcal{T}\}\cap\{\max\{a,b\}
\leq
\nu_w\}$ we have
%
%e37 ###
%
\begin{eqnarray} \label{etape3}\qquad
A&=&\sum_{w\in\U} \E_\mu\bigl[ \ind_{\{w\in\mathcal{T}\}}
\w{J}_2\bigl({\E}_{\cdot}'\bigl[\varphi\bigl(t+{\tau'}, \bigl[\w{x}_{[0,t)};
\w{Y}'_{[0,{\tau'})}\bigr] \bigr) \rme^{rm{\tau'}} \bigr]_{|t=\beta(w),
\w{x}=\w{X}^w}\nonumber\\
&&\hspace*{79.5pt}{}\otimes{\E}_{\cdot}'\bigl[\psi\bigl(t+{\tau'}, \bigl[\w{x}_{[0,t)};
\w{Y}'_{[0,{\tau'})}\bigr] \bigr)\rme^{rm{\tau'}}\bigr]_{|t=\beta(w), \w
{x}=\w{X}^w}
\bigr)\\
&&\hspace*{250.1pt}{}\times\bigl(\w{X}^w_{\beta(w)_-}\bigr)\bigr]
\nonumber
\end{eqnarray}
with $\w{J}_2$ defined by (\ref{defJ2}). The function under the expectation
in (\ref{etape3}) depends on $\beta(w)$ and
$\w{X}^w_{[0,\beta(w))}$. Equality
(\ref{eq:formuledualite3}) then gives the
result.
\end{pf}

We shall give a version of Proposition \ref{fourche2bis} when the
functions of the path depend only\vspace*{1pt} on the terminal value of the path. We
shall define $J_2$ as a~simpler version of $\w{J}_2$ [see definition
(\ref{defJ2-2})] acting\vspace*{1pt} only on the spatial motion; for all nonnegative
measurable function $f$ from $E^2$ to $\R$,
%
%e38 ###
%
\begin{equation}
\label{eq:J2-0}
J_2 f(x)= m \int_0^1 \E\bigl[ (H-1)f\bigl(F^{(H)}_I(x, \theta),F^{(H)}_K(x,
\theta)\bigr)\bigr]\,d\theta,
\end{equation}
where $(H,I,K)$ are as in (\ref{defJ2-2}).

As a direct consequence of Proposition
\ref{fourche2bis} and of the fact that $Y$ is c\`{a}dl\`{a}g, we have
the following corollary.
\begin{corol}[(Many-to-One formula for forks over the whole tree)]
\label{cor:fourche2bis}
Let $({Q}_t, t\geq0)$ be the transition semi-group of ${Y}$.
For all nonnegative measurable functions $f,g\in\mathcal{D}$, we
have
%
%e39 ###
%
\begin{eqnarray}\label{tatalo2bis}
&&\E_\mu\biggl[\sum_{(u,v) \in\mathcal{FT}}
f\bigl(\beta(u),{X}^{u}_{\beta(u)- }\bigr) g\bigl( \beta(v),{X}^v_{\beta(v)-}\bigr)\biggr]
\nonumber\\[-8pt]\\[-8pt]
&&\qquad=r^3\int_{[0,\infty)^3} \rme^{r(m-1) (s+t+t')} \,ds \,dt \,dt'\, \mu{Q}_s
\bigl(J_2({Q}_t
f_{t+s} \otimes{Q}_{t'}
g_{t'+s})\bigr),
\nonumber
\end{eqnarray}
where $f_t(x)=f(t,x)$ and $g_t(x)=g(t,x)$ for $t\geq0$ and $x\in
E$.
\end{corol}

We can also derive a Many-to-One formula for forks at fixed
time.
\begin{prop}[(Many-to-One formula for forks at fixed time)]
\label{fourche2}
Let $t\in\R_+$ and $\varphi, \psi$ be two nonnegative measurable
functions on $\mathbb{D}([0,t],E)$. Then
%
%e40 ###
%
\begin{eqnarray}\label{tatalo2}
\hspace*{1pt}&&\E_\mu\biggl[\mathop{\sum_{(u,v) \in V_t^2}}_{u \not= v} \varphi
\bigl(\w{X}^{u}_{[0,t]}\bigr) \psi\bigl(\w{X}^v_{[0,t]}\bigr)\biggr]
\nonumber\\
\hspace*{1pt}&&\qquad=r\rme^{2r(m-1)t} \int_0^t\hspace*{-2pt} da \bigl( \rme^{-r(m-1)a}\E_\mu\bigl[
\w{J}_2\bigl({\E}_{\cdot}'\bigl[\varphi\bigl(\bigl[\w{y}_{[0,a)};
\w{Y}'_{[0,t-a]}\bigr] \bigr) \bigr]_{|\w{y}=\w{Y}} \\
\hspace*{1pt}&&\qquad\quad\hspace*{149pt}{}
\otimes{\E}_{\cdot}'\bigl[\psi\bigl(\bigl[\w{y}_{[0,a)};
\w{Y}'_{[0,t-a]}\bigr]\bigr) \bigr]_{| \w{y}=\w{Y}}
\bigr)(\w{Y}_{a})\bigr]\bigr),
\nonumber
\end{eqnarray}
where, under $\E'_{x,\lambda}$, $\w{Y}'$ is distributed as
$( Y,\Lambda+\lambda)$ under $\E_x$.
\end{prop}

The left-hand side of (\ref{tatalo2}) approximates the distribution of a
pair of individuals uniformly chosen from the population at time $t$.
Indeed, we have in the right-hand side of
(\ref{tatalo2}) an exponential weight $\rme^{2r(m-1)t}$ and thanks to
Lemma \ref{lemmeNt}, we know that $\E[N_t(N_t-1)]\sim C \rme^{2r(m-1)t}$.
%which is the dominant term in the average number of pairs of
% individuals that can be picked at time $t$.
The distribution of the
paths associated with a random pair is described by the law of
forks constituted of independent portions of the auxiliary process
$\w{Y}$ and splitted at a time $a\in[0,t]$. Notice that (\ref{tatalo2})
indicates that the fork splits at an exponential random time with
mean $1/r(m-1)$, conditioned to be less than~$t$.

\begin{pf*}{Proof of Proposition \ref{fourche2}}
The proof is similar to the proof of Proposition \ref{fourche2bis}
except that we use Proposition \ref{propprocessusauxiliaire} instead
of Proposition \ref{proptraj} to obtain an analogue of
(\ref{etapediffpreuve}).
\end{pf*}

%s4 ###
\section{Law of large numbers}\label{sectionlgn}
In this section, we are interested in averages over the population living
at time $t$ for large $t$. When the Galton--Watson tree is not supercritical
we have almost sure extinction and thus, we assume here that $m>1$.

%s4.1 ###
\subsection{Results and comments}\label{sec4.1}

Notice that $N_t=0$ implies $Z_t=0$ and by convention we set
$Z_t/N_t=0$ in this case.\vadjust{\goodbreak}
For $t\in\R_+$ and $f$ a real function defined on~$E$, we derive
laws of large numbers for
%
%e41 ###
%
\begin{equation}
\frac{\langle
Z_t,f\rangle}{N_t}= \frac{\sum_{u\in V_t} f(X_t^u)}{N_t}
\quad\mbox{and}\quad \frac{\langle Z_t,f\rangle}{\E[N_t]}=\frac{\sum
_{u\in
V_t} f(X_t^u)}{\E[N_t]},
\end{equation}
provided the
auxiliary process introduced in the previous section satisfies some
ergodic conditions.

Let $(Q_t, t\geq0) $ be the semigroup of the auxiliary process $Y$
of Definition~\ref{def:Y},
%
%e42 ###
%
\begin{equation}\label{eq:Qt}
\E_\mu[f(Y_t)]=\mu Q_t f
\end{equation}
for all $\mu\in
\cp(E)$ and $f$ nonnegative. Recall the operator $J_2$ defined in
(\ref{eq:J2-0}).

We shall consider the following ergodicity and integrability
assumptions on $f$, a real
measurable function defined on $E$, and on $\mu\in\cp(E)$.
\begin{longlist}[(H4)]
\item[(H1)] There exists a nonnegative finite measurable function $g$
such that $Q_t|f|(x) \le g(x)$ for all $t \ge0$ and $x \in E$.
\item[(H2)] There exists $\pi\in\cp(E)$ such that $\langle\pi, |f|
\rangle<+\infty$ and for all $x\in E$, $\lim_{t \rightarrow
+\infty} Q_tf(x) = \langle\pi
, f \rangle$.
\item[(H3)] There exists $\alpha< r(m-1)$ and $c_1>0$ such that
$\mu Q_tf^2 \le c_1 \rme^{\alpha t}$ for every $t \ge0$.
\item[(H4)] There exists $\alpha< r(m-1)$ and $c_2>0$ such that
$\mu Q_tJ_2(g \otimes g) \le c_2 \rme^{\alpha t}$ for every $t \ge0$
with $g$ defined in (H1).
\end{longlist}
Notice that in (H3) and (H4), the constants $\alpha, c_1$ and $c_2$ may depend
on $f$ and $\mu$.
\begin{rque}\label{rquecontientcontinuborne}
When the auxiliary process $Y$ is ergodic [i.e., $Y$ converges in
distribution to $\pi\in\cp(E)$], the class of continuous bounded
functions satisfies (H1)--(H4) with $g$ constant and $\alpha=0$. In some
applications, one may have to consider polynomially growing functions. This
is why we shall consider hypotheses (H1)--(H4) instead of the ergodic
property in Theorem \ref{thLGN} or in Proposition \ref{cvL2ssarbre}.
\end{rque}

The next theorem states the law of large numbers; the asymptotic
empirical measure is distributed as the stationary distribution $\pi$
of $Y$.
\begin{theorem}\label{thLGN}
For any $\mu\in\mathcal{P}(E)$ and any real measurable function $f$
defined on $E$ satisfying \textup{(H1)}--\textup{(H4)}, we have
%
%e44 ###
%e43 ###
%
\begin{eqnarray}
\label{eq:lln-l2}
\lim_{t\rightarrow+\infty} \frac{\langle
Z_t,f\rangle}{\E[N_t]}
&=&\langle\pi, f \rangle W
\qquad\mbox{in $L^2(\PP_\mu)$,} \\
\label{eq:lln-P}
\lim_{t\rightarrow+\infty} \frac{\langle
Z_t,f\rangle}{N_t}
&=& \langle\pi, f \rangle\ind_{\{W
\neq0\}} \qquad\mbox{in $\PP_\mu$-probability}
\end{eqnarray}
with $W$ defined by (\ref{eq:CV-N-W}) and $\pi$ defined in \textup{(H2)}.
\end{theorem}

For the proof which is postponed to Section \ref{sec:proof}, we use
ideas developed in~\cite{delmasmarsalle} in a discrete time setting.
We give an intuition of the result. According to Proposition
\ref{propprocessusauxiliaire}, an individual chosen at random at time
$t$ is heuristically distributed as~$Y_t$, that is, as $\pi$ for large
$t$ thanks to the ergodic property of~$Y$ [see (H2)]. Moreover, two
individuals chosen at random among the living individuals at time $t$
have a MRCA who died early which implies that they behave almost
independently. Since Lemma \ref{lemmeNt} implies that the number of
individuals alive at time $t$ grows to infinity on $\{W\neq0\}$,
this yields the LLN stated in Theorem \ref{thLGN}.

Notice that Theorem \ref{thLGNannonceintro} is a direct consequence of
Theorem \ref{thLGN} and Remark~\ref{rquecontientcontinuborne}.

We also present a LLN when summing over the set of
all individuals who died before time $t$. Recall that
$D_t=\sum_{u\in\mathcal{T}}\ind_{\{ \beta(u)<
t\}} $ denotes its cardinal.

Recall $S$ in Definition \ref{def:Y}. Notice that $\E[S_t]=rmt$. We
shall consider a~slightly stronger hypothesis than (H3):
\begin{longlist}[(H5)]
\item[(H5)] There exists $\alpha< r(m-1)$ and $c_3>0$ such that
$\E_\mu[f^2(Y_t) S_t] \le c_3 \rme^{\alpha t}$ for every $t
\ge0$.\vspace*{-3pt}
\end{longlist}
\begin{prop} \label{cvL2ssarbre}
For any $\mu\in\mathcal{P}(E)$ and any nonnegative measurable
function $f$
defined on $E$ satisfying \textup{(H1)}--\textup{(H5)}, we have
%
%e46 ###
%e45 ###
%
\begin{eqnarray}\quad
\label{eq:path-lln-l2}
\lim_{t\rightarrow+\infty}\frac{\sum_{u\in
\mathcal{T}}f(X^u_{\beta(u)-} )\ind_{\{ \beta(u)<
t\}}}{\E[D_t]}
&=&\langle\pi, f \rangle W \qquad\mbox{in $
L^2(\PP_\mu)$}, \\[-2pt]
\lim_{t\rightarrow+\infty}\ind_{\{N_t>0\}} \frac{\sum_{u\in
\mathcal{T}}f(X^u_{\beta(u)-} )\ind_{\{ \beta(u)<
t\}}}{D_t}
&=&\langle\pi, f \rangle\ind_{\{W\neq0\}} \nonumber\\[-10pt]\\[-10pt]
&&\eqntext{\mbox{in $\PP
_\mu
$-probability}}
\end{eqnarray}
with $W$ defined by (\ref{eq:CV-N-W}) and $\pi$ defined in
\textup{(H2)}.\vspace*{-3pt}
\end{prop}

We can then extend these results to path dependent functions. In
particular, the next theorem describes the asymptotic distribution of
the motion and lineage of an individual taken at random in the tree. In
order to avoid a set of complicated
hypothesis we shall assume that $Y$ is ergodic with limit distribution
$\pi$ and consider bounded functions.\vspace*{-3pt}
\begin{theorem}
\label{thLGNfonctionnel}
We assume that there exists $\pi\in\cp(E)$ such that for all $x\in E$,
and all real-valued bounded measurable function $f$ defined on $E$,
$\lim_{t\rightarrow\infty} Q_tf(x)=\langle\pi, f \rangle$.
Let $T > 0$. For any real bounded measurable function $\varphi$ on
$\D([0,T],E\times\R_+)$, we have
\begin{eqnarray*}
&& \lim_{t\rightarrow\infty} \frac{1}{\E[N_{t}]}\sum_{u\in
V_{t}}\varphi\bigl(X^u_{[t-T,t]}, \Lambda^u_{[t-T,t]} -\Lambda^u_{t-T}\bigr)
\\[-2pt]
&&\qquad=
\E_\pi\bigl[\varphi\bigl(\w{Y}_{[0,T]}\bigr) \bigr] W
\qquad\mbox{in $L^2(\PP_\mu)$,}\\
&& \lim_{t\rightarrow\infty} \frac{1}{N_{t}}\sum_{u\in
V_{t}}\varphi\bigl(X^u_{[t-T,t]}, \Lambda^u_{[t-T,t]} -\Lambda^u_{t-T}\bigr)\\
&&\qquad=
\E_\pi\bigl[\varphi\bigl(\w{Y}_{[0,T]}\bigr) \bigr]\ind_{\{W \neq0\}}
\qquad\mbox{in $\PP_\mu$-probability}
\end{eqnarray*}
with $W$ defined by (\ref{eq:CV-N-W}).
\end{theorem}

Let ${J}_1$ be the following operator associated with the possible
jumps of $Y$: for all nonnegative
measurable function $f$ from $E$ to $\R$,
%
%e47 ###
%
\begin{equation}
\label{eq:J1}
J_1f(x)= m \int_0^1 \E\bigl[f\bigl(F^{(H)}_I(x, \theta)\bigr) \bigr]\,d \theta,
\end{equation}
where $H$ has the size-biased offspring distribution and, conditionally
on $H$, $I$ is uniform on $\{1, \ldots, H\}$.
\begin{prop} \label{prop:cvL2ssarbre-path}
We assume that there exists $\pi\in\cp(E)$ such that for all $x\in E$,
and all real-valued bounded measurable function $f$ defined on $E$,
$\lim_{t\rightarrow\infty} Q_tf(x)=\langle\pi, f \rangle$.

Let $\varphi$ be a real bounded measurable function
defined on $E$-valued paths. We set for $x\in E$, $f(x)=\E_x[\varphi(Y_{[0,
\tau_1)})]$, with $\tau_1$ from Definition \ref{def:Y}. Then,
%
%e48 ###
%
\begin{eqnarray}\label{eq:all-path-lln-l2}
\lim_{t\rightarrow+\infty}\frac{\sum_{u\in
\mathcal{T}}\varphi(X^u_{[\alpha(u), \beta(u))} )\ind
_{\{ \beta(u)<
t\}}}{\E[D_t]}
&=& \langle\pi, J_1 f \rangle W \quad\mbox{in $
L^2(\PP_\mu)$},\nonumber\\
\hspace*{25pt}\lim_{t\rightarrow+\infty}\ind_{\{N_t>0 \}} \frac{\sum_{u\in
\mathcal{T}}\varphi(X^u_{[\alpha(u), \beta(u))} )\ind
_{\{ \beta(u)<
t\}}}{D_t}
&=& \langle\pi, J_1 f \rangle\ind_{\{W\neq0\}} \\
\eqntext{\mbox{in $\PP_\mu$-probability}}\quad
\end{eqnarray}
with $W$ defined by (\ref{eq:CV-N-W}).
\end{prop}
\begin{rque}
The hypothesis on $Y$ in Theorem \ref{thLGNfonctionnel} and
Proposition \ref{prop:cvL2ssarbre-path} is slightly stronger than the
ergodic condition (i.e., $Y$ converges in distribution to $\pi$), but it
is fulfilled if $Y$
converges to $\pi$ for the distance in total variation [i.e., for all
$x\in E$,
$\lim_{t\rightarrow\infty} \sup_{A\in\ce} |\PP_x(Y_t\in A)
-\pi(A)|=0$]. This property is very common for ergodic processes.
\end{rque}

%s4.2 ###
\subsection{Proofs}
\label{sec:proof}

\mbox{}

\begin{pf*}{Proof of Theorem \ref{thLGN}}
We assume (H1)--(H4). We shall first pro\-ve~(\ref{eq:lln-l2}) for $f$ such
that $\langle\pi, f\rangle=0$.
We have
\[
\E_\mu\biggl[ \frac{\langle
Z_t,f\rangle^2}{\E[N_t]^2}\biggr]=A_t+B_t,
\]
where
\[
A_t= \E[N_t]^{-2} \E_\mu\biggl[\sum_{u\in V_t}
f^2(X_t^u)\biggr]\vadjust{\goodbreak}
\]
and
\[
B_t = \E[N_t]^{-2} \E_\mu\biggl[\mathop{\sum_{(u,v) \in V_t^2}}_{u
\ne
v} f(X_t^u)f(X_t^v)\biggr].
\]
Notice that
%
%e49 ###
%
\begin{equation}\label{carre}
A_t
= \rme^{-r(m-1)t} \E_\mu[ f^2(Y_t) ]
=\rme^{-r(m-1)t} \mu Q_t f^2 \mathop{\longrightarrow
}_{t\rightarrow\infty} 0
\end{equation}
thanks to (\ref{equivdeterministe}) and (\ref{defYt}) for the first
equality and (H3) for the convergence. We focus now on $B_t$.
Proposition \ref{fourche2} and then (H1) and (H4) imply that
\begin{eqnarray*}
&&\E[N_t]^{-2} \E_\mu\biggl[ \mathop{\sum_{(u,v) \in V_t^2}}_{u \ne
v} |f(X_t^u)f(X_t^v)| \biggr]\\
&&\qquad= r \int_0^t \mu Q_sJ_2(Q_{t-s}|f| \otimes Q_{t-s}|f|)
\rme^{-r(m-1)s}\,ds
\end{eqnarray*}
is finite and
%
%e50 ###
%
\begin{equation}
\label{eq:Btcv}
B_t = r \int_0^t \mu Q_sJ_2(Q_{t-s}f \otimes Q_{t-s}f)
\rme^{-r(m-1)s}\,ds.
\end{equation}
Now, since $\langle\pi,f \rangle=0$, we deduce from (H2) that for $s$
fixed and \mbox{$y,z\in E$}, $\lim_{t \rightarrow\infty}(Q_{t-s}f \otimes
Q_{t-s}f)(y,z)=0$. Thanks to (H1), there exists $g$ such that
$\ind_{\{s\leq t\}} |(Q_{t-s}f \otimes Q_{t-s}f)| \le(g \otimes
g)$ and (H4) implies\vspace*{1pt} the finiteness of\break $\int_0^\infty ds\, \rme
^{-r(m-1) s} \mu
Q_sJ_2(g \otimes g) $. Lebesgue's theorem entails that
\[
\lim_{t \rightarrow\infty} B_t=\lim_{t \rightarrow\infty} r
\int_0^t
\mu Q_sJ_2(Q_{t-s}f \otimes Q_{t-s}f)\rme^{-r(m-1)s}\,ds =0.
\]
This ends the proof of (\ref{eq:lln-l2}) when $\langle
\pi,f \rangle=0$.

In the general case we have
%
%e51 ###
%
\begin{equation}\label{othercases}
\frac{\langle
Z_t,f\rangle}{\E[N_t]} - \langle\pi,f \rangle W = \frac{\langle
Z_t,f-\langle\pi, f \rangle\rangle}{\E[N_t]} +\langle\pi,f
\rangle
\biggl( \frac{N_t}{\E[N_t]}-W \biggr).
\end{equation}
Notice that if $f$ and $\mu$ satisfy (H1)--(H4) then so do $f-\langle
\pi
,f \rangle$ and $\mu$. The first term of the sum in the right-hand side of
(\ref{othercases}) converges
to $0$ in $L^2$ thanks to the first part of the proof. The second term
converges to $0$ in~$L^2$ thanks to Lemma \ref{lemmeNt}. Hence, we get
(\ref{eq:lln-l2}) if $f$ and $\mu$ satisfy (H1)--(H4). Equation~(\ref
{eq:lln-P}) stems from (\ref{eq:lln-l2}) and (\ref{eq:CV-N-W}).
\end{pf*}
\begin{pf*}{Proof of Proposition \ref{cvL2ssarbre}}
We assume (H1)--(H5). We shall first prove (\ref{eq:path-lln-l2}) for
$f$ such
that $\langle\pi, f\rangle=0$. We have
\[
\E[D_t]^{-2}\E_\mu\biggl[\biggl(\sum_{u\in
\mathcal{T}}f\bigl(X^u_{\beta(u)-} \bigr)\ind_{\{ \beta(u)<
t\}}\biggr)^2\biggr]
= A_t +B_t+C_t,
\]
where
\begin{eqnarray*}
A_t &=& \E[D_t]^{-2} \E_\mu\biggl[\sum_{u\in
\mathcal{T}}f^2\bigl(X^u_{\beta(u)-} \bigr)\ind_{\{ \beta(u)<
t\}} \biggr],\\
B_t&=& \E[D_t]^{-2} \E_\mu\biggl[\sum_{(u,v)\in\mathcal{FT}}
f\bigl(X^u_{\beta(u)-} \bigr)f\bigl(X^v_{\beta(v)-} \bigr)\ind_{\{
\beta(u)<
t, \beta(v)< t\}} \biggr],\\
C_t&=& 2\E[D_t]^{-2} \E_\mu\biggl[\sum_{u\prec v, v\in\ct}
f\bigl(X^u_{\beta(u)-} \bigr)f\bigl(X^v_{\beta(v)-} \bigr)\ind_{\{
\beta(v)< t\}} \biggr].
\end{eqnarray*}
The terms $A_t$ and $B_t$ will be handled similarly as in the proof of
Proposition~\ref{thLGN}.
%
%e52 ###
%
\begin{eqnarray}\label{carre2}
A_t &=& r \E[D_t]^{-2} \int_0^t ds\, \rme^{r(m-1)s} \E_\mu
[f^2(Y_{s-})]\nonumber\\[-8pt]\\[-8pt]
&=& \frac{r (m-1)^2}{ (\rme^{r(m-1)t} -1)^{2}} \int_0^t ds\, \rme
^{r(m-1)s} \mu Q_s f^2
\mathop{\longrightarrow}_{t\rightarrow\infty} 0\nonumber
\end{eqnarray}
thanks to (\ref{formuledualite2}) for the first
equality, (\ref{eq:EDt}) for the second and (H3) for the convergence.

Notice that Corollary \ref{cor:fourche2bis} and then (H1) and (H4)
imply that
\begin{eqnarray*}
&&\E_\mu\biggl[ \sum_{(u,v)\in
\mathcal{FT}}\bigl|f\bigl(X^u_{\beta(u)-} \bigr)\bigr| \bigl|f\bigl(X^v_{\beta(v)-}
\bigr)\bigr|\ind_{\{ \beta(u)<
t, \beta(v)< t\}} \biggr] \\
&&\qquad= r^3 \int_{[0,+\infty)^3 } \mu Q_sJ_2(Q_{s'}|f| \otimes
Q_{s''}|f|)
\rme^{r(m-1)(s+s'+s'')}\\
&&\qquad\quad\hphantom{r^3 \int_{[0,+\infty)}}\hspace*{1pt}
{}\times\ind_{\{s+s'<t,s+ s''<t\}} \,ds\,ds'\,ds''
\end{eqnarray*}
is finite and
\begin{eqnarray*}
B_t
& = &\frac{r^3 (m-1)^2}{ (\rme^{r(m-1)t} -1)^{2}}\\
&&{} \times
\int_{[0,+\infty)^3 } \mu Q_sJ_2(Q_{s'}f \otimes Q_{s''}f)
\rme^{r(m-1)(s+s'+s'')}\\
&&\hphantom{{}\times \int_{[0,+\infty)}}\hspace*{2pt}
{}\times\ind_{\{s+s'<t,s+ s''<t\}} \,ds\,ds'\,ds'' \\
& = &\frac{r^3 (m-1)^2\rme^{2r(m-1)t}}{ (\rme^{r(m-1)t} -1)^{2}}\int
_{[0,+\infty)^3 } ds \,dt' \,dt''\mu Q_sJ_2(Q_{t-t'}f \otimes Q_{t-t''}f)\\
&&\hspace*{129.7pt}{} \times
\rme^{r(m-1)(s-t'-t'')}\ind_{\{s<t'<t,s<t''<t\}}.
\end{eqnarray*}
Now, since $\langle\pi,f \rangle=0$, we deduce from (H2) that for
$t', t''$ fixed and
$y,z\in E$, $\lim_{t \rightarrow\infty}(Q_{t-t'}f \otimes
Q_{t-t''}f)(y,z)=0$. Thanks to (H1), there exists $g$ such that
$ |(Q_{t-t'}f \otimes Q_{t-t''}f)| \le(g \otimes
g)$ and (H4) implies that $\int_{[0,+\infty)^3 } \mu
Q_sJ_2(g \otimes g)\times
\rme^{r(m-1)(s-t'-t'')}\ind_{\{s<t',s<t''\}} \,ds\,dt'\,dt''$ is finite.
Then, by Lebesgue's theorem
%
%e53 ###
%
\begin{equation}
\label{eq:cvBt2}
\lim_{t \rightarrow\infty} B_t =0.
\end{equation}
Let us now consider $C_t$. We have $C_t\leq C'_t+ C''_t$ where
\begin{eqnarray*}
C'_t
& = & \E[D_t]^{-2} \E_\mu\biggl[\sum_{u\prec v, v\in\ct}
f^2\bigl(X^v_{\beta(v)-} \bigr)\ind_{\{ \beta(v)< t\}} \biggr],\\
C''_t & = & \E[D_t]^{-2} \E_\mu\biggl[\sum_{u\prec v, v\in\ct}
f^2\bigl(X^u_{\beta(u)-} \bigr)\ind_{\{ \beta(v)< t\}} \biggr].
\end{eqnarray*}
We deduce from (\ref{formuledualite2}) that
\begin{eqnarray*}
C'_t
&=&\E[D_t]^{-2} \E_\mu\biggl[\sum_{ v\in\ct}
|v|f^2\bigl(X^v_{\beta(v)-} \bigr)\ind_{\{ \beta(v)< t\}} \biggr]\\
&=& \E[D_t]^{-2}r\int_0^{t } ds\, \rme^{r(m-1)s} \E_\mu[S_{s-}
f^2(Y_{s-})]\\
&=& \frac{r (m-1)^2}{ (\rme^{r(m-1)t} -1)^{2}} \int_0^{t } ds\,
\rme^{r(m-1)s} \E_\mu[S_{s}
f^2(Y_{s})].
\end{eqnarray*}
We deduce from (H5) that
%
%e54 ###
%
\begin{equation}
\label{eq:cvC'}
\lim_{t \rightarrow\infty} C'_t =0.\
\end{equation}
Using the conditional expectation w.r.t. $X^u$, (\ref{eq:EDt}) and
(\ref{formuledualite2}), we get
\begin{eqnarray*}
C''_t
&=&\E[D_t]^{-2} \E_\mu\biggl[\sum_{u\in\ct}
f^2\bigl(X^u_{\beta(u)-} \bigr)\ind_{\{ \beta(u)< t\}} m\E\biggl[\sum
_{v\in
\ct} \ind_{\{\beta(v)<t'\}}\biggr]_{|t'=t-\beta(u)} \biggr]\\
&=& \frac{m}{m-1} \E[D_t]^{-2} \E_\mu\biggl[\sum_{u\in\ct}
f^2\bigl(X^u_{\beta(u)-} \bigr)\ind_{\{ \beta(u)< t\}} \bigl(\rme
^{r(m-1) (t-\beta(u))} -1\bigr) \biggr]\\
&=& \frac{m}{m-1} \E[D_t]^{-2} \int_0^t ds\, \rme^{r(m-1)s} \E_\mu
\bigl[
f^2(Y_{s-}) \bigl(\rme^{r(m-1) (t-s)} -1\bigr) \bigr]\\
&\leq&\frac{m(m-1)\rme^{r(m-1)t}}{ (\rme^{r(m-1)t} -1)^{2}} \int_0^t
ds\, \mu Q_s f^2.
\end{eqnarray*}
We deduce from (H3) [or (H5)] that
%
%e55 ###
%
\begin{equation}
\label{eq:cvC''}
\lim_{t \rightarrow\infty} C''_t =0.
\end{equation}
The proof of (\ref{eq:path-lln-l2}), when $\langle
\pi,f \rangle=0$, results from (\ref{carre2})--(\ref{eq:cvC'}) and (\ref{eq:cvC''}).

In the general case, we have
%
%e56 ###
%
\begin{eqnarray}
\label{othercases2}\hspace*{21pt}
&&\E[D_t]^{-1} \sum_{u\in
\mathcal{T}}f\bigl(X^u_{\beta(u)-} \bigr)\ind_{\{ \beta(u)<
t\}} - \langle\pi,f \rangle W \nonumber\\[-10pt]\\[-10pt]
&&\quad= \E[D_t]^{-1} \sum_{u\in
\mathcal{T}}\bigl(f\bigl(X^u_{\beta(u)-} \bigr)-\langle\pi, f
\rangle\bigr) \ind_{\{ \beta(u)<
t\}} + \langle\pi,f \rangle\biggl( \frac{D_t}{\E[D_t]}-W \biggr).
\nonumber
\end{eqnarray}
Notice that if $f$ and $\mu$ satisfy (H1)--(H5) then so do $f-\langle
\pi
,f \rangle$ and $\mu$. Thanks to the first part of the proof and to
Lemma \ref{lemmeDt}, we get
(\ref{eq:path-lln-l2}) if $f$ and $\mu$ satisfy (H1)--(H5). The convergence
in probability is thus obtained thanks to (\ref{eq:path-lln-l2}) and
(\ref{eq:CV-N-W}).\vspace*{-3pt}
\end{pf*}
\begin{pf*}{Proof of Theorem \ref{thLGNfonctionnel}}
The proof is similar to the proof of Theorem~\ref{thLGN}. Some
arguments are shorter as we assume that $\varphi$ is bounded.

We shall first consider the case $\E_\pi[\varphi(\w
{Y}_{[0,T]})]=0$. We
assume that $t>T$.
\[
\E[N_t]^{-2}\E_\mu\biggl[\biggl(\sum_{u\in
V_t}
\varphi\bigl(X^u_{[t-T, t]}, \Lambda^u_{[t-T, t]} -\Lambda
^u_{t-T}\bigr)
\biggr)^2\biggr]
= A_t +B'_t+B''_t,
\]
where
\begin{eqnarray*}
A_t& =&\E[N_t]^{-2} \E_\mu\biggl[\sum_{u\in V_t}\varphi^2\bigl(X^u_{[t-T,
t]}, \Lambda^u_{[t-T, t]} -\Lambda^u_{t-T}\bigr)
\biggr],\\[-2pt]
B'_t&=& \E[N_t]^{-2} \E_\mu\biggl[\mathop{\sum_{(u,v) \in V_t^2}}_{u
\ne v} \varphi\bigl(X^u_{[t-T, t]}, \Lambda^u_{[t-T, t]}
-\Lambda^u_{t-T}\bigr)\\[-2pt]
&&\hspace*{89.2pt}{} \times\varphi
\bigl(X^v_{[t-T, t]}, \Lambda^v_{[t-T, t]}
-\Lambda^v_{t-T}\bigr)
\ind_{\{\beta(u\wedge v) \geq t-T\}}\biggr],\\[-2pt]
B''_t&=& \E[N_t]^{-2} \E_\mu\biggl[\mathop{\sum_{(u,v) \in V_t^2}}_{u
\ne v} \varphi\bigl(X^u_{[t-T, t]}, \Lambda^u_{[t-T, t]}
-\Lambda^u_{t-T}\bigr)\\[-2pt]
&&\hspace*{89.1pt}{} \times\varphi
\bigl(X^v_{[t-T, t]}, \Lambda^v_{[t-T, t]}
-\Lambda^v_{t-T}\bigr) \ind_{\{\beta(u\wedge v) < t-T\}}\biggr].
\end{eqnarray*}
We assume that $\varphi$ is bounded by a constant, say $c$.
We have $A_t\leq c^2 \E[N_t]^{-1}$ so that $\lim_{t\rightarrow\infty}
A_t=0$. We have, using Proposition \ref{fourche2},
\begin{eqnarray*}
|B'_t|
&\leq& c^2 \E[N_t]^{-2} \E_\mu\biggl[\mathop{\sum_{(u,v) \in V_t^2}}_{
u \ne v}\ind_{\{\beta(u\wedge v) \geq t-T\}}\biggr]\\[-2pt]
&=&c^2 r \int_0^t \rme^{-r(m-1)a}
\ind_{\{a\geq t-T\}} \,da ,
\end{eqnarray*}
so that $\lim_{t\rightarrow\infty}
B'_t=0$.\vadjust{\goodbreak}

We set $f(x)=\E_x[\varphi(\w{Y}_{[0,T]})]$. Using Proposition \ref
{fourche2} once
more, we get
\begin{eqnarray*}
B''_t
&=& r \int_0^t \rme^{-r(m-1)a}\\[-2pt]
&&\hspace*{18.1pt}{}\times\E_\mu\bigl[
\w{J}_2\bigl({\E}_{\cdot}'\bigl[\varphi\bigl(Y'_{[t-a-T, t-a]},
\Lambda'_{[t-a-T, t-a]} -\Lambda'_{t-a-T} \bigr) \bigr]\\[-2pt]
&&\hspace*{61pt}{}\otimes
{\E}_{\cdot}'\bigl[\varphi\bigl(Y'_{[t-a-T, t-a]},
\Lambda'_{[t-a-T, t-a]} -\Lambda'_{t-a-T} \bigr) \bigr]
\bigr)(\w{Y}_{a})\bigr]\\[-2pt]
&&\hspace*{18.1pt}{}\times\ind_{\{a<t-T\}}\,da\\[-2pt]
&=& r \int_0^{t-T} \rme^{-r(m-1)a} \mu Q_a
J_2(Q_{t-a-T} f \otimes Q_{t-a-T}f
)\,da.
\end{eqnarray*}
By hypothesis on $Y$, we have
that, for fixed $a$, $\lim_{t\rightarrow\infty}
Q_{t-a-T}f=\langle\pi, f\rangle=0$. Using Lebesgue's theorem, we get
$\lim_{t\rightarrow\infty} B''_t=0$.
This gives the result for the $L^2(\PP_\mu)$ convergence when
$\langle
\pi, f \rangle=0$. We conclude in the general case and for the
convergence in probability
as in the proof of Theorem \ref{thLGN}.
\end{pf*}
\begin{pf*}{Proof of Proposition \ref{prop:cvL2ssarbre-path}}
The proof is similar to the proof of Proposition~\ref{cvL2ssarbre}. Some
arguments are shorter as we assume that $\varphi$ is bounded.

We shall first prove (\ref{eq:all-path-lln-l2}) for $\varphi$ such
that $\langle\pi, J_1 f \rangle=0$.
We have
\[
\E[D_t]^{-2}\E_\mu\biggl[\biggl(\sum_{u\in
\mathcal{T}}
\varphi\bigl(X^u_{[\alpha(u), \beta(u))} \bigr)
\ind_{\{ \beta(u)< t\}}\biggr)^2\biggr]
= A_t +B_t+C_t,
\]
where
\begin{eqnarray*}
A_t&=&\E[D_t]^{-2} \E_\mu\biggl[\sum_{u\in
\mathcal{T}}\varphi\bigl(X^u_{[\alpha(u), \beta(u))} \bigr)
^2\ind_{\{ \beta(u)<
t\}} \biggr],\\[-2pt]
B_t&=& \E[D_t]^{-2} \E_\mu\biggl[\sum_{(u,v)\in\mathcal{FT}}
\varphi\bigl(X^u_{[\alpha(u), \beta(u))} \bigr)
\varphi\bigl(X^v_{[\alpha(v), \beta(v))} \bigr)
\ind_{\{ \beta(u)<
t, \beta(v)< t\}} \biggr],\\[-2pt]
C_t&=& 2\E[D_t]^{-2} \E_\mu\biggl[\sum_{u\prec v, v\in\ct}
\varphi\bigl(X^u_{[\alpha(u), \beta(u))} \bigr)
\varphi\bigl(X^v_{[\alpha(v), \beta(v))} \bigr)
\ind_{\{
\beta(v)< t\}} \biggr].
\end{eqnarray*}
We assume that $\varphi$ is bounded by a constant, say $c$.
We have $A_t\leq c^2 /\E[D_t]$ so that $\lim_{t\rightarrow\infty}
A_t=0$. Thanks to Corollary \ref{cor:fourche2bis}, we have
\begin{eqnarray*}
|C_t|
&\leq& 2 c^2\E[D_t]^{-2} \E_\mu\biggl[\sum_{v\in\ct} |v| \ind_{\{
\beta(v)< t\}} \biggr] \\[-2pt]
&=& 2 c^2\E[D_t]^{-2} \int_0^t ds\, \rme^{r(m-1)s}
\E_\mu[S_s] \\[-2pt]
&=& 2 c^2\E[D_t]^{-2} \int_0^t ds\, srm \rme^{r(m-1)s}.
\end{eqnarray*}
This implies that\vadjust{\goodbreak} $\lim_{t\rightarrow\infty}
C_t=0$.

We set $h_t(x)=\E_x[\varphi(X_{[0,\tau)}) \ind_{\{\tau<t\}}]$, where
$\tau$ is an exponential random variable with mean 1, independent of
$X$.

Using the conditional expectation w.r.t. $X^{u'}$, where $u'$ is the
ancestor of~$u$, and $X^{v'}$, where
$v'$ is the ancestor of $v$, we have, according to $u'=v'$ or $u'\neq v'$,
\[
B_t= B'_t+B''_t,
\]
where
\begin{eqnarray*}
B'_t
&=&\E[D_t]^{-2}\E_\mu\biggl[\sum_{u'\in\ct}
J_2\bigl(h_{t-\beta(u')}\otimes h_{t-\beta(u')}\bigr)\bigl(X^{u'}_{\beta(u')-}\bigr)
\ind_{\{ \beta(u')<
t\}} \biggr], \\
B''_t
&=&\E[D_t]^{-2}\E_\mu\biggl[\sum_{(u',v')\in\mathcal{FT}}
J_1\bigl(h_{t-\beta(u')}\bigr) \bigl(X^{u'}_{\beta(u')-}\bigr)\\
&&\hphantom{\E[D_t]^{-2}\E_\mu\biggl[\sum_{(u',v')\in\mathcal{FT}}}
{} \times
J_1\bigl(h_{t-\beta(v')}\bigr) \bigl(X^{v'}_{\beta(v')-}\bigr)
\ind_{\{ \beta(u')<
t, \beta(v')< t\}} \biggr].
\end{eqnarray*}
Using the definition of $J_2$, (\ref{eq:J2-0}), we get $|B'_t|\leq c^2
\E[D_t]^{-1} (\varsigma^2+m^2-m)$ and thus, $\lim_{t\rightarrow
\infty}
B'_t=0$.

We deduce from Corollary \ref{cor:fourche2bis}, that
\begin{eqnarray*}
B''_t
&=& \frac{r^3 (m-1)^2}{ (\rme^{r(m-1)t} -1)^{2}}\\
&&{}\times\int_{[0,+\infty)^3 } ds\,ds'\,ds''
\mu Q_sJ_2(Q_{s'}J_1h_{t-s-s'} \otimes
Q_{s''}J_1h_{t-s-s''} )\\
&&\hphantom{{}\times\int_{[0,+\infty)^3 }}
{}\times\rme^{r(m-1)(s+s'+s'')}\ind_{\{s+s'<t,s+ s''<t\}}
\\
&=& \frac{r^3 (m-1)^2\rme^{2 r(m-1)t}}{ (\rme^{r(m-1)t} -1)^{2}}\\
&&{}\times\int_{[0,+\infty)^3 } ds\,dv'\,dv''
\mu Q_sJ_2(Q_{t-s-v'}J_1h_{v'} \otimes
Q_{t-s-v''}J_1h_{v''} )\\
&&\hphantom{{}\times\int_{[0,+\infty)^3 }}
{}\times\rme^{-r(m-1)(s+v'+v'')}\ind_{\{v'<t-s, v''<t-s\}}.
\end{eqnarray*}
By hypothesis on $Y$, we have
that, for fixed $s$ and $v$, $\lim_{t\rightarrow\infty}
Q_{t-s-v}J_1h_{v}=\langle
\pi, J_1 h_{v} \rangle$. Using Lebesgue's theorem, we get
\[
\lim_{t\rightarrow\infty} B''_t
= r^3(m-1)^2 \int_{[0,+\infty)^3 } ds\,dv'\,dv'' \langle
\pi, J_1 h_{v'} \rangle\langle
\pi, J_1 h_{v''} \rangle\rme^{-r(m-1)(s+v'+v'')}.
\]
Notice that $h_t(x)=\E_x[\varphi(Y_{[0,
\tau_1)})\rme^{r(m-1)\tau_1} \ind_{\{\tau_1<t\}}]/m$ so that
\[
% \label{eq:h-X-Y}
r(m-1) \int_0^{+\infty} dt \,h_t(x)\rme^{-r(m-1)t}
=\frac{1}{m}\E_x\bigl[\varphi\bigl(Y_{[0, \tau_1)}\bigr)\bigr].
\]
Recall $f(x)=\E_x[\varphi(Y_{[0, \tau_1)})]$. We get $\lim
_{t\rightarrow
\infty} B''_t =\frac{1}{(m-1)m^2} \langle\pi, J_1 f \rangle^2=0$.
The\-refore, we get that
\[
\lim_{t\rightarrow
\infty} A_t+B_t+C_t=0,
\]
which gives the result for the $L^2(\PP_\mu)$ convergence when
$\langle
\pi, J_1 f \rangle=0$. We conclude in the general case and for the
convergence in probability
as in the proof of Proposition \ref{cvL2ssarbre}.
\end{pf*}

%s5 ###
\section{Examples} \label{sectionexamples}
We now investigate several examples. In Section
\ref{sectionexspitteddiffusion}, splitted diffusions are considered as
scholar examples. In Section \ref{exaging}, we give a biological
application to ``cellular aging'' when cells divide in continuous time,
which is one of the motivations of this work. In Section
\ref{exkangourous}, we give a~central limit theorem for nonlocal
branching L\'{e}vy processes.

%s5.1 ###
\subsection{Splitted real diffusions}\label{sectionexspitteddiffusion}

A first example consists in binary branching:
the continuous tree $\T$ is a Yule tree. For the Markov process $X$, we
consider a real diffusion with generator
%
%e57 ###
%
\begin{equation}
\label{eq:defLexple}
Lf(x)=b(x)f'(x)+\frac{\sigma^2(x)}{2}f''(x).
\end{equation}
We assume that $b$ and $\sigma$ are such that there exists a unique
strong solution to the corresponding SDE (see,
e.g., \cite{ikedawatanabe}, Theorem 3.2, page 182).

When a branching occurs, each daughter inherits a random fraction of the
value of the mother
\[
F^{(1)}(x,\theta)= G^{-1}(\theta)x,\qquad F^{(2)}(x,\theta)=
\bigl(1-G^{-1}(\theta)\bigr)x ,
\]
where $G$ is the cumulative distribution function of the random
fraction in $[0,1]$ associated with the branching event. We assume the
distribution of the random fraction is symmetric: $G(x)=1-G(1-x)$.

The infinitesimal generator of $Y$ is characterized for $f\in
\Co^2_b(\R,\R)$ by
%
%e58 ###
%
\begin{eqnarray}
Af(x)& = & b(x)f'(x)+\sigma(x)f''(x)\nonumber\\
&&{} +  2r \int_0^1 \biggl(\frac{1}{2}\bigl(f(G^{-1}(\theta)x)-f(x)
\bigr)\nonumber\\[-8pt]\\[-8pt]
&&\hspace*{45pt}{} +\frac{1}{2}\bigl(f\bigl(\bigl(1-G^{-1}(\theta)\bigr)x\bigr)-f(x)\bigr)\biggr)\,d\theta
\nonumber\\
& = & b(x)f'(x)+\sigma(x)f''(x)+2r \int_0^1
\bigl(f(qx)-f(x)\bigr)G(dq).\nonumber
\end{eqnarray}
Particular choices for the functions $b$ and $\sigma$ are the
following ones:
\begin{longlist}
\item If $b(x)=0$ and $\sigma(x)=\sigma$, we obtain the
splitted Brownian process.
\item If $b(x)=-\beta(x-\alpha)$ and $\sigma(x)=\sigma$, we obtain
the splitted Ornstein--Uhlenbeck process.
\item If $b(x)=1$ and $\sigma(x)=0$, the deterministic process
$X$ can represent the linear growth of some biological content of the
cell (nutriments, proteins, parasites$\ldots$) which is shared randomly in
the two daughter cells when the cell divides. More precisely here, each
daughter inherits
random fraction of this biological content.
% Moreover, the auxiliary process which then grows linearly and
% undergoes multiplication by idependant random variables with the same
% cumulative distribution function $G$ is classical in TCP (??????????)
% theory (e.g., \cite{dgr, grd, jac}).
\end{longlist}
Let us note that if $b(x)=\beta x$ and $\sigma(x)^2=\sigma^2 x$, we obtain
the splitted Feller branching diffusion. But in this case, almost
surely, the auxiliary process either
becomes extinct or goes to infinity as $t\rightarrow
\infty$. The assumption (H2) is not satisfied. This process is studied
in \cite{bansayetran}
as a model for parasite infection.
%in a more general setting as the branching rate $r$ may depend on the
%state of the cell.

The following results give the asymptotic limit of the splitted
diffusion under some condition which is satisfied by the
examples (i)--(iii). For this we use results due to Meyn and Tweedie
\cite{meyntweedieII,meyntweedie}.
\begin{prop}\label{propapplimeyntweedie}
Assume that $Y$ is Feller and irreducible (see \cite{meyntweedie},
pa\-ge~520) and that there exists $K\in\R_+$, such that for every
$|x|>K$, $b(x)/x<r'$ with $r'<r$.
Then the auxiliary process $Y$ with generator $A$ is
ergodic with stationary probability $\pi$. Furthermore, $\sum_{u\in
V_t}\delta_{X^u_t}(dx)/N_t$ converges weakly to $\pi$ as
$t\rightarrow\infty$ and this convergence holds in probability.
\end{prop}

\begin{pf}
Once we check that $Y$ is ergodic, then Corollary \ref{corolintro} and
the fact that
$W$ defined by (\ref{eq:CV-N-W}) is a.s. positive readily imply the
weak convergence of the
proposition.
To prove the ergodicity of $Y$, we use Theorem 4.1 of \cite
{meyntweedieII} and Theorem~6.1 of \cite{meyntweedie}. Since $Y$ is Feller and
irreducible, the process $Y$ admits a~unique invariant probability
measure $\pi$ and is exponentially ergodic provided the condition
(CD3) in \cite{meyntweedie} is satisfied. Namely, if there exists a~%
positive measurable function $V \dvtx x\mapsto V(x)$ such that $\lim
_{x\rightarrow\pm\infty}V(x)=+\infty$ and for which
%
%e59 ###
%
\begin{equation}\label{conditiondeLyapounov}
\exists c>0, d\in\R, \forall x\in\R\qquad AV(x)\leq
-cV(x)+d.
\end{equation}
For $V(x)=|x|$ regularized on an $\varepsilon$-neighborhood of 0
($0<\varepsilon<1$), we have
%
%e60 ###
%
\begin{eqnarray}\label{applifosterlyapunov}
\forall|x|>\varepsilon\qquad AV(x)&=& \operatorname{sign}(x) b(x)+2r |x| \int
_0^1 (q-1)G(dq)\nonumber\\[-8pt]\\[-8pt]
&=&\operatorname{sign}(x)b(x)-r |x|\nonumber
\end{eqnarray}
as the distribution of $G$ is symmetric. By assumption, there exists
$\eta>0$ and $K>\varepsilon$, such that (\ref{applifosterlyapunov}) implies
%
%e61 ###
%
\begin{equation}\label{critereLyapounovexemple}
\forall x\in\R\qquad AV(x)\leq-\eta V(x)+\Bigl({\sup_{|x|\leq
K}}|b(x)|+rK\Bigr)\ind_{\{|x|\leq K\}}.
\end{equation}
This implies (\ref{conditiondeLyapounov}) and finishes the proof; the
geometric ergodicity expresses here as
%
%e62 ###
%
\begin{eqnarray}\label{exponentielltgeom1}
&&\exists\beta>0, B<+\infty, \forall t\in\R_+, \forall x\in\R
\nonumber\\[-10pt]\\[-10pt]
&&\qquad{\sup_{g / |g(u)|\leq
1+|u|}}|Q_t g(x)-\langle\pi,g\rangle|\leq B(1+|x|)\rme
^{-\beta t}.
\nonumber\vspace*{-2pt}
\end{eqnarray}
\upqed\end{pf}
\begin{rque}\label{rquemeyntweedie}The examples (i)--(iii) satisfy the
assumptions of Proposition~\ref{propapplimeyntweedie}. If $b$ and
$\sigma$ are bounded Lipschitz functions, $X$ is Feller (e.g., \cite
{stroockvaradhan}, Theorem 6.3.4, page 152) and thus, $Y$ is also
Feller. The Feller property also holds for Ornstein--Uhlenbeck
processes. The irreducibility property is well known for diffusions as
(i) and (ii) and trivial for (iii).\vspace*{-3pt}
\end{rque}
\begin{rque}\label{rquemeyntweedie2}
If there exists $K>0$ in Proposition \ref{propapplimeyntweedie} such
that for every $|x|>K$, $2 b(x)/x+6 \sigma(x)/x^2<r'$ with $r'<r\int
_0^1 (1-q^4)^2 G(dq)$, then we can use similar arguments as in the
proof of Proposition \ref{propapplimeyntweedie}. We get that the
auxiliary process $Y$ is geometrically ergodic with
%
%e63 ###
%
\begin{eqnarray}\label{exponentielltgeom2}
&&\exists\beta>0, B<+\infty, \forall t\in\R_+, \forall x\in\R
\nonumber\\[-10pt]\\[-10pt]
&&\qquad\sup_{g / |g(u)|\leq
1+|u|^4}|Q_t g(x)-\langle\pi,g\rangle|\leq B(1+|x|^4)\rme
^{-\beta t}
\nonumber
\end{eqnarray}
instead of (\ref{exponentielltgeom1}). This will be used for the proof
of the central limit theorem.\vspace*{-3pt}
\end{rque}

%s5.2 ###
\subsection{Cellular aging process}\label{exaging} We now
present a generalization to the continuous time of Guyon \cite{guyon} and
Delmas and Marsalle \cite{delmasmarsalle} about cellular aging. When a
rod shaped cell divides, it
produces a new end per progeny cell. So each new cell has a pole (or
end) which is new and another one which was created one or more
generations ago. This number of generations is the age of the cell.
Since each
cell has a new pole and an older one, at the next division one of the
two daughters will inherit the new pole and
the other one will inherit the older pole. Experiments indicate that the
first one has a~larger growth rate than the second one (see Stewart
et al. \cite{tad} for details), which indicates aging.

To detect this aging effect, \cite{guyon,delmasmarsalle} used discrete
time Markov models by looking at cells of a given generation.
Considering continuous time genealogical trees may allow to take into
account the asynchrony of cell divisions.

We consider the following model. Cells are characterized by a type
$\eta\in\{0,1\}$ (type 0 corresponds to a cell of age 1 and type 1 to
cell of greater age) and a quantity $\zeta$ (growth rate, quantity of
damage in the cell)
that evolves according to a
Markov process depending on the type of
the cell. Cells may die, which leads us to the following model. At rate
$r$, each cell is replaced
by one cell of type 0 (resp., $1$) with
probability $p_0\geq0$ (resp.,
$p_1\geq0$), by two cells of type 0 and 1 with probability
$p_{0,1}\geq0$, or by no cell with probability $1-p_0-p_1-p_{0,1}\geq
0$. The way the quantity $\zeta$ is given to a daughter depends on its
type and on the fact that it has or has not a sister.\vadjust{\goodbreak}

This can be stated in the framework of Sections
\ref{sectiondescriptiontreeindexedmarkovprocess} and \ref{sec:MT1}.
For the sake of simplicity, we shall assume that $\zeta$ evolves as
a real
diffusion between two branching times.

Let $L^0$ and $L^1$ be two diffusion generators; for $f\in
\cc^2(\R\times\{0,1\}, \R)$,
%
%e64 ###
%
\begin{equation}\label{exempledelmasmarsalle0}
L^\eta f(\zeta,\eta)=b(\zeta,\eta)\,\partial_\zeta
f(\zeta,\eta)+\sigma(\zeta,\eta) \,\partial^2_{\zeta,\zeta}f
(\zeta,\eta), \qquad\eta\in\{0,1\}.
\end{equation}
We assume there exists a unique strong solution to the corresponding two
SDEs (see, e.g., \cite{ikedawatanabe}, Theorem 3.2, page 182). We
consider the underlying process $X=((\zeta_t,\eta_t) , t\geq0)$ with
generator
\[
Lf(\zeta,\eta)=\ind_{\{\eta=0\}}L^0 f(\zeta,0)+\ind_{\{\eta=1\}
}L^1 f(\zeta,1).
\]
Notice the process $(\eta_t, t\geq0)$ is constant between two
branching times.
The offspring distribution is
%
%e65 ###
%
\begin{equation}
\label{exempledelmasmarsalle1}\quad
p(dk)=(1-p_0-p_1-p_{0,1}) \delta_0(dk) + (p_0+p_1) \delta_1(dk) +
p_{0,1} \delta_2(dk).
\end{equation}
The offspring position is given by
%
%e66 ###
%
\begin{eqnarray}\label{exempledelmasmarsalle2}
F_{1}^{(1)}((\zeta,\eta),\theta) & = & (g_0(\zeta,
\eta,\theta_1),0)\ind_{\{\theta_2\leq
p_0/(p_0+p_1)\} }\nonumber\\
&&{}+(g_1(\zeta,\eta,\theta_1),1)\ind_{\{\theta_2>
p_0/(p_0+p_1)\}},\\
F_{i+1}^{(2)}((\zeta,\eta),\theta) & = &
\bigl(g_i^{(2)}(\zeta,\eta,\theta),i\bigr) \qquad\mbox{for $i\in\{0,1\}$,}
\nonumber
\end{eqnarray}
for some functions $g_0, g_1, g_0^{(2)}, g_1^{(2)}$ and $(\theta
_1,\theta_2)$ a
function of $\theta$ such that if $\theta$ is uniform on $[0,1]$, then
$\theta_1$ and $\theta_2$ are independent and uniform on $[0,1]$. The
division is asymmetric if $g_0^{(2)}\neq g_1^{(2)}$. One important issue
is, using the LLN (Section~\ref{sectionlgn}) and
fluctuation results, to test if the division is asymmetric, which means
aging, or not. Let us mention that a natural question would
be to give the test in a more general model in which the division rate
depends on the state of the cell and of the quantity of interest $\zeta$
(which is realistic if, e.g., $\zeta$ describes the quantity of damage
of the
cell).

Let us consider a test function $f \dvtx (t,\zeta,\eta)\mapsto
f_t(\zeta,\eta)$ in $\Co^{1,2}_b(\R_+\times(\R\times\{0,1\}),\R
)$ and
let $(B^u)_{u\in\mathcal{U}}$ be a family of independent standard
Brownian motions. The SDE describing the evolution of the population of
cells then becomes, with the notation of (\ref{exempledelmasmarsalle0}),
(\ref{exempledelmasmarsalle1}) and (\ref{exempledelmasmarsalle2}),
\begin{eqnarray*}
\langle Z_t,f_t\rangle &=& \langle Z_0,f_0\rangle\\[-2pt]
&&{} +\int_0^t \int_{\mathcal{U}\times\{0, 1,2\}\times[0,1]}
\mathbf{1}_{\{u\in
V_{s_-}\} }\\[-2pt]
&&\hspace*{92.5pt}{}\times\Biggl(\sum_{j=1}^{k}
f_s\bigl(F_j^{(k)}((\zeta^u_{s_-},\eta_{s_-}^u),\theta)\bigr)
-f_s(\zeta^u_{s_-},\eta_{s_-}^u)\Biggr)\\[-2pt]
&&\hspace*{92.5pt}{}\times\rho(ds,du,dk,d\theta)\\[-2pt]
&&{} + \int_0^t\int_{\R\times\{0,1\}}\bigl(L^\eta
f_s(\zeta,\eta)+\partial_s
f_s(\zeta,\eta)\bigr)Z_s(d\zeta,d\eta) \,ds \\[-2pt]
&&{} +\int_0^t \sum_{u\in V_{s}}
\sqrt{2\sigma(\zeta^u_{s},\eta^u_{s})}\,\partial_{\zeta}f
(\zeta^u_s,\eta^u_s) \,dB^u_s.
\end{eqnarray*}

If $Y$ is ergodic, then $\ind_{\{N_t>0\}}\sum_{u\in V_t}\delta
_{X^u_t}(dx)/N_t$
converges to a deterministic nondegenerated measure on $\R_+\times\{
0,1\}$.
Given a particular choice for the parameters $g_0$, $g_1$, $g_0^{(2)}$,
$g_1^{(2)}$, $L^0$, $L^1$, $p_0$, $p_1$ and $p_{0,1}$ of the model, one
can use arguments similar to the ones used in Proposition \ref
{propapplimeyntweedie} and Remark~\ref{rquemeyntweedie} to prove the
ergodicity of $Y$. Let us give an example that can be viewed as a
direct generalization of the model of Delmas and Marsalle
\cite{delmasmarsalle}.

The quantity $\zeta\in\R$ models the cell growth rate which is
assumed constant during the cell's life: $b(\zeta,\eta)=0$ and
$\sigma(\zeta,\eta)=0$.
For the functions $g_0$, $g_1$, $g_0^{(2)}$, $g_1^{(2)}$ which describe
the daugthers' growth rates, as functions of their mothers'
characteristics, we set
%
%e67 ###
%
\begin{eqnarray}\label{exemple:growthrate}
g_0(\zeta,\eta,\theta)&=&\alpha_0 \zeta+\beta_0+ \varepsilon_0,\qquad
g_1(\zeta,\eta,\theta)=\alpha_1 \zeta+\beta_1+
\varepsilon_1, \nonumber\\[-10pt]\\[-10pt]
g_0^{(2)}(\zeta,\eta,\theta)&=&\alpha'_0 \zeta+\beta'_0+
\varepsilon'_0,\qquad g_1^{(2)}(\zeta,\eta,\theta)= \alpha'_1
\zeta+\beta'_1+ \varepsilon'_1,\nonumber
\end{eqnarray}
where $\alpha_0$, $\alpha_1$, $\alpha'_0$, $\alpha'_1\in(-1,1)$
and where $\beta_0$, $\beta_1$, $\beta'_0$, $\beta'_1\in\R$. The
random variables $\varepsilon_0$, $\varepsilon_1$, $\varepsilon'_0$
and $\varepsilon'_1$ generated thanks to the uniform variable $\theta
$ and their distributions are as follows: $\varepsilon'_0$ and
$\varepsilon'_1$ are Gaussian centered r.v. with variances $\sigma
_0^2>0$ and $\sigma_1^2>0$, respectively, while $(\varepsilon
_0,\varepsilon_1)$ is a vector of Gaussian centered r.v. with covariance
\[
\sigma^2 \pmatrix{
1 & \rho\cr
\rho& 1}, \qquad\sigma^2>0,\qquad \rho\in(-1,1).
\]
In \cite{delmasmarsalle}, this model is used to test aging phenomena,
for instance, which correspond to $(\alpha_0,\beta_0)=(\alpha
_1,\beta_1)$. Delmas and Marsalle in discrete time prove that the
auxiliary process, which correspond to the Markov chain associated here
to the continuous time pure jump process $Y$, is ergodic. As a
consequence, $Y$ is recurrent, admits an invariant probability
distribution [since the jump rate $r(p_0+p_{1}+2p_{0,1})$ is a
constant] and is hence ergodic (see, e.g., Norris \cite{norris}).

%we can take $b(\zeta,\eta)=-b_\eta\ind_{\zeta>0}$ and $\sigma(\zeta,
%decreases linearly with time, and faster for cells of type 1. The
%constant $\sigma$ stands for a noise due to the environment. For the
%functions $g_0$, $g_1$, $g_0^{(2)}$, $g_1^{(2)}$ which model the
%growth rates at birth of the daughters, as functions of their mothers'
%characteristics, we set for a positive $\varepsilon>0$:
%& g_0(\zeta,\eta,\theta)=g_0^{(2)}(\zeta,\eta,\theta)=\zeta-
%& g_1(\zeta,\eta,\theta)=g_1^{(2)}(\zeta,\eta,\theta)=\zeta+
%example, the growth rates of new offspring depend only on their types,
%and not on the number of children left by the mother. Since the
%division rate is $r$, the mean time between two successive divisions
%is $1/r$. The average variation of the growth rate on such interval is
%$-b_{\eta}/r$. Thus, $\zeta+b_\eta/r$ is an approximation of the
%mother's growth rate at birth, with a random penalty for cells of type
%$0$ and a random bonus for cells of type 1.

%s5.3 ###
\vspace*{-4pt}\subsection{\texorpdfstring{Branching L\'{e}vy process}{Branching
Levy process}} \label{exkangourous}

We consider particles moving independent\-ly on $\R$ following a L\'{e}vy
process $X$ and reproducing with constant rate $r$. Each child jumps
from the location of the mother when the branching occurs. We are
interested in the rescaled population location at large time.

The generator of the underlying process $X$ is given by
\begin{eqnarray*}
Lf(x)&=& b
f'(x)+\frac{\sigma^2}{2}f''(x)\\[-2pt]
&&{}+\int_{\R\setminus\{0\}}
\bigl(f(x+y)-f(x)-yf'(x)\ind_{\{ |y|<1\}}\bigr)h(dy)\vadjust{\goodbreak}
\end{eqnarray*}
with\vspace*{1pt} $b\in\R$, $\sigma\in\R_+$ and $h$ a measure on $\R\setminus
\{0\}$
such that $\int_{\R\setminus\{0\}}y^2 h(dy)<+\infty$. The
particles reproduce at rate $r$ in
a random number of offspring distributed as $p=(p_k, k\in\N)$ such
that $\sum_{k\geq1} kp_k>1$ (supercritical case).
The offspring position is defined as
%
%e68 ###
%
\begin{equation}
F_j^{(k)}(x,\theta)=x+\Delta_j^k(\theta),\qquad j\in\{ 1,\ldots,k\},
\end{equation}
where we recall that $x$ is the location just before branching time and
$k$ is the number of offspring.
We assume\vspace*{1pt} the following second moment condition:
$\sum_{k\in\N} p_k \sum_{j=1}^k \E[\Delta_j^k(\Theta)^2]<\infty$,
where $\Theta$ is uniform on $[0,1]$.
\begin{prop}
We have the following weak convergence in $\mathcal{M}_F(\R)$:
%
%e69 ###
%
\begin{equation}
\label{CLTlev}
\lim_{t\rightarrow+\infty}\frac{1}{N_t}\sum_{u\in
V_t}\delta_{({X_t^u-\beta
t})/{\sqrt{t}}}(dx)=\pi_\Sigma(dx)\ind_{\{W>0\}}
\qquad\mbox{in probability,}
\end{equation}
where $\pi_\Sigma$ is the centered Gaussian probability measure with
variance $\Sigma$ and
%
%e71 ###
%e70 ###
%
\begin{eqnarray}
\beta &=& b + \int_{\R\setminus\{0\}}y\ind_{\{|y|\geq1\}}
h(dy)+r\sum_{k=1}^{+\infty}p_k\sum_{j=1}^k \E[\Delta_j^k(\Theta
)],\\
\Sigma &=& \sigma^2 + \int_{\R\setminus\{0\}}y^2 h(dy)
+r\sum_{k=1}^{+\infty}p_k\sum_{j=1}^k \E[\Delta_j^k(\Theta)^2].
\end{eqnarray}
\end{prop}
\begin{pf}
The auxiliary process $Y$ is a L\'{e}vy process with generator
\begin{eqnarray*}
Af(x) &=& b f'(x)
+\frac{\sigma^2}{2}f''(x)\\
&&{}+\int_{\R}\bigl(f(x+y)-f(x)-yf'(x)
\ind_{\{|y|<1\}}\bigr) h(dy)\\
&&{} + rm\sum_{k=1}^{+\infty}\frac{kp_k}{m}\int_0^1 \sum_{j=1}^k\frac{1}{k}
\bigl(f\bigl(x+\Delta_j^{k}(\theta)\bigr)-f(x)\bigr) \,d\theta.
\end{eqnarray*}
In particular, we have for all $x \in\R$,
\begin{eqnarray*}
\E_x[Y_t] &=& x +t\Biggl(b + \int_{\R\setminus\{0\}}y\ind_{\{|y|\geq
1\}} h(dy)+r\sum_{k=1}^{+\infty}p_k\sum_{j=1}^k
\E[\Delta_j^k(\Theta)]\Biggr)\\
&=&x+\beta t,\\
\E_x[Y_t^2]-\E_x[Y_t]^2 &=& t \Biggl(\sigma^2 +
\int_{\R\setminus\{0\}}y^2 h(dy)
+r\sum_{k=1}^{+\infty}p_k\sum_{j=1}^k
\E[\Delta_j^k(\Theta)^2]\Biggr)\\
&=&\Sigma t.
\end{eqnarray*}
Then, we deduce from the central limit theorem for L\'{e}vy processes
or directly from L\'{e}vy Khintchine formula, that $((Y_t-\beta
t)/\sqrt{t}, t\geq0)$ converges in distribution to $\pi_\Sigma$.
This implies that for any fixed $s$, $((Y_{t-s}-\beta
t)/\sqrt{t}, t\geq0)$ converges in distribution to $\pi_\Sigma$.

Let $\varphi$ be a continuous bounded real function and define
\[
f_t(x):=\varphi\bigl((x-\beta t)/\sqrt{t}\bigr) \qquad\mbox{for
$t\geq0, x
\in\R$.}
\]
Let $(Q_t, t\geq0)$ be the transition semigroup of $Y$. We get that
for any fixed~$s$ and $x\in\R$,
%
%e72 ###
%
\begin{equation}
\label{eq:cvftpif}
\lim_{t\rightarrow+\infty} Q_{t-s} f_t(x)= \langle
\pi_\Sigma,\varphi\rangle.
\end{equation}
It is then very easy to adapt the proof of
Theorem \ref{thLGN} with $f$ replaced by $f_t- \langle\pi_\Sigma,
\varphi\rangle$; (\ref{carre}) holds since
$f_t$ is uniformly bounded; (\ref{eq:Btcv}) holds using similar arguments
with (\ref{eq:cvftpif}) instead of (H2)
and $f_t$ uniformly
bounded instead of (H1) and~(H4). Similar arguments, as in the end of the proof of
Theorem~\ref{thLGN}, imply that for any continuous bounded real function
$\varphi$, the following convergence in probability holds:
\begin{eqnarray*}
\lim_{t\rightarrow+\infty}\frac{1}{N_t}\sum_{u\in
V_t} \varphi\biggl(\frac{X_t^u-\beta
t}{\sqrt{t}}\biggr)&=&\lim_{t\rightarrow+\infty}\frac{1}{N_t}\langle Z_t,
f_t \rangle\\
&=& \langle\pi_\Sigma, \varphi\rangle \ind_{\{W>0\}}.
\end{eqnarray*}
This gives (\ref{CLTlev}).
\end{pf}

%s6 ###
\section{Central limit theorem}\label{sec6}

%s6.1 ###
\subsection{Fluctuation process}\label{sec6.1}

In order to study the fluctuations associated to the LLNs, Theorem
\ref{thLGN}, we shall use the martingale associated to $Z_t$ [see~%
(\ref{martingalegdepop})]. We focus on the simple case of splitted
diffusions developed in Section~%
\ref{sectionexspitteddiffusion}. Our
main result for this section is stated as Proposition
\ref{propconvergencefluctuations}.

In the sequel, $C$ denotes a constant that may change from line to
line. We work in the framework of Section \ref
{sectionexspitteddiffusion}.

We consider the following family, indexed by $T > 0$, of fluctuation
processes. For $f\in\mathcal{B}_b(\R_+,\R)$ and $t\geq0$
%
%e73 ###
%
\begin{equation}\label{fluctuationsdef}
\langle\eta^T_t,f\rangle= \sqrt{\E[N_{t+T}]}\biggl(\frac{\langle
Z_{t+T},f\rangle}{\E[N_{t+T}]}-\frac{\langle Z_T,Q_t f\rangle}{\E
[N_T]} \biggr),
\end{equation}
where we recall that $N_t=\mbox{Card}(V_t)=\langle Z_t,1\rangle$ and
$Q_t$ has been defined in~(\ref{eq:Qt}). The family $Q_t$ is the
transition semigroup of the auxiliary process $Y$ which is given by
%
%e74 ###
%
\begin{equation}\label{processusauxilieredernierepartie}
Y_t=X_0+\int_0^t b(Y_s)\,ds+\int_0^t \sigma(Y_s)\,dB_s- \int_0^t
(1-q)Y_{s_-} \rho(ds,dq),
\end{equation}
where $X_0$ is an initial condition with distribution $\mu$, where
$(B_t , t\geq 0)$ is a standard real Brownian\vadjust{\goodbreak} motion and where $\rho
(ds,dq)$ is a Poisson point measure with intensity $2r \,ds\otimes
\widetilde{G}(dq)$ with $\widetilde{G}$ such that\vspace*{1pt} $\int_{[0,1]}
f(q)\widetilde{G}(dq)=\int_{[0,1]}(f(q)/2+f(1-q)/2) G(dq)$. As in
Section \ref{sectionexspitteddiffusion}, we will assume in the sequel
that $G$ is symmetric. In this case, $\widetilde{G}(dq)=G(dq)$.

The idea in (\ref{fluctuationsdef}) is to compare the independent
trees that have grown from the particles of $Z_T$ between times $T$ and
$t+T$, with the positions of independent auxiliary processes at time
$t$ and started at the positions $Z_T$. We recall that $L$ is the
generator defined in (\ref{eq:defLexple}) and let $J$ be the operator
defined on the space of locally integrable functions by
%
%e75 ###
%
\begin{eqnarray}\label{eq:defJ}
Jf(x) &=& -\frac{3r}{2} f(x)+r \int_{0}^1 \bigl(f(qx)+f\bigl((1-q)x\bigr)
\bigr)G(dq)\nonumber\\[-9pt]\\[-9pt]
&=& -\frac{3r}{2} f(x)+2r \int_{0}^1 f(qx)G(dq).\nonumber
\end{eqnarray}
This operator will naturally appear when computing the equation
satisfied by $\eta^T$ by applying (\ref{pbmsplitteddiff}) with
$f_t(x)=\rme^{-rt/2}f(x)$.
\begin{prop}\label{propeqfluctu}The fluctuation process (\ref
{fluctuationsdef}) satisfies the following evolution equation:
%
%e76 ###
%
\begin{equation}
\label{fluctuations}
\langle\eta^T_t,f\rangle= \int_0^t \int_{\R} \bigl(Lf(x)+Jf(x)
\bigr)\eta^T_s(dx) \,ds +
M_t^{T}(f),
\end{equation}
where $M^{T}_t(f)$ is a square integrable martingale with quadratic variation,
%
%e77 ###
%
\begin{eqnarray}\label{crochetfluctuations}\quad
\langle M^{T}(f)\rangle_t
&=& \int_0^tds
\int_\R\frac{Z_{s+T}(dx)}{\E[N_{s+T}]}\nonumber\\[-2pt]
&&\hphantom{\int_0^tds \int_\R}
{}\times\biggl[
r\int_0^1
\bigl(f(qx)+f\bigl((1-q)x\bigr)-f(x)\bigr)^2 G(dq)\\[-2pt]
&&\hspace*{178.5pt}{}+2\sigma^2(x)f'(x)^2 \biggr].
\nonumber
\end{eqnarray}
\end{prop}

The proof of this proposition is given in Section
\ref{sectionprooftcl}. In the following, we are interested in the
behavior of the fluctuation process when $T\rightarrow+\infty$. The
processes $\eta^T$ take their values in the space $\mathcal{M}_S(\R
)$ of
signed measures. Since this space endowed with the topology of weak
convergence is not metrizable, we follow the approach of M\'{e}tivier
\cite{metivierIHP} and M\'{e}l\'{e}ard \cite{meleardfluctuation} (see also
\cite{ferrieretran,chithese}) and embed $\mathcal{M}_S(\R)$ in weighted
distribution spaces. This is described in the sequel. We then prove the
convergence of the fluctuation processes to a~distribution-valued
diffusion driven by a Gaussian white noise (Proposition~\ref{propconvergencefluctuations}).

%s6.2 ###
\subsection{Convergence of the fluctuation process: The central limit theorem}\label{sec6.2}

Let us introduce the Sobolev spaces that we will use (see, e.g., Adams
\cite{adams}). We follow in this the steps of
\cite{metivierIHP,meleardfluctuation}. To obtain estimates\vadjust{\goodbreak} of our fluctuation processes,
the following additional regularities for $b$ and $\sigma$ are
required as well as assumptions on our auxiliary process.
\begin{hyp}\label{hyptcl}
We assume the following:

\begin{longlist}
\item
$b$ and $\sigma$ are in $\Co^{8}(\R,\R)$ with bounded
derivatives.
\item
There exists $K>0$ such that for every $|x|>K$, $2 b(x)/x+6 \sigma
(x)/\break x^2<r'$ with $r'<r\int_0^1 (1-q^4)^2 G(dq)$.
\item
$Y$ is ergodic with stationary measure $\pi$ such that $\langle
\pi, |x|^8\rangle<+\infty$.
\item
For every initial condition $\mu$ such that $\langle\mu
,|x|^8\rangle<+\infty$,\break $\sup_{t\in\R_+}\E_\mu[Y_t^{8}]<+\infty$.
\end{longlist}
\end{hyp}
\begin{rque}\label{rquehyptcl}
\begin{longlist}
\item
Notice that under Assumption \ref{hyptcl}(i), there exist $ \bar
{b}$ and
$\bar{\sigma}>0$ s.t. for all $x\in\R$, we have $|b(x)|\leq
\bar{b}(1+|x|)$ and $|\sigma(x)|\leq\bar{\sigma}(1+|x|)$.
\item
 Conditions for the ergodicity of $Y$ have been provided in
Proposition~\ref{propapplimeyntweedie} and Remarks \ref
{rquemeyntweedie} and \ref{rquemeyntweedie2}. Under Assumption \ref
{hyptcl}(ii), Remark \ref{rquemeyntweedie2} applies and we have
geometrical ergodicity with (\ref{exponentielltgeom2}).
\item
 The moment hypothesis of Assumption \ref{hyptcl}(iv) is
fulfilled for the examples (i)--(iii) of Section \ref
{sectionexspitteddiffusion} provided the initial condition satisfies
$\langle\mu, |x|^8\rangle<+\infty$. This can be seen by using It\^
{o}'s formula (e.g., \cite{ikedawatanabe}, Theorem 5.1, page~67) and
Gronwall's lemma. Moreover, for every $p\in\{1,\ldots,8\}$, $\E_\mu
[|Y_t|^p] <+\infty$.
\item
Assumption \ref{hyptcl}(iii) and (iv) imply that $\forall p\in\{
1,\ldots,7\}$, $\int_{\R}|x|^p \pi(dx)<+\infty$ and $\lim
_{t\rightarrow+\infty}\E_\mu[|Y|^p]=\int_{\R}|x|^p \pi(dx)$.
This is a consequence of the equi-integrability of $(|Y_t|^p)_{t\geq
0}$ for $p\in\{1,\ldots,7\}$.
\end{longlist}
\end{rque}

For $j\in\N$ and $\alpha\in\R_+$, we denote by $W^{j,\alpha}$ the
closure of $\Co^\infty(\R,\R)$ with respect to the norm
%
%e78 ###
%
\begin{equation}
\|g\|_{W^{j,\alpha}}:=\biggl(\sum_{k\leq j}\int_{\R}\frac
{|g^{(k)}(x)|^2}{1+|x|^{2\alpha}}\,dx\biggr)^{1/2},
\end{equation}
where $g^{(k)}$ is the $k$th derivative of $g$. The
space $W^{j,\alpha}$ endowed with the norm \mbox{$\|\cdot\|_{W^{j,\alpha}}$}
defines a Hilbert space. We denote by $W^{-j,\alpha}$ the dual
space. Let~$C^{j,\alpha}$ be the space of functions $g$ with $j$
continuous derivatives and such that\looseness=-1
\[
\forall k\leq j\qquad \lim_{|x|\rightarrow+\infty}\frac
{|g^{(k)}(x)|}{1+|x|^\alpha}=0.
\]\looseness=0
When endowed with the norm
%
%e79 ###
%
\begin{equation}
\|g\|_{C^{j,\alpha}}:=\sum_{k\leq j}\sup_{x\in\R}\frac
{|g^{(k)}(x)|}{1+|x|^\alpha},
\end{equation}
these spaces are Banach spaces and their dual spaces are denoted
by\vadjust{\goodbreak}
$C^{-j,\alpha}$.

In the sequel, we will use the following embeddings (see
\cite{adams,meleardfluctuation}):
%
%e80 ###
%
\begin{eqnarray}\label{emboitement}
&& C^{7,0}\hookrightarrow W^{7,1}\hookrightarrow_{\mathrm{H.S.}}
W^{5,2}\hookrightarrow C^{4,2}\hookrightarrow W^{4,3}\hookrightarrow
C^{3,3}\nonumber\\
&&\qquad\hookrightarrow W^{3,4}\hookrightarrow C^{2,4}
C^{-2,4}\hookrightarrow W^{-3,4}\hookrightarrow
C^{-3,3}\\
&&\qquad\hookrightarrow W^{-4,3}\hookrightarrow C^{-4,2}\hookrightarrow
W^{-5,2}\hookrightarrow_{\mathrm{H.S.}} W^{-7,1}\hookrightarrow
C^{-7,0},\nonumber
\end{eqnarray}
where $\mathrm{H.S.}$ means that the corresponding embedding is
Hilbert--Schmidt (see \cite{adams}, page 173). Let us explain briefly
why we
use these embeddings. Following\vspace*{1pt} the preliminary estimates
of~\cite{meleardfluctuation} (Proposition 3.4), it is possible to choose~%
$W^{-3,4}$ as a reference space for our study. We control the norm of
the martingale part in $W^{-4,3}$ using the embeddings
$W^{4,3}\hookrightarrow C^{3,3}\hookrightarrow W^{3,4}$. We obtain
uniform estimate for the
norm of $\eta^T_t$ in $C^{-4,2}$. The spaces~$W^{-5,2}$ and~$W^{-7,1}$
are used to apply the tightness criterion in
\cite{meleardfluctuation} (see our Lem\-ma~\ref{lemmemetiviermeleard}).
The space $C^{-7,0}$ is used for proving
uniqueness of the accumulation point of the family $(\eta^T)_{T\geq0}$.
\begin{prop}\label{propconvergencefluctuations}
Let $\Upsilon>0$. The sequence $(\eta^T)_{T\in\R_+}$ converges in\break
$\D([0, \Upsilon],C^{-7,0})$ when $T\rightarrow+\infty$ to the
unique solution in $C([0,\Upsilon],C^{-7,0})$ of the following
evolution equation:
%
%e81 ###
%
\begin{equation}\label{limitefluctuations}
\langle\eta_t,f\rangle= \int_0^t \int_{\R} \bigl(Lf(x)+Jf(x)
\bigr)\eta_s(dx) \,ds + \sqrt{W}\mathcal{W}_t(f),
\end{equation}
where $\mathcal{W}(f)$ is a Gaussian martingale independent of $W$ and
which bracket is $V(f) \times t$ with
%
%e82 ###
%
\begin{eqnarray}\label{defV}
V(f)&=&\int_{\R} \biggl(r\int_0^1 \bigl(f(qx)+f\bigl((1-q)x\bigr)-f(x)\bigr)^2
G(dq)\nonumber\\[-8pt]\\[-8pt]
&&\hspace*{139pt}{} +
2\sigma^2(x) f'(x)^2\biggr) \pi(dx).\nonumber
\end{eqnarray}
\end{prop}

Notice that unlike the discrete case treated in \cite{delmasmarsalle},
our fluctuation process here has a finite variational part.

%s6.3 ###
\subsection{Proofs}\label{sectionprooftcl}

We begin by establishing the evolution equation for $\eta^T$ that is
announced in Proposition \ref{propeqfluctu}.
\begin{pf*}{Proof of Proposition \ref{propeqfluctu}}
From Lemma \ref{lemmeNt} and applying (\ref{pbmsplitteddiff}) with
$f_t(x)=\rme^{-rt/2}f(x)$, we obtain
%
%e83 ###
%
\begin{eqnarray}\label{tcl1}
&&
\langle Z_{t+T},f\rangle\rme^{-r(t+T)/2} \nonumber\\
&&\qquad = \langle Z_T,f\rangle
\rme^{-rT/2}+M^{T}_t(f) \\
%+ \int_T^{t+T}\int_{\R}(Lf(x)-\frac{r}{2}f(x)-r\int_0^1
%(f(qx)+f((1-q)x)-f(x)) G(dq) )\expp{-rs/2}Z_{s}(dx) \,ds\\
%= \langle Z_T,f\rangle\expp{-rT/2}+ M^{f,T}_t
&&\qquad\quad{} + \int_0^{t}\int_{\R}\bigl(Lf(x)+Jf(x)\bigr)\rme
^{-r(s+T)/2}Z_{s+T}(dx)
\,ds,\nonumber
\end{eqnarray}
where $M^{T}_t(f)$ is a square integrable martingale with quadratic variation
%
%e84 ###
%
\begin{eqnarray}\qquad
\langle M^{T}(f)\rangle_t&=& \int_T^{t+T} ds
\int_{\R}\rme^{-rs}Z_{s}(dx) \nonumber\\
&&\hspace*{50.1pt}{}\times\biggl[r\int_0^1
\bigl(f(qx)+f\bigl((1-q)x\bigr)-f(x)\bigr)^2
G(dq)\\
&&\hspace*{189.1pt}{}+2\sigma^2(x)f'(x)^2\biggr],
\nonumber
\end{eqnarray}
which is the bracket announced in (\ref{crochetfluctuations}).
Computing $\langle Z_{t},f\rangle\rme^{-rt/2}$
in the same way and taking the
expectation gives, with (\ref{eq:Qt}) and Proposition
\ref{propprocessusauxiliaire},
\[
Q_tf(x)\rme^{rt/2}= f(x)+\int_0^t Q_s( Lf+Jf)(x)\rme^{rs/2}\,ds.
\]
Integrating with respect to $Z_T$ and multiplying by $\rme^{-rT/2}$ implies
%
%e85 ###
%
\begin{eqnarray}
\label{tcl2}
&&\langle Z_{T},Q_tf\rangle\rme^{-r(T-t)/2}\nonumber\\[-8pt]\\[-8pt]
&&\qquad= \langle Z_T, f\rangle
\rme^{-rT/2}
+ \int_0^t \rme^{-r(T-s)/2}\,ds \langle Z_T, Q_s( Lf+Jf
)\rangle.\nonumber
\end{eqnarray}
We deduce the announced result from (\ref{fluctuationsdef}), (\ref
{tcl1}) and (\ref{tcl2}).
\end{pf*}

We now prove that our fluctuation process $\eta^T$ can be viewed as a
process with values in $W^{-3,4}$ by following the preliminary
estimates of \cite{meleardfluctuation} (Proposition 3.4). This space
$W^{-3,4}$ is then chosen as reference space and in all the spaces
appearing in the second line of (\ref{emboitement}) that contain
$W^{-3,4}$, the norm of~$\eta^T_t$ is finite and well defined.
\begin{lemme}\label{lemmecontroleprelim}Let $\Upsilon>0$. There
exists a finite constant $C$ that does not depend on $T$ nor on
$\Upsilon$ such that
%
%e86 ###
%
\begin{equation}
\sup_{t\in[0,\Upsilon]}\E_\mu[ \|\eta^T_t\|_{W^{-3,4}}^2
]\leq C\rme^{r(\Upsilon+T)}.
\end{equation}
\end{lemme}
\begin{pf}
Let $(\varphi_p)_{p\in\N^*}$ be a complete orthonormal basis of
$W^{3,4}$ that are~$\Co^{\infty}$ with compact support. We have by
Riesz representation theorem and Parseval's identity
%
%e87 ###
%
\begin{eqnarray}\label{etape444}\quad
\rme^{-r(t+T)}\|\eta^T_t\|_{W^{-4,3}}^2 & = & \rme^{-r(t+T)} \sum
_{p\geq1} \langle\eta^T_t,\varphi_p\rangle^2 \nonumber\\
& = & \rme^{-r(t+T)}\E[N_{t+T}]\sum_{p\geq1}\biggl(\frac{\langle
Z_{t+T},\varphi_p\rangle}{\E[N_{t+T}]}-\frac{\langle Z_T,Q_t\varphi
_p\rangle}{\E[N_T]}\biggr)^2\\
&\leq& 2 \sum_{p\geq1}\biggl(\frac{\langle Z_{t+T},\varphi_p\rangle
^2}{\E[N_{t+T}]^2}+\frac{\langle Z_T,Q_t\varphi_p\rangle^2}{\E
[N_T]^2}\biggr).\nonumber
\end{eqnarray}
Under the Assumption \ref{hyptcl}(iii) and thanks to Remark \ref
{rquecontientcontinuborne} and Example \ref{ex:YuleNt}, we use the
same proof as in Theorem \ref{thLGN}, especially (\ref{carre}) and
(\ref{eq:Btcv}), % and obtain that for any bounded continuous function
%$f$ such that $\langle\pi,f\rangle=0$:
% \lim_{t\rightarrow+\infty}\E_\mu[\frac{\langle Z_t,f\rangle^2}{
%For bounded continuous functions $f$, since $
% \E_\mu[\frac{\langle Z_t,f\rangle^2}{\E[N_t]^2}]\leq& 2 \E_
%]+2 \frac{\E[N_t^2]}{\E[N_t]^2} \langle\pi,f\rangle^2.
%Th. \ref{thLGN} and the second term to $4\langle\pi,f\rangle^2$ by (
%
%e88 ###
%
\begin{eqnarray} \label{etape4}\hspace*{20pt}
0&< &\E_\mu\biggl[\frac{\langle Z_{t+T},\varphi_p\rangle^2}{\E
[N_{t+T}]^2}+\frac{\langle Z_T,Q_t\varphi_p\rangle^2}{\E
[N_T]^2}\biggr]\nonumber\\[-2pt]
& = & \rme^{-r(t+T)}\mu Q_{t+T}\varphi_p^2 +r\int_0^{t+T} \mu
Q_sJ_2(Q_{t+T-s}\varphi_p\otimes Q_{t+T-s}\varphi_p)\rme
^{-rs}\,ds\nonumber\\[-10pt]\\[-10pt]
&&{} + \rme^{-rT}\mu Q_{T}(Q_t\varphi_p)^2 +r\int_0^{T} \mu
Q_sJ_2(Q_{T-s}Q_t\varphi_p\otimes Q_{T-s}Q_t\varphi_p)\rme
^{-rs}\,ds\nonumber\\[-2pt]
&\leq& 2\rme^{-rT} \mu Q_{t+T} \varphi_p^2 + 4r \int_0^{t+T} \int
_0^1 \int_\R\varphi_p^2(qx)\rme^{-rs} \mu Q_{t+T}(dx) G(dq)
\,ds,\nonumber
\end{eqnarray}
since by (\ref{eq:J2-0}), the Cauchy--Schwarz inequality and symmetry of $G$,
\begin{eqnarray*}
&& J_2(Q_{t+T-s}|\varphi_p|\otimes Q_{t+T-s}|\varphi_p|)(x)\\[-2pt]
&&\qquad= 2\int_0^1 \bigl( Q_{t+T-s}|\varphi_p|(qx) Q_{t+T-s}|\varphi
_p|\bigl((1-q)x\bigr)\bigr) G(dq)\\[-2pt]
&&\qquad\leq 2 \int_0^1 Q_{t+T-s}\varphi_p^2(qx) G(dq).
\end{eqnarray*}
We deduce from (\ref{etape444}) and (\ref{etape4}) that
%
%e89 ###
%
\begin{eqnarray}\label{etape5}
&&\rme^{-r(t+T)}\E_\mu[\|\eta^T_t\|_{W^{-4,3}}^2]\nonumber\\[-2pt]
&&\qquad\leq
4\rme^{-rT} \int_{\R} \sum_{p\geq1}\varphi_p^2(x) \mu
Q_{t+T}(dx)\\[-2pt]
&&\qquad\quad{} + 8r \int_0^{t+T} \int_0^1 \int_\R\sum_{p\geq1}\varphi
_p^2(qx)\rme^{-rs} \mu Q_{t+T}(dx) G(dq) \,ds.\nonumber
\end{eqnarray}
Let us consider the linear forms $D_{x,q}(g)=g(qx)$ for $q\in[0,1]$,
$x\in\R$ and $g\in W^{3,4}\hookrightarrow C^{2,4}$,
\[
|D_{x,q}(g)|= |g(qx)|\leq(1+|x|^4)\|g\|_{C^{2,4}}\leq C(1+|x|^4)\|g\|
_{W^{3,4}.}
\]
Using Riesz representation theorem and Parseval's identity, we get
%
%e90 ###
%
\begin{equation}\label{limiteetape2}
\sum_{p\geq1}D_{x,q}(\varphi_p)^2=\|D_{x,q}\|^2_{W^{-3,4}}\leq
C(1+|x|^4).
\end{equation}
We deduce from (\ref{eq:Qt}), (\ref{etape5}) and Assumption \ref
{hyptcl}(iv) that
%
%e91 ###
%
\begin{eqnarray}\label{etape6}
&&\rme^{-r(t+T)}\E_\mu[\|\eta^T_t\|_{W^{-4,3}}^2]\nonumber\\[-10pt]\\[-10pt]
&&\qquad\leq
C \E_\mu[1+|Y_{t+T}|^4]\biggl(\rme^{-rT}+ \frac{1- \rme
^{-r(t+T)}}{r}\biggr)\leq C,\nonumber
\end{eqnarray}
where the constant $C$ is finite and does not depend on $\Upsilon$ nor
$T$. The proof is complete.\vadjust{\goodbreak}
\end{pf}

We now turn to the proof of the central limit theorem stated in
Proposition \ref{propconvergencefluctuations}. To achieve this aim, we
first prove Lemma \ref{lemmeestimeemomenttcl}.\vspace*{-3pt}

\begin{lemme}\label{lemmeestimeemomenttcl}Suppose that Assumption
\ref{hyptcl} is satisfied and let $\Upsilon\in\R_+$.

\begin{longlist}
\item
We have
%
%e92 ###
%
\begin{equation}\label{estimeeuniforme}
\sup_{T\in\R_+}\sup_{t\leq\Upsilon}\E_\mu[\|\eta^T_t\|
^2_{C^{-4,2}}]<+\infty.
\end{equation}
\item Let us denote by $M^T_t$ the operator that associates
$M^{T}_t(f)$ to $f$. Then
%
%e93 ###
%
\begin{equation}
\label{estimeeuniformemartingale}
\sup_{T\in\R_+}\sup_{t\leq
\Upsilon}\E_\mu[\|M^T_t\|^2_{W^{-4,3}}]<+\infty.\vspace*{-3pt}
\end{equation}
\end{longlist}
\end{lemme}

\begin{pf}
Let us first deal with (\ref{estimeeuniformemartingale}).
We consider the following linear forms:
$D_{x,\sigma}(g)=\sigma(x) g'(x)$ and
$D_{x,q}(g)=g(qx)+g((1-q)x)-g(x)$. Notice
that for $g\in W^{4,3}\hookrightarrow C^{3,3}$, $x\in\R$ and $q\in[0,1]$,
%
%e94 ###
%
\begin{eqnarray} \label{etape4545}
|D_{x,\sigma}(g)| &=& |\sigma(x)g'(x)|\leq\bar{\sigma
}(1+|x|)|g'(x)|\nonumber\\[-2pt]
& \leq&
C(1+|x|^4)\|g\|_{C^{3,3}}\leq
C(1+|x|^4)\|g\|_{W^{4,3}},\nonumber\\[-10pt]\\[-10pt]
|D_{x,q}(g)| &=& \bigl|g(qx)+g\bigl((1-q)x\bigr)-g(x)\bigr|\nonumber\\[-2pt]
& \leq& 3(1+|x|^3)\|g\|_{C^{3,3}}\leq
C(1+|x|^3)\|g\|_{W^{4,3}},\nonumber
\end{eqnarray}
where $C$ is not dependent on $x$ nor on $q$. This implies that
$D_{x,\sigma}$ and $D_{x,q}$
are continuous from $W^{4,3}$ into $\R$ and their norms in $W^{-4,3}$
are upper bounded by $C(1+|x|^4)$ and $C(1+|x|^3)$, respectively.
Let us consider a sequence of functions $(\varphi_p)_{p\in\N^*}$
constituting a complete orthonormal basis
of $W^{4,3}$ and that are $\Co^\infty$ with compact support. Using Riesz
representation theorem and Parseval's identity, we get
%
%e95 ###
%
\begin{eqnarray}\label{eq:op-bound1}
\sum_{p\geq1} D_{x,\sigma}(\varphi_p)^2&=&\|D_{x,\sigma}\|
_{W^{-4,3}}^2\leq C
(1+|x|^8),\nonumber\\[-10pt]\\[-10pt]
\sum_{p\geq1} D_{x,q}(\varphi_p)^2&=&\|D_{x,q}\|_{W^{-4,3}}^2\leq C
(1+|x|^6).\nonumber
\end{eqnarray}
%
%where $C$ does not dependent on $x$ or $q$.
We have
%
%e96 ###
%
\begin{eqnarray}\label{majorationnormemartingale}
&& \E_\mu\Bigl[ \sup_{t\leq\Upsilon} \|M_t^T\|_{W^{-4,3}}^2
\Bigr]\nonumber\\[-3pt]
&&\qquad
\leq\E_\mu\biggl[\sum_{p\geq1}\sup_{t\leq\Upsilon}M^T_t(\varphi
_p)^2\biggr]\nonumber\\[-3pt]
&&\qquad\leq
4 \sum_{p\geq1}\E_\mu[\langle
M^T(\varphi_p)\rangle_{\Upsilon}]\nonumber\\[-3pt]
&&\qquad=4 \int_T^{T+\Upsilon}ds\, \E_\mu\biggl[\int_\R\frac{Z_s(dx)}{\E[N_s]}
\biggl(r\int_0^1 \sum_{p\geq1} D_{x,q}(\varphi_p)^2 G(dq) \\[-3pt]
&&\hspace*{190.8pt}{} + 2
\sum_{p\geq1} D_{x,\sigma}(\varphi_p)^2 \biggr)\biggr]\nonumber\\[-3pt]
&&\qquad\leq C \int_T^{T+\Upsilon}ds\, \E_\mu\biggl[\int_\R\frac
{Z_s(dx)}{\E[N_s]}(1+|x|^8)\biggr]\nonumber\\[-3pt]
&&\qquad= C \int_T^{T+\Upsilon} ds\, \E_\mu[(1+|Y_s|^8)],\nonumber
\end{eqnarray}
where the first inequality comes from \cite{adams}, Lemma 6.52, the
second is Doob's inequality, the third line is a consequence of
(\ref{crochetfluctuations}), the fourth inequality comes from the
bounds (\ref{eq:op-bound1}) and
the last equality comes from (\ref{defYt}). The proof is then finished
since by Assumption \ref{hyptcl}(iv), $\sup_{t\geq0} \E_\mu
[Y^8_t]<\infty$.

%***************
%A) Mq $\sup_{t\geq0} \E_\mu[Y^4_t]<\infty$.
%*************

Let us now consider the proof of (\ref{estimeeuniforme}). Recall $J$
defined by (\ref{eq:defJ}). It is clear that $J$ is a bounded operator
from $C^{4,2}$ into itself
%
%e97 ###
%
\begin{equation}
\label{eq:Jbound}
\|J\varphi\|_{C^{4,2}}\leq C \|\varphi\|_{C^{4,2}},
\end{equation}
where $C$ does not depend on $\varphi\in C^{4,2}$.

Let us denote\vspace*{2pt}
by $U(t)$ the semigroup of the diffusion with generator $L$ given by
(\ref{eq:defLexple}). Proposition 3.9 in
\cite{meleardfluctuation} and Assumption \ref{hyptcl}
yield that for $\varphi\in C^{4,2}$ and $\psi\in C^{3,3}$
%
%e98 ###
%
\begin{eqnarray}
\label{eq:Ubound}
\sup_{t\leq\Upsilon} \|U(t)(\varphi)\|_{C^{4,2}}&\leq& C \|\varphi\|
_{C^{4,2}}\quad\mbox{and}\nonumber\\[-8pt]\\[-8pt]
\sup_{t\leq\Upsilon} \|U(t)(\psi)\|_{C^{3,3}}&\leq& C
\|\psi\|_{C^{3,3}},\nonumber
\end{eqnarray}
where $C$ does not depend on $\varphi$ nor on $\psi$.

Let us consider the test function $\psi_t \dvtx(s,x)\mapsto
U(t-s)\varphi(x)$ with $\varphi\in C^{4,2}$. Using It\^{o}'s formula
\[
\langle\eta^T_t,\varphi\rangle= \int_0^t \langle
\eta^T_s,JU(t-s)\varphi\rangle \,ds+\int_0^t \langle
\,dM^T_s,U(t-s)\varphi\rangle,
\]
that is, $
\eta^T_t= \int_0^t U(t-s)^* J^*
\eta^T_s \,ds+\int_0^t U(t-s)^* \,dM^T_s$, where $U(t-s)^*$ and~%
$J^*$ stand for the adjoint operators of $U(t-s)$ and $J$.
Then, for $t\leq\Upsilon$
%
%e99 ###
%
\begin{eqnarray}
\label{eq:majoeta}
&&\E_\mu[\|\eta^T_t\|^2_{C^{-4,2}}]\nonumber\\
&&\qquad\leq
2\Upsilon\int_0^t \E_\mu[\|U(t-s)^* J^* \eta^T_s\|
^2_{C^{-4,2}}]\,ds\\
&&\qquad\quad{}+2\E_\mu\biggl[\biggl\|\int_0^t U(t-s)^* \,dM^T_s\biggr\|^2_{C^{-4,2}}\biggr].
\nonumber
\end{eqnarray}
Thanks to (\ref{eq:Jbound}) and (\ref{eq:Ubound}), we have for $s\leq
t\leq\Upsilon$,
%
%e100 ###
%
\begin{equation}\label{etape4646}
\E_\mu[\|U(t-s)^* J^* \eta^T_s\|^2_{C^{-4,2}}]
\leq C \E_\mu[\|\eta^T_s\|^2_{C^{-4,2}}].
\end{equation}
The second term of the right-hand side  of (\ref{eq:majoeta}) is upper bounded
by considering the norm in $W^{-4,3}$. To prove that
%
%e101 ###
%
\begin{equation}\label{etape4647}
\sup_{T\in\R_+} \sup_{t\leq\Upsilon} \E_\mu\biggl[\biggl\|\int_0^t U(t-s)^*
\,dM^T_s\biggr\|^2_{W^{-4,3}}\biggr] <+\infty,
\end{equation}
we use similar arguments as those used for the proof of
(\ref{estimeeuniformemartingale}) and
(\ref{eq:Ubound}). In the proof below, we replace\vspace*{1pt} the linear forms
$D_{x,\sigma}$ and
$D_{x,q}$ by $\bar D_{x,t-s,\sigma}$ and $\bar D_{x,t-s,q}$ with
$\bar D_{x,t-s,\sigma} (\varphi)=D_{x,\sigma}(U(t-s)\varphi) $ and
$\bar D_{x,t-s,q} (\varphi)=D_{x,q}(U(t-s) \varphi)$.\vadjust{\goodbreak}
%DETAIL DE LA PREUVE
Notice that by (\ref{etape4545}) for $g\in W^{4,3} \hookrightarrow
C^{3,3}$, $x\in\R$ and $q\in[0,1]$,
\begin{eqnarray*}
|\bar D_{x,t,\sigma}(g)| & = &|D_{x,\sigma}(U(t)g)|\leq C(1+|x|^4)\|
U(t)g\|_{C^{3,3}}\\
&\leq& C(1+|x|^4)\|g\|_{C^{3,3}}
\leq C(1+|x|^4)\|g\|_{W^{4,3}},\\
|\bar D_{x,t,q}(g)| &=& |D_{x,q}(U(t)g)|\leq C(1+|x|^3)\|U(t)g\|
_{C^{3,3}}\\
&\leq&
C(1+|x|^3)\|g\|_{W^{4,3}},
\end{eqnarray*}
where $C$ is not dependent on $x$.
%Using again that $(\varphi_p)_{p\in\N^*}$ is a complete orthonormal
%basis of $W^{4,3}$ that are $\Co^\infty$ with compact support.
Using again Riesz
representation theorem and Parseval's identity, we get
\begin{eqnarray*}
\sum_{p\geq1} \bar D_{x,t,\sigma}(\varphi_p)^2&=&\|\bar D_x\|
_{W^{-4,3}}^2\leq C
(1+|x|^8),\\
\sum_{p\geq1} \bar D_{x,t,q}(\varphi_p)^2&=&\|\bar D_{x,q}\|
_{W^{-4,3}}^2\leq C
(1+|x|^6),
\end{eqnarray*}
where $C$ is not dependent on $x$ nor on $q$.
We have with the same arguments as in~(\ref{majorationnormemartingale}),
\begin{eqnarray*}
&&\E_\mu\biggl[ \sup_{t\leq\Upsilon} \biggl\|\int_0^t U(t-s)^* \,dM_s^T\biggr\|
_{W^{-4,3}}^2 \biggr]\\[-2pt]
&&\qquad \leq\E_\mu\biggl[\sum_{p\geq1}\sup_{t\leq\Upsilon}\int_0^t
\bigl(U(t-s)\varphi_p\bigr)^2 \,dM_s^T\biggr]\\[-2pt]
&&\qquad\leq
4 \sum_{p\geq1}\E_\mu\biggl[\int_0^\Upsilon\bigl(U(t-s)\varphi_p\bigr)^2\,
d\langle
M^T\rangle_{s}\biggr]\\[-2pt]
&&\qquad=4 \int_T^{T+\Upsilon} ds\, \E_\mu\biggl[\int_\R
\biggl(r\int_0^1 \sum_{p\geq1} \bar D_{x,t-s,q}(\varphi_p)^2
G(dq) \\[-2pt]
&&\qquad\quad\hspace*{119.2pt}{} + 2
\sum_{p\geq1} \bar D_{x,t-s,\sigma}(\varphi_p)^2 \biggr)\frac
{Z_s(dx)}{\E[N_s]}\biggr]\\[-2pt]
%&\leq C \int_T^{T+\Upsilon} ds \E_\mu[\int_\R\frac{Z_s(dx)}{
&&\qquad\leq C \int_T^{T+\Upsilon} ds \,\E_\mu[(1+|Y_s|^8)].
\end{eqnarray*}
The proof\vspace*{2pt} is then complete as $\sup_{t\geq0} \E_\mu[Y^8_t]<\infty
$ by Assumption \ref{hyptcl}(iv).\vadjust{\goodbreak}

%**************

%Mq $\sup_{t\geq0} \E_\mu[Y^8_t]<\infty$.

%**************

%*****************

%Or la majoration dans $W^{-4,2}$ ne donne pas la majoration dans
%$C^{-3,2}$ mais dans $C^{-4,1}$.
%Il faut donc faire la majoration \reff{estimeeuniforme} dans $C^{-4,1}$
%et donc supposer $b$ et $\sigma$ $C^5$. On utilise alors
%terme de droite de \reff{eq:majoeta} et \reff{eq:Ubound} dans $C^{3,2}$
%pour le second terme de droite.

Thus, we get from (\ref{eq:majoeta}), (\ref{etape4646}) and (\ref{etape4647})
\[
\E_\mu[\|\eta^T_t\|^2_{C^{-4,2}}]\leq C\biggl(1+\int_0^t
\E_\mu[\|\eta^T_s\|^2_{C^{-4,2}}]\,ds\biggr).
\]
We use Gronwall's lemma and the fact that
$\E_\mu[\|\eta^T_t\|^2_{C^{-4,2}}]$ is locally bounded
(see Lemma \ref{lemmecontroleprelim}) to conclude.\vspace*{-2pt}
\end{pf}

We now prove the tightness of the fluctuation process.\vspace*{-2pt}
\begin{prop}\label{proptension}Let $\Upsilon>0$. The sequence $(\eta
^T)_{T\in\R_+}$ is tight in $\D([0,\Upsilon]$, $W^{-7,1})$.\vspace*{-2pt}
\end{prop}

We use a tightness criterion from \cite{joffemetivier}, which we
recall (see \cite{meleardfluctuation}, Lemma~C, page 217).\vspace*{-2pt}
\begin{lemme}\label{lemmemetiviermeleard}
A sequence $(\Theta^T)_{T\in\R_+}$ of Hilbert $H$-valued c\`{a}dl\`
{a}g processes is tight in $\D([0,\Upsilon],H)$ if the following
conditions are satisfied:

\begin{longlist}
\item
There exists a Hilbert space $H_0$ such that $H_0\hookrightarrow
_{\mathrm{H.S.}} H$ and $\forall t\leq\Upsilon$, $\sup_{T\in\R_+}\E[\|
\Theta^T_t\|_{H_0}^2]<+\infty$.
\item
(Aldous condition). For every $\varepsilon>0$, there exists
$\delta>0$ and $T_0\in\R_+$ such that for every sequence of stopping
time $\tau_T\leq\Upsilon$,
\[
\sup_{T>T_0}\sup_{\varsigma<\delta}\PP(\|\Theta^T_{\tau
_T+\varsigma}-\Theta^T_{\tau_T}\|_{H} > \varepsilon)<
\varepsilon.\vspace*{-2pt}
\]
\end{longlist}
\end{lemme}
\begin{pf*}{Proof of Proposition \ref{proptension}}
We shall use Lemma \ref{lemmemetiviermeleard} with $H_0=W^{-5,2}$ and
$H=W^{-7,1}$.
Condition (i) is a direct consequence of the uniform estimates
obtained in
(\ref{estimeeuniforme}) and of the fact that
$\|\eta^T_t\|^2_{W^{-5,2}}\leq C \|\eta^T_t\|^2_{C^{-4,2}}$.\vspace*{1pt}

Let us now turn to condition (ii). By the Rebolledo criterion (see,
e.g.,~\cite{joffemetivier}), it is sufficient to show the Aldous condition
for the finite variation part and for the trace of the martingale\vspace*{2pt}
part of~(\ref{fluctuations}). Let $(\varphi_p)_{p\geq1}$ be
a complete orthonormal system of $W^{7,1}\hookrightarrow C^{6,1}$. We
recall that the trace of the martingale part is defined as $
\operatorname{tr}_{W^{-7,1}}\dlangle M^T \drangle_{t}=\sum_{p\geq1}\langle
M^T(\varphi_p)\rangle_t$ (see, e.g.,~\cite{joffemetivier}). Let
$\varepsilon>0$ and let $\tau_T\leq
\Upsilon$ be a sequence of stopping times. For $T_0>0$ and
$\delta>0$, following the steps of (\ref{majorationnormemartingale}),
we get
%
%e102 ###
%
\begin{eqnarray} \label{etape4547}
&&
\sup_{T>T_0}\sup_{\varsigma<\delta}\PP(|{\operatorname{tr}}_{W^{-7,1}}\dlangle M^T \drangle_{\tau_T+\varsigma}-
\operatorname{tr}_{W^{-7,1}}\dlangle M^T
\drangle_{\tau_T}|>\varepsilon)\nonumber\\[-2pt]
&&\qquad
\leq \sup_{T>T_0}\sup_{\varsigma<\delta}\frac{1}{\varepsilon}\E
\biggl[\int_{\tau_T}^{\tau_T+\varsigma} \biggl\langle\frac{Z_{s+T}}{\E
[N_{s+T}]},\nonumber\\[-10pt]\\[-10pt]
&&\hspace*{130pt}r\int_0^1 \sum_{p\geq1} D_{x,q}(\varphi_p)^2 G(dq)\nonumber\\
&&\hspace*{162.2pt}{}+2
\sum_{p\geq1}D_{x,\sigma}(\varphi_p)^2 \biggr\rangle \,ds\biggr].\nonumber
\end{eqnarray}
Using the embedding $W^{7,1}\hookrightarrow C^{6,1}$ and computations
similar to (\ref{etape4545}),
\begin{eqnarray*}
\sum_{p\geq1} D_{x,q}(\varphi_p)^2&=&\|D_{x,q}\|^2_{W^{-7,1}}\leq
C(1+|x|^2),\\
\sum_{p\geq1} D_{x,\sigma}(\varphi_p)^2&=&\|D_{x,\sigma}\|
^2_{W^{-7,1}}\leq C(1+|x|^4).
\end{eqnarray*}
%
%Indeed, for $g\in W^{7,1}\hookrightarrow C^{6,1}$:
% & |\sigma(x)g'(x)|\leq\bar{\sigma}(1+|x|)|g'(x)|\leq C(1+|x|^2)\|g
% & |g(qx)+g((1-q)x)-g(x)|\leq3(1+|x|)\|g\|_{C^{6,1}}\leq C(1+|x|)\|g
%& 2\sigma^2(x) \sum_{p\geq1}\varphi'_p(x)^2 = \|g \mapsto\sqrt{2}
% \mbox{and } &
%)^2 = \|g\mapsto g(qx)+g((1-q)x)-g(x)\|^2_{W^{-7,1}}\leq
%C(1+|x|^2).
Thus, (\ref{etape4547}) gives
%
%e103 ###
%
\begin{eqnarray}\label{etape4546}
&&\sup_{T>T_0}\sup_{\varsigma<\delta}\PP_\mu( |{\operatorname{tr}}_{W^{-7,1}}\dlangle M^T \drangle_{\tau_T+\varsigma}-
\operatorname{tr}_{W^{-7,1}}\dlangle M^T \drangle_{\tau_T}|>\varepsilon
)\nonumber\\[-1pt]
&&\qquad\leq \frac{C}{\varepsilon} \sup_{T>T_0}\E_\mu\biggl[\int_{\tau
_T}^{\tau_T+\delta} \biggl\langle\frac{Z_{s+T}}{\E
[N_{s+T}]},1+|x|^4\biggr\rangle\biggr] \,ds \nonumber\\[-8pt]\\[-8pt]
&&\qquad\leq \frac{C}{\varepsilon} \sup_{T>T_0} \E_\mu\biggl[ \int
_0^\delta\langle Z_{s+\tau_T+T},1+|x|^4\rangle\rme^{-r(s+\tau
_T+T)}\biggr] \,ds\nonumber\\
&&\qquad\leq \frac{C}{\varepsilon} \sup_{T>T_0}\int_0^\delta\E_\mu
\bigl[ \E_{Z_{\tau_T}}\bigl[\langle Z_{s+T},1+|x|^4\rangle\rme
^{-r(s+T)}\bigr] \bigr] \,ds\nonumber
\end{eqnarray}
by using the strong Markov property of $(Z_t, t\geq0)$. Now, using the
branching property gives
%
%e104 ###
%
\begin{eqnarray}\label{etape4548}
&&\E_{Z_{\tau_T}}\bigl[\langle Z_{s+T},1+|x|^4\rangle\rme
^{-r(s+T)}\bigr] \nonumber\\
&&\qquad=  \int_{\R} \E_y\bigl[\langle Z_{s+T},1+|x|^4\rangle\rme
^{-r(s+T)}\bigr] Z_{\tau_T}(dy)\nonumber\\[-8pt]\\[-8pt]
&&\qquad= \int_{\R} \E_y[1+|Y_{s+T}|^4] Z_{\tau_T}(dy)\nonumber
\\
&&\qquad\leq \int_{\R} \bigl(\langle\pi,1+|x|^4\rangle+ B(1+|y|^4)\rme
^{-\beta(s+T)}\bigr) Z_{\tau_T}(dy)
\nonumber
\end{eqnarray}
for some $\beta$ and $B>0$ given by (\ref{exponentielltgeom2}) [see
Remark \ref{rquehyptcl}(ii)]. Since we have a Yule tree, $\E[N_{\tau
_T}]\leq\E[N_\Upsilon]=\exp(r\Upsilon)$. Moreover, using (\ref
{pbmsplitteddiff}) where the integrand in the second term of the right-hand side
is negative for our choice $f(x)=|x|^4$ and noticing that $Z_s$ is a
positive measure, we obtain with localizing arguments that for any
$t\in\R_+$,
%neglection of the negative terms, we obtain for any $t\in\R_+$:
%
\begin{eqnarray*}
&&\E_\mu[\langle Z_{t\wedge\tau_T}, 1+|x|^4\rangle]\\
&&\qquad\leq
\langle\mu,1+|x|^4\rangle+\int_0^t
(8\bar{b}+24\bar{\sigma})\E_\mu[\langle Z_{s\wedge\tau
_T},1+|x|^4\rangle]\,ds.
\end{eqnarray*}
We deduce from Gronwall's lemma that
%
%e105 ###
%
\begin{equation}\label{etape4549}
\E_\mu[\langle Z_{\tau_T}, 1+|x|^4\rangle]\leq\langle
\mu,1+|x|^4\rangle\rme^{(8\bar{b}+24\bar{\sigma})\Upsilon
}.
\end{equation}
Then (\ref{etape4546}), (\ref{etape4548}) and (\ref{etape4549})
imply that
%
%e106 ###
%
\begin{eqnarray}
&&\sup_{T>T_0}\sup_{\varsigma<\delta}\PP_\mu( |
{\operatorname{tr}}_{W^{-7,1}}\dlangle M^T
\drangle_{\tau_T+\varsigma}-\operatorname{tr}_{W^{-7,1}}\dlangle M^T
\drangle_{\tau_T}|>\varepsilon)\nonumber\\[-8pt]\\[-8pt]
&&\qquad\leq\frac{C \delta}{\varepsilon}
\bigl(\rme^{r\Upsilon}+\rme^{(8\bar{b}+24\bar{\sigma})\Upsilon
}\bigr),\nonumber
\end{eqnarray}
which ends the proof of the Aldous inequality for the trace of the martingale.

\begin{rque}\label{rquedanslapreuve}This also shows the tightness of
$(M^T)_{T\geq0}$ in $W^{-7,1}$.
\end{rque}

%***********
%Le sup n'est pas dans l'$\E$. Pour ce qui suit utiliser ``similar
%arguments as in the proof of \reff{estimeeuniformemartingale}''.
%************

%Since the integrand in the right-hand side  of (\ref{aldoustracemart})
%converges in $L^2$ to $\langle\pi,1+|x|^4\rangle$, the right-hand side  of
%The (\ref{aldoustracemart}) can be made as small as possible with a
%proper
%choice of $\delta$.

%************************
%la Cv et la finitude de $(\pi, x^4)$ n'est pas claire
%**********************

For the finite variation part,
%
%e107 ###
%
\begin{eqnarray}\label{aldousvarfinie}
&&
\sup_{T>T_0}\sup_{\varsigma<\delta}\PP\biggl(
\biggl\|\int_0^{\tau_T+\varsigma}(L+J)^*\eta^T_s \,ds-\int_0^{\tau
_T}(L+J)^*\eta^T_s \,ds\biggr\|_{W^{-7,1}}>\varepsilon\biggr)\nonumber\\
&&\qquad
\leq\sup_{T>T_0}\sup_{\varsigma<\delta}\frac{1}{\varepsilon
^2}\E\biggl[\biggl\|\int_{\tau_T}^{\tau_T+\varsigma}(L+J)^*\eta^T_s
\,ds\biggr\|^2_{W^{-7,1}} \biggr]\nonumber\\
&&\qquad\leq
\sup_{T>T_0}\sup_{\varsigma<\delta}\frac{\varsigma}{\varepsilon
^2}\E\biggl[\int_{\tau_T}^{\tau_T+\varsigma}\|(L+J)^*\eta^T_s
\|^2_{W^{-7,1}} \,ds\biggr]\\
&&\qquad\leq\sup_{T>T_0}\frac{C\delta}{\varepsilon^2}\int_0^{\Upsilon
+\delta}\E[\|\eta^T_t \|^2_{C^{-4,2}}
]\,dt\nonumber\\
&&\qquad\leq
\frac{C\delta(\Upsilon+\delta)}{\varepsilon^2}\sup_{T>T_0}
\sup_{t\leq\Upsilon}\E[\|\eta^T_t \|^2_{C^{-4,2}}].\nonumber
\end{eqnarray}
%
%where we use \cite{yosida} p. 133.
%****************
%cette ref est-elle utile?
%|\int f|^2\leq(\int|f|)^2\leq(\int|f|^2) (\int1).
%*****************
%For the second inequality.
We use Cauchy--Schwarz's inequality for the second inequality. For the
third inequality, we notice that under the Assumption \ref{hyptcl} and
for $\varphi\in W^{7,1}$,
%
%e108 ###
%
\begin{equation}
\|L\varphi\|_{C^{4,2}}\leq C\|\varphi\|_{C^{6,1}}\leq C\|\varphi\|_{W^{7,1}}
\end{equation}
as $W^{7,1}\hookrightarrow C^{6,1}$. We can make the right-hand side  of (\ref
{aldousvarfinie}) as small as we wish thanks to (\ref
{estimeeuniforme}) and this ends the proof of the tightness.
\end{pf*}

Then we identify the limit by showing that the limiting values solve an
equation for which uniqueness holds.
This will prove the central limit theorem.

\begin{pf*}{Proof of Proposition \ref{propconvergencefluctuations}}
First of all, by Remark \ref{rquedanslapreuve}, the sequence of
martingales $(M^T)_{T\geq0}$ is tight in $W^{-7,1}$ and thus also in
$C^{-7,0}$ by (\ref{emboitement}). Let us prove that in the latter
space it is moreover $C$-tight in the sense of Jacod and Shiryaev
\cite{jacodshiryaev}, page 315. Using the Proposition 3.26(iii) of this
reference, it remains to prove the convergence of ${\sup_{t\leq
\Upsilon}}\|\Delta M^T_t\|_{C^{-7,0}}$ to $0$ where $\Delta
M^T_t=M^T_t-M^T_{t_-}$. Since the finite variation part of (\ref
{fluctuations}) is continuous, $\Delta M^T_t=\Delta\eta^T_t$
and\vadjust{\goodbreak}
since in (\ref{fluctuationsdef}) $t\mapsto\langle Z_t,Q_tf\rangle$
is continuous, we have for $f\in C^{7,0}$,
%
%e109 ###
%
\begin{eqnarray}\label{etape7}\quad
&&\sup_{t\leq\Upsilon}|\Delta M^T_t(f)|\nonumber\\[-2pt]
&&\qquad= \sup_{t\leq\Upsilon} \rme
^{-{r(t+T)}/{2}}\bigl|f\bigl(q(\omega,t+T)X^{u(\omega,t+T)}_{t+T}
\bigr)\\[-2pt]
&&\qquad\quad\hspace*{65.4pt}{} +
f\bigl(\bigl(1-q(\omega,t+T)\bigr)X^{u(\omega,t+T)}_{t+T}\bigr)-f
\bigl(X^{u(\omega,t+T)}_{t+T}\bigr)\bigr|,\hspace*{-15pt}
\nonumber
\end{eqnarray}
where $u(\omega,t+T)\in V_{t+T}$ is the label of the particle that
undergoes division at $t+T$ and where $q(\omega,t+T)$ is the fraction
which appears in the splitting. By convention, if there is no splitting
at $t+T$, the term in the supremum of the right-hand side  of (\ref{etape7}) is $0$.
Thus, $
\sup_{t\leq\Upsilon}|\Delta M^T_t(f)|\leq3 \rme^{-rT/2}\|f\|
_\infty\leq3 \rme^{-rT/2} \|f\|_{C^{7,0}}$. This proves that
%
%e110 ###
%
\begin{equation}\label{saut->0}
\sup_{t\leq\Upsilon}\|\Delta M^T_t\|_{C^{-7,0}}\leq3\rme
^{-rT/2},
\end{equation}
which converges a.s. to $0$ when $T\rightarrow+\infty$. This finishes
the proof of the $C$-tightness of $(M^T)_{T\geq0}$ in $C^{-7,0}$. The
inequality\vspace*{2pt} (\ref{saut->0}) also ensures that the sequence $\sup
_{t\leq\Upsilon}\|\Delta M^T_t\|_{W^{-7,1}}$ is uniformly integrable.
From the LLN\vspace*{2pt} of Proposition \ref{thLGN}, the integrand of (\ref
{crochetfluctuations}) converges to $W\times V(f)$ which does not
depend on $s$ any more. Thus,\vspace*{2pt} using Theorem 3.12, page 432 in
\cite{jacodshiryaev}, we obtain that $(M^T)_{T\geq0}$ converges in
distribution in $\D([0,\Upsilon],C^{-7,0})$ to a~Gaussian process
$\mathcal{W}$ with the announced quadratic variation. Since $W$ is
$\bigcap_{\varepsilon>0}\sigma(\eta_\varepsilon)$-measurable, it
follows that $W$ and $\mathcal{W}$ are independent.

By Proposition \ref{proptension}, the sequence $(\eta^T, T\geq 0)$
is tight in $W^{-7,1}$ and hence, also in $C^{-7,0}$ by (\ref
{emboitement}). Let $\eta$ be an accumulation point in $\D
([0,\Upsilon], C^{-7,0})$. Because of (\ref{fluctuations}) and (\ref
{saut->0}), $\eta$ is almost surely a continuous process. Let us call
again by $(\eta^T)_{T\geq0}$, with an abuse of notation, the
subsequence that converges in law to $\eta$. Since $\eta$ is
continuous, we get from (\ref{fluctuations}) that it
solves (\ref{limitefluctuations}). Using Gronwall's inequality, we
obtain that
this equation admits in $\Co([0,\Upsilon],C^{-7,0})$ a unique solution
for a given Gaussian white noise $W$ which is in $C^{-7,0}$.
This achieves the proof.
%Our purpose is to prove that $\eta$ is the solution of (
%a unique solution for a given Gaussian white noise $W$. This
%uniqueness can be obtained by using Gronwall's inequality. Let us
%consider the functional $\Psi$ that associates to $\nu\in
%As $(\eta^T)_{T\geq0}$ converges to the continuous $\eta$, we have by
%using Assumption \ref{hyptcl} (i), (\ref{eq:Ubound}) and computation
%similar to those in the proof of Proposition \ref{proptension} (i),
%that $\forall\varphi\in C^{7,0}, \forall t\in[0,\Upsilon],
\vspace*{-4pt}\end{pf*}

%suskaldyti doi

% imsref loaded by lrinkeviciute, 2011-03-14 16:32:06
% imsref loaded by lrinkeviciute, 2011-03-15 08:11:01
% imsref loaded by lrinkeviciute, 2011-03-15 08:36:12
%

%
\printaddresses

\end{document}